\documentclass{article}

\usepackage[nonatbib,final]{neurips_2024}

\usepackage[utf8]{inputenc} 
\usepackage[T1]{fontenc}    
\usepackage[ngerman, english]{babel} 
\usepackage{csquotes} 

\usepackage{enumitem} 
\usepackage[unicode, colorlinks=true, linkcolor=.]{hyperref} 
\usepackage{url}            

\usepackage{xcolor}
\usepackage{nicefrac}       
\usepackage{microtype}      

\usepackage{booktabs}       

\usepackage{wrapfig} 
\usepackage{tikz, tikz-cd, pgfplots} 
\pgfplotsset{compat=1.16}
\usetikzlibrary{arrows,positioning}

\usepackage{mathtools} 
\usepackage{amsthm, amssymb} 
\usepackage{thmtools} 
\usepackage{aligned-overset} 

\usepackage{cancel} 

\usepackage[final]{fixme} 

\usepackage[
	style=numeric,
	backend=biber,
	seconds=true,
	date=iso,
	urldate=iso,
	maxcitenames=2,
	natbib
]{biblatex}
\newcommand*{\grad}{\nabla}

\newcommand*{\real}{\mathbb{R}}
\newcommand*{\nat}{\mathbb{N}}
\newcommand*{\dims}{d}

\newcommand*{\E}{\mathbb{E}}
\renewcommand{\Pr}{\mathbb{P}}
\newcommand*{\iid}{\text{\normalfont iid}}
\newcommand*{\cdf}{\Phi}

\newcommand*{\identity}{\mathbb{I}}
\newcommand*{\ind}{\mathbf{1}}
\DeclareMathOperator{\diag}{diag}

\newcommand*{\Loss}{\mathcal{L}}
\newcommand*{\Cost}{\mathbf{J}}
\newcommand*{\cost}{J}
\newcommand*{\rf}{\mathbf{f}}
\newcommand*{\rg}{\mathbf{g}}
\newcommand*{\loss}{\ell}
\newcommand*{\param}{w}
\newcommand*{\Param}{W}

\DeclareMathOperator*{\argmin}{argmin}
\DeclareMathOperator*{\argmax}{argmax}
\DeclareMathOperator*{\Root}{Root}

\newcommand{\normal}{\mathcal{N}}

\newcommand*{\step}{\mathbf{d}}
\newcommand*{\lr}{h}
\newcommand*{\stepsize}{\eta}

\newcommand*{\scale}{s}
\DeclareMathOperator{\LUE}{\hyperref[def: LUE]{LUE}} 
\DeclareMathOperator{\BLUE}{\hyperref[def: BLUE]{BLUE}} 

\newcommand*{\modifiedBessel}{K_\nu}
\newcommand*{\bigO}{\mathcal{O}}

\newcommand*{\kernel}{\kappa}
\DeclareMathOperator{\linHull}{span}
\DeclareMathOperator{\Span}{span}
\DeclareMathOperator{\Cov}{Cov}
\DeclareMathOperator{\Var}{Var}
\newcommand*{\noise}{\varsigma}
\newcommand*{\C}{\mathcal{C}}
\newcommand*{\sqC}{C} 
\newcommand*{\ikernel}{C} 
\newcommand*{\signal}{\theta}

\newcommand*{\mcoeff}[1][p]{\hyperref[eq: matern abstractions]{a_{#1}}}
\newcommand*{\dcoeff}[1][p]{\hyperref[def: dcoeff]{\tilde{a}_{#1}}}
\newcommand*{\mf}[1][\nu]{\hyperref[eq: matern abstractions]{f_{#1}}}

\newcommand*{\proj}[1]{P_{#1}}

\newcommand*{\batchsize}{b}
\newcommand*{\Batchsize}{B}

\newcommand*{\cX}{\mathcal{X}}

\definecolor{magenta}{HTML}{D81B66}
\definecolor{blue}{HTML}{1F76C1}
\definecolor{yellow}{HTML}{FFC107}
\definecolor{teal}{HTML}{00B981}

\newcommand*{\blue}[1]{{\color{blue} #1}}
\newcommand*{\magenta}[1]{{\color{magenta} #1}}
\newcommand*{\teal}[1]{{\color{teal} #1}}

\theoremstyle{plain}
\newtheorem{prop}{Proposition}[section]
\newtheorem{lemma}[prop]{Lemma}
\newtheorem{corollary}[prop]{Corollary}
\newtheorem{theorem}[prop]{Theorem}

\newtheorem{extension}[prop]{Extension}

\theoremstyle{definition}
\newtheorem{definition}[prop]{Definition}
\newtheorem{example}[prop]{Example}

\theoremstyle{remark}
\newtheorem{remark}[prop]{Remark}

\title{Random Function Descent}

\author{
	Felix Benning\\
	University of Mannheim\\
	\texttt{felix.benning@uni-mannheim.de}
	\And
	Leif Döring\\
	University of Mannheim\\
	\texttt{leif.doering@uni-mannheim.de}
}

\begin{document}

	\maketitle

	\begin{abstract}
Classical worst-case optimization theory neither explains the success of
optimization in machine learning, nor does it help with step size selection.
In this paper we demonstrate the viability and advantages of replacing the
classical `convex function' framework with a `random function' framework.
With complexity \(\bigO(n^3\dims^3)\), where \(n\) is the number of steps and
\(\dims\) the number of dimensions, Bayesian optimization with gradients has
not been viable in large dimension so far. By bridging the gap between
Bayesian optimization (i.e. random function optimization theory) and
classical optimization we establish viability. Specifically, we use a
`stochastic Taylor approximation' to rediscover gradient descent, which is
scalable in high dimension due to \(\bigO(n\dims)\) complexity. This
rediscovery yields a specific step size schedule we call Random Function Descent
(RFD). The advantage of this random function framework is that RFD is scale
invariant and that it provides a theoretical foundation for common step size
heuristics such as gradient clipping and gradual learning rate warmup.
\end{abstract}
	\section{Introduction}

Cost function minimization is one of the most fundamental mathematical problems in
machine learning. Gradient-based methods, popular for this task, require a step
size, typically chosen using established heuristics. This article aims to deepen
the theoretical understanding of these heuristics and proposes a new algorithm
based on this insight. 

Classical optimization theory uses \(L\)-smoothness,
which limits the rate of change of the gradient by \(L\), to provide some
convergence guarantees for learning rates smaller than \(1/L\)
\parencite[e.g.][]{nesterovLecturesConvexOptimization2018}. As this
theory is based on an upper bound (the worst case), the learning rate \(1/L\) is
naturally much more conservative than necessary on average. Even if \(L\) was
known, this learning rate would therefore be impractical. Since line search algorithms typically
require access to full cost function evaluations,
the field of machine learning (ML) therefore
relies heavily on step size heuristics
\parencite[e.g.][]{smithDisciplinedApproachNeural2018,smithCyclicalLearningRates2017,pascanuDifficultyTrainingRecurrent2013,goyalAccurateLargeMinibatch2018}.
To investigate these heuristics, we introduce new ideas based on a `random function'
perspective.

While automatic step size selection in the convex function framework
is possible \parencite{defazioLearningRateFreeLearningDAdaptation2023},
convexity is generally only satisfied asymptotically and locally. So the
understanding of the initial stages of optimization, which includes the warmup
heuristic \parencite{goyalAccurateLargeMinibatch2018}, greatly benefits from a
framework which also admits non-convex functions. This objective is achieved by the
`random function' framework we investigate. 

Many successful algorithms in computer science are significantly slower in the
worst case than in the average case based on a probabilistic framework (e.g. Quicksort
\parencite{hoareQuicksort1962} or the simplex algorithm
\parencite[e.g.][]{borgwardtSimplexMethodProbabilistic1986}). On random quadratic
functions the average case behavior of first order optimizers is already being
investigated by the ML community
\citep[e.g.][]{zhangWhichAlgorithmicChoices2019,
	pedregosaAccelerationSpectralDensity2020,
	lacotteOptimalRandomizedFirstOrder2020,
	deiftConjugateGradientAlgorithm2021,
	cunhaOnlyTailsMatter2022,
	paquetteHaltingTimePredictable2022,
	paquetteUniversalityConjugateGradient2022}.
Interested in the landscape of high dimensional random functions
as a model for `spin glasses`, the physics community independently started
studying the average case of optimization as well \parencite[e.g.][]{
	auffingerComplexityGaussianRandom2023,
	elalaouiOptimizationMeanfieldSpin2021,
	montanariOptimizationSherringtonKirkpatrick2021,
	subagFollowingGroundStates2021,huangTightLipschitzHardness2022},
albeit not geared for ML algorithms.

Average case analysis fundamentally requires a prior distribution over possible cost functions.
The evaluations seen so far then result in a posterior over the cost of other parameter inputs.
Using this posterior for optimization is called ``Bayesian optimization'' (BO)
\citep[e.g.][]{kushnerNewMethodLocating1964,
	shahriariTakingHumanOut2016,
	frazierBayesianOptimization2018,
	agnihotriExploringBayesianOptimization2020},
which is best known in the context of low dimensional optimization (e.g.
hyperparameter tuning) in the ML community. BO is treated like a
zero order method for low dimensional problems due to the \(\bigO(n^3)\)
complexity for the covariance matrix inversion of the \(n\) evaluations seen so
far, which increases to \(\bigO(n^3\dims^3)\) when gradient information is included
\citep[e.g.][]{lizottePracticalBayesianOptimization2008,wuBayesianOptimizationGradients2017},
where \(\dims\) is the input dimension of
our cost function. This limits classic BO to relatively small dimensions even under
sparsity considerations
\citep[e.g.][]{roosHighDimensionalGaussianProcess2021,padidarScalingGaussianProcesses2021}.

While the BO algorithms developed in the `random function framework' might not
have been viable in high dimension so far, due to their computational
complexity, this framework is already used to explain the high relative
frequency of saddle points in high dimension
\citep{dauphinIdentifyingAttackingSaddle2014} and to explain the highly
predictable progress optimizers make on high dimensional cost
functions \citep{benningGradientSpanAlgorithms2024}.

In this work we bridge the gap between BO and (computationally viable) gradient
based methods, derived from the first Taylor approximation, with the
introduction of a stochastic Taylor approximation based on a forgetful BO
posterior.  The optimization method ``Random Function Descent'' (RFD), resulting
from the minimization of this stochastic Taylor approximation, coincides with a
specific form of gradient descent which establishes its viability in high
dimension. The advantages of its BO heritage are scale invariance and an
explicit step size schedule, which illuminates the inner workings of step size
heuristics such as gradient clipping
\parencite{pascanuDifficultyTrainingRecurrent2013} and gradual learning rate
warmup \parencite{goyalAccurateLargeMinibatch2018}.

\paragraph*{Our contributions and outline}
The main goal of this paper is to demonstrate the \textbf{viability} and
\textbf{advantages} of replacing the classical ``convex function'' framework with
a ``random function'' framework. Theorem~\ref{thm: explicit rfd} is the main
theoretical result establishing \textbf{viability} (computatability and scalable
complexity) for a given covariance model. Section~\ref{sec: stochastic loss and
covariance estimation} is concerned with practical estimation of the covariance
model and viability is demonstrated with a practical example in the MNIST case
study (Section~\ref{sec: mnist case study}). The \textbf{advantages} of this
approach are scale invariance (Advantage~\ref{advant: scale invariance}) and an
explicit step size schedule, which does not require expensive tuning and
explains existing ML heuristics such as warmup (cf. Section~\ref{subsec: step
size heuristics}). This explanation of the initial stage of optimization could
never be delivered by the convex framework, because the convexity assumption is
not fulfilled initially so it can at best explain asymptotic behavior.

\begin{description}
	\item[Sec.~\ref{sec: rfd}] We motivate a stochastic Taylor approximation and RFD
	and prove its scale-invariance.

	\item[Sec.~\ref{sec: distribution over functions}] We briefly motivate
	and discuss the common distributional assumptions in BO.
	
	\item[Sec.~\ref{subsec: explicit rfd}] We establish the connection between
	RFD and gradient descent.

	\item[Sec.~\ref{subsec: rfd step sizes}] We investigate the step size schedule
	suggested by RFD. In particular we
	\begin{enumerate}[start=0,leftmargin=0cm]
		\item calculate explicit formulas for the step size schedules resulting
		from common covariance models (Table~\ref{table: optimal step size}, Sec.~\ref{appendix: covariance models}),
		\item analyze the general asymptotic behavior (Sec.~\ref{subsec: asymptotic RFD step sizes}),
		\item discuss how RFD explains gradient clipping and learning rate warmup (Sec.~\ref{subsec: step size heuristics}),
	\end{enumerate}

	\item[Sec.~\ref{sec: stochastic loss and covariance estimation}]
	We develop a non-parametric variance estimation method, which is robust with
	respect to the choice of covariance kernel. Finally, we present an extension of RFD to
	mini-batch losses.
	
	\item[Sec.~\ref{sec: mnist case study}] We conduct a case study on the MNIST dataset.
	\item[Sec.~\ref{sec: extensions and limitations}] We discuss extensions (see also Sec.~\ref{sec: extensions}) and limitations.
\end{description}

	\section{The random function descent algorithm}\label{sec: rfd}

The classic derivation of gradient descent \citep[e.g.][p.
29]{nesterovLecturesConvexOptimization2018},
adds an \(L\)-smoothness based trust bound to the first Taylor approximation,
\(T[\cost(\theta) \mid \cost(\param),\nabla \cost(\param)]\), of the cost
function \(\cost\) around \(\param\) resulting in the gradient step
\[
	\param - \tfrac1L \nabla \cost(\param)
	= \argmin_\theta T[ \cost(\theta)\mid \cost(\param),\nabla \cost(\param)]
	+ \tfrac{L}2 \|\theta - \param\|^2.
\]
Our unusual notation for the Taylor approximation \(T[ \cost(\theta)\mid
\cost(\param),\nabla \cost(\param)]\) is meant to highlight the connection to
the stochastic Taylor approximation we define below.

\begin{wrapfigure}{R}{0.5\linewidth}
	\centering
	\def\svgwidth{\linewidth}
	\vspace*{-1.2cm}
	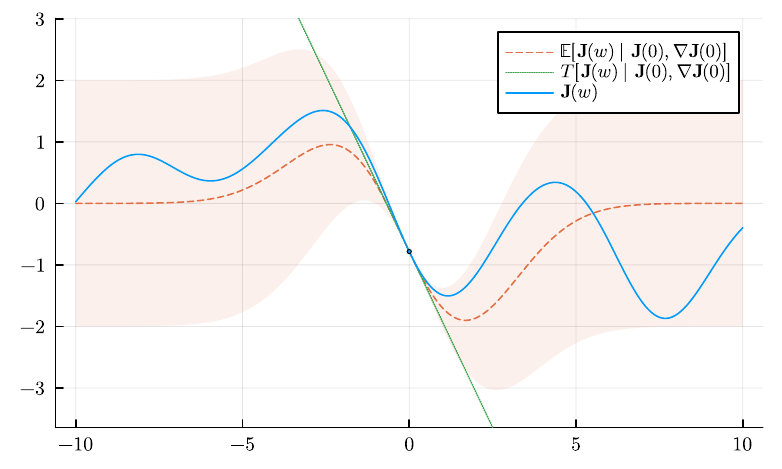
	\caption{
		The stochastic Taylor approximation naturally contains a trust bound
		in contrast to the classical one. Here \(\Cost\) is a Gaussian
		random function (with covariance as in Equation~\eqref{eq: sqExp
		covariance model}, with length scale \(\scale=2\) and variance \(\sigma^2=1\)).
		The ribbon represents two conditional standard
		deviations around the conditional expectation.
	}
	\label{fig: visualize conditional expectation}
\end{wrapfigure}

\begin{definition}[Stochastic Taylor approximation]
	We define the first order stochastic Taylor approximation of a random
	(cost) function\footnote{
		\textbf{Remark on terminology}: ``stochastic process''
		\parencite[e.g.][]{williamsBayesianClassificationGaussian1998}, ``random
		field'' \parencite[e.g.][]{adlerRandomFieldsGeometry2007} and
		``random function''
		\parencite[e.g.][]{matheronIntrinsicRandomFunctions1973} are all synonyms.
		However the latter seems most descriptive of random variables in the set
		of functions. ``Gaussian processes'' are naturally Gaussian stochastic processes,
		i.e. Gaussian random
		functions. To better distinguish random functions from deterministic
		functions, we use bold letters to denote random functions (as the usual
		convention of capitalizing random variables often clashes with other
		conventions for functions).
	} \(\Cost\) around \(\param\) by the
	conditional expectation
	\[
		\E[\Cost(\theta) \mid \Cost(\param), \nabla\Cost(\param)].
	\]
	This is the best \(L^2\) approximation
	\parencite[Cor.~8.17]{klenkeProbabilityTheoryComprehensive2014} of
	\(\Cost(\theta)\) provided first order knowledge of \(\Cost\) at
	\(\param\).
\end{definition}

We call this the `stochastic Taylor approximation' because this approximation
only makes use of derivatives in a single point. While the standard Taylor
approximation is a polynomial approximation, the `stochastic Taylor
approximation' is the best approximation in an \(L^2\) sense and already
mean-reverting by itself, i.e. it naturally incorporates covariance-based trust
(cf.~Figure~\ref{fig: visualize conditional expectation}).
While \(L\)-smoothness-based trust \emph{guarantees} that the gradient still
points in the direction we are going (for learning rates smaller \(1/L\)), covariance
based trust tells us whether the derivative is still negative \emph{on average}.
Minimizing the stochastic Taylor approximation is therefore optimized for the
average case. Since convergence proofs for gradient descent typically rely on an
improvement \emph{guarantee}, proving convergence is significantly harder in the average
case and we answer this question only partially in
Corollary~\ref{prop: convergence}.

\begin{definition}[Random Function Descent -- RFD]
	\label{def: rfd}
	Select \(\param_{n+1}\) as the minimizer\footnote{
		we ignore throughout the main body that \(\argmin\) could be set-valued and
		that the \(\param_n\) would be random variables (cf.~Section~\ref{sec: formal rfd}
		for a formal approach).
	}
	of the first order stochastic Taylor
	approximation
	\begin{equation*}
		\param_{n+1}
		:= \argmin_{\param}\E[\Cost(\param)\mid \Cost(\param_n), \nabla\Cost(\param_n)].
	\end{equation*}
\end{definition}

\paragraph*{Properties of RFD}
Before we make RFD more explicit in Section~\ref{subsec: explicit rfd}, we
discuss some properties which are easier to see in the abstract form.

First, observe that RFD is greedy and forgetful in the same way gradient descent is greedy
and forgetful when derived as the minimizer of the regularized first Taylor
approximation, or the Newton method as the minimizer of the second Taylor
approximation. This is because the Taylor approximation only uses derivatives
from the last point \(\param_n\) (forgetful), and we minimize this approximation
(greedy). Since momentum methods retain some information about past gradients,
they are not as forgetful. We therefore expect a similar improvement
could be made for RFD in the future.

Second, it is well known that classical gradient descent with exogenous step sizes (and
most other first order methods) lack the scale invariance property of the Newton
method \parencite[e.g.][]{hanssonOptimizationLearningControl2023,deuflhardAffineInvariantConvergence1979}. Scale invariance means that scaling the input parameters \(\param\) or the
cost itself (e.g. by switching from the mean squared error to the sum squared
error) does not change the points selected by the optimization method.

\begin{restatable}[Scale invariance]{advantage}{scaleInvariance}
	\label{advant: scale invariance}
	RFD is invariant to additive shifts and positive scaling of the cost
	\(\Cost\). RFD is also invariant with respect to transformations of the
	parameter input of \(\Cost\) by differentiable bijections whose Jacobian is
	invertible everywhere (e.g. invertible linear maps).
\end{restatable}
While invariance to bijections of inputs is much stronger than the affine
invariance offered by the Newton method, non-linear bijections will typically
break the `isotropy' assumption of the following section which makes RFD explicit.
This invariance should therefore be viewed as an opportunity to look for the
bijection of inputs which ensures isotropy (e.g. a whitening
transformation). The discussion of geometric anisotropy in Section~\ref{sec:
geometric anisotropy} is conducive to build an understanding of this.

	\section{A distribution over cost functions}
\label{sec: distribution over functions}


It is impossible to make average case analysis explicit without a distribution over
functions, so we use the canonical distributional
assumption of Bayesian optimization
\parencite[e.g.][]{frazierBayesianOptimization2018,wuBayesianOptimizationGradients2017,rasmussenGaussianProcessesMachine2006},
`isotropic Gaussian random functions'. This assumption was also used in the high
dimensional setting by \citet{dauphinIdentifyingAttackingSaddle2014} to argue
that saddle points are much more common than minima in high dimension, which is
often cited to explain why second order methods are uncommon in machine
learning.

To motivate isotropy, we note that in average case analysis the uniform
distribution is popular, since it weighs all problem instances equally (e.g.
all possible permutations in sorting). Isotropy is such a uniformity assumption,
which essentially requires ``\(\Pr(\Cost=\cost) = \Pr(\Cost=\cost\circ
\phi)\)``, for all isometries \(\phi\). In other words, the probability that our
cost function is equal to \(\cost\) is equal to the probability that it is
equal to a shifted and turned version of \(\cost\), given by \(\cost\circ\phi\).

Since the probability of any single realization of a cost function \(\cost\)
is zero, the equation we put in quotes is mathematically unsound. The formal
definition follows below.


\begin{definition}[Isotropy]
    A random function \(\Cost\) is called isotropic if its distribution stays
    the same under isometric transformations of its input, i.e.
    for any isometry \(\phi\) we have
    \[
        \Pr_{\Cost} = \Pr_{\Cost \circ \phi}.
    \]

    If \(\Cost\) is Gaussian, isotropy is well known
    \citep[e.g.][]{rasmussenGaussianProcessesMachine2006,adlerRandomFieldsGeometry2007}
    to be equivalent to the condition that there exists \(\mu\in \real\) and a
    function \(\ikernel:\real\to\real\)
    such that for all \(\param,\tilde{\param}\in \real^\dims\) the expectation and covariance
    are
    \[
        \E[\Cost(\param)] = \mu,
        \qquad
        \Cov(\Cost(\param), \Cost(\tilde{\param})) = \ikernel\bigl(\tfrac{\|\param-\tilde{\param}\|^2}{2}\bigr).
    \]
    For these isotropic Gaussian random functions we use the notation \(\Cost\sim\normal(\mu, \ikernel)\).
\end{definition}

We discuss generalizations to isotropy in Section~\ref{sec: appendix: input
invariance} and \ref{sec: geometric anisotropy}, but for ease of exposition we
retain the (stationary) isotropy assumption throughout the main body.
Note that the Gaussian assumption can be statistically tested in practice
(cf.~Figure~\ref{fig:
covariance estimation}), but it is also straightforward to reproduce our results
with the ``best linear unbiased estimator'' (BLUE) (Section~\ref{sec: BlUE})
 in place of
the conditional expectation to remove the Gaussian assumption.
We finally want to highlight that, in contrast
to the uniformity assumption on finite sets, `isotropic Gaussian random
functions' leave us with a family of plausible distributions. It is therefore
necessary to estimate \(\mu\) and \(\ikernel\), which is the topic of
Section~\ref{sec: stochastic loss and covariance estimation}.

	\section{Relation to gradient descent}
\label{subsec: explicit rfd}

While we were able to define RFD abstractly without any assumptions on the
distribution \(\Pr_{\Cost}\) of the random cost \(\Cost\), an explicit
calculation requires distributional assumptions and we have motivated isotropic
Gaussian random functions in Section~\ref{sec: distribution over functions} for
this purpose.
The assumption of isotropy allows for an explicit version of the stochastic
Taylor approximation which then immediately leads to an explicit version of RFD.

\begin{restatable}[Explicit first order stochastic Taylor approximation]{lemma}{firstStochTaylor}
	\label{lem: first stoch Taylor}
	For \(\Cost\sim\normal(\mu, \ikernel)\), the first order stochastic Taylor
	approximation is given by
	\[
		\E[\Cost(\param-\step)\mid \Cost(\param),\nabla\Cost(\param)]
		= \mu + \frac{\ikernel\bigl(\frac{\|\step\|^2}2\bigr)}{\ikernel(0)}
		(\Cost(\param)-\mu) - \frac{\ikernel'\bigl(\frac{\|\step\|^2}2\bigr)}{\ikernel'(0)}
		\langle \step, \nabla\Cost(\param)\rangle.
	\]
\end{restatable}

The explicit version of RFD follows by fixing the step size
\(\stepsize = \|\step\|\) and optimizing over the direction first.

\begin{restatable}[Explicit RFD]{theorem}{explicitRFD}
	\label{thm: explicit rfd}
	Let \(\Cost\sim\normal(\mu, \ikernel)\), then RFD coincides with gradient
	descent
	\[
		\param_{n+1}
		= \param_n  -\stepsize_n^*\tfrac{\nabla\Cost(\param_n)}{\|\nabla\Cost(\param_n)\|},
	\]
	where the RFD step sizes are given by
	\begin{equation}
		\label{eq: rfd step size}
		\stepsize_n^*
		:= \argmin_{\stepsize\in\real} \frac{\ikernel\bigl(\frac{\stepsize^2}2\bigr)}{\ikernel(0)}(\Cost(\param_n)-\mu)
		- \stepsize \frac{\ikernel'\bigl(\frac{\stepsize^2}2\bigr)}{\ikernel'(0)}\|\nabla\Cost(\param_n)\|.
	\end{equation}
\end{restatable}

\begin{table*}[t]
	\caption{
		RFD step size (cf.~Figure~\ref{fig: rfd step sizes} and
		Eq.~\eqref{eq: sqExp covariance model}, \eqref{eq: rational
		quadratic}, \eqref{eq: matern model} for the formal definitions of the models).
		In particular, \(\scale\) is the length scale in all covariance models.
	} \label{table: optimal step size}
	\centering
	\begin{tabular}{p{1.7cm} llll l}
	  	\toprule

	  	\multicolumn{2}{c}{Model} & \multicolumn{2}{c}{RFD step size \(\stepsize^*\) for \(\Cost(\param) \le \mu\)} 
		& A-RFD
		\\
		\cmidrule(r){3-5}
		& & General case \(\bigl(\text{with }\Theta=\tfrac{\|\nabla\Cost(\param)\|}{\mu-\Cost(\param)}\bigr)\)
		& \(\Cost(\param) = \mu\)
		& \(\Theta\to 0\)
		\\
	  	\midrule
		Matérn & \(\nu\) \\
	  	\cmidrule(r){1-2}
		& \(3/2\) 
		& \(
			\frac{\scale}{\sqrt{3}}\frac{1}{\left(
				1 + \frac{\sqrt{3}}{s\Theta}
			\right)}
		\)
		& \(\approx 0.58\scale\) 
		& \(\frac{1}{3}\scale^2 \Theta\)
		\\
		& \(5/2\) 
		& \(
			\frac{\scale}{\sqrt{5}}\frac{
				(1-\zeta)+\sqrt{4 + (1+\zeta)^2}
			}{2(1+\zeta)}\)
			with \(\zeta := \frac{\sqrt{5}}{3\scale\Theta}.
		\)
		& \(\approx 0.72\scale\)
		& \(\frac{3}{5} \scale^2\Theta\)
		\\
		Squared-exponential
		& \(\infty\)
		& \(
			\frac{
				\scale^2
			}{
				\sqrt{\bigl(\frac{\mu-\Cost(\param)}{2}\bigr)^2+\scale^2\|\nabla\Cost(\param)\|^2}
				+ \frac{\mu-\Cost(\param)}{2}
			}\|\nabla\Cost(\param)\|
		\)
		& \(\scale\)
		& \(\scale^2\Theta\)
		\\
	  	\cmidrule(r){1-2}
		Rational quadratic & \(\beta\)
		& \(
		\scale \sqrt{\beta} \Root\limits_\stepsize\left(
			- 1 + \frac{\sqrt{\beta}}{\scale\Theta}\stepsize
			+ (1+\beta)\stepsize^2 + \frac{\sqrt{\beta}}{\scale\Theta}\stepsize^3
		\right)
		\) &
		\(
			\scale \sqrt{\frac{\beta}{1+\beta}}
		\)
		& \(\scale^2\Theta\)
		\\
		\bottomrule
	\end{tabular}
\end{table*}

While the descent direction is a universal property for all isotropic
Gaussian random functions, it follows from \eqref{eq: rfd step size} that the
step sizes depend much more on the specific covariance structure. In particular
it depends on the decay rate of the covariance acting as the trust bound.

\begin{remark}[Scalable complexity]
	While Bayesian optimization typically has computational complexity \(\bigO(n^3\dims^3)\)
	in number of steps \(n\) and dimensions \(\dims\)
	\citep{wuBayesianOptimizationGradients2017,roosHighDimensionalGaussianProcess2021},
	RFD under the isotropy assumption has the same computational complexity as
	gradient descent (i.e. \(\bigO(n\dims)\)).
\end{remark}

\begin{remark}[Step until the given information is no longer informative]
While \(L\)-smoothness-based trust prescribes step sizes that \emph{guarantee}
the slope to point downwards over the entire step, RFD prescribes steps which
are exactly large enough that the gradient is no longer correlated to the
previously observed evaluation. This is because the first order condition
demands
\[
	0 \overset{!}= \nabla \E[\Cost(\param)\mid \Cost(\param_n), \nabla\Cost(\param_n)]
	= \E[\nabla \Cost(\param)\mid \Cost(\param_n), \nabla\Cost(\param_n)].
\]
And for measurable functions \(\phi:\real^{\dims+1}\to\real\) such that \(\Phi=\phi(\Cost(\param_n),
\nabla\Cost(\param_n))\) is sufficiently integrable, \(\Phi\) is then uncorrelated
from \(\partial_{i}\Cost(\param)\) by the first order condition
\[
	\Cov(\partial_{i}\Cost(\param), \Phi)
	= \E\Bigl[
		\underbrace{\E[\partial_{i}\Cost(\param) \mid \Cost(\param_n), \nabla\Cost(\param_n)]}_{=0}
		(\Phi-\E[\Phi]) \Bigr] = 0.
\]
\end{remark}

\section{The RFD step size schedule}
\label{subsec: rfd step sizes}

While classical theory leads to
`learning rates', RFD suggests `step sizes' applied to normalized gradients
representing the actual length of the step size. In the following we thus make the distinction
\[
	\param_{n+1}
	= \param_n - \underbrace{\lr_n}_{\mathclap{\text{`learning rate'}}} \nabla\Cost(\param_n)
	= \param_n - \underbrace{\stepsize_n}_{\text{\text{`step size'}}}
	\tfrac{\nabla\Cost(\param_n)}{\|\nabla\Cost(\param_n)\|}.
\]
To get a better feel for the step sizes suggested by RFD, it is enlightening to
divide \eqref{eq: rfd step size} by \(\mu-\Cost(\param_n)\) which results in a
minimization problem
\begin{equation}
	\label{eq: simplified step size opt}
	\stepsize^* : =\stepsize^*(\Theta) := \argmin_\stepsize q_\Theta(\stepsize)
	\qquad\text{for} \qquad
	q_\Theta(\stepsize) := -\frac{\ikernel\bigl(\frac{\stepsize^2}2\bigr)}{\ikernel(0)}
		- \stepsize \frac{\ikernel'\bigl(\frac{\stepsize^2}2\bigr)}{\ikernel'(0)}\Theta,
\end{equation}
which is only parametrized by the ``gradient cost quotient''
\[
	\Theta_n = \frac{\|\nabla\Cost(\param_n)\|}{\mu-\Cost(\param_n)},
\]
\begin{wrapfigure}{R}{0.5\linewidth}
	\def\svgwidth{\linewidth}
	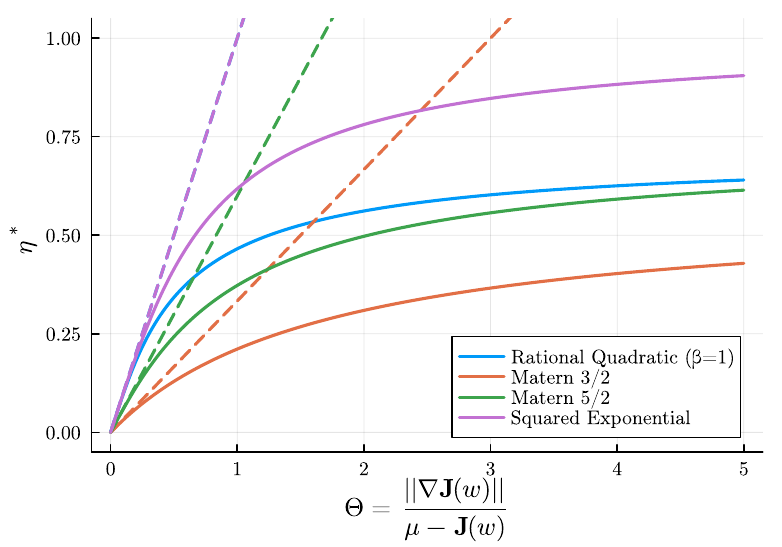	
	\caption{
		RFD step sizes as a function of \(\Theta=\frac{\|\nabla\Cost(\param)\|}{\mu -
		\Cost(\param)}\) assuming scale \(\scale=1\) (cf.~Table~\ref{table: optimal step
		size}). A-RFD (Definition~\ref{def: a-rfd}) is plotted as dashed lines.
		A-RFD of the rational quadratic coincides with A-RFD of the squared
		exponential covariance.
	}
	\label{fig: rfd step sizes}
\end{wrapfigure}
i.e. \(\stepsize_n^* = \stepsize^*(\Theta_n)\).
This minimization problem can be solved explicitly for the most common
\citep[ch.~4]{rasmussenGaussianProcessesMachine2006} differentiable isotropic
covariance models, see Table~\ref{table:
optimal step size}, Figure~\ref{fig: rfd step sizes} and Appendix~\ref{appendix:
covariance models} for details.

Figure~\ref{fig: rfd step sizes} can be interpreted as follows: At the start of
optimization, the cost should be roughly equal to the average cost \(\mu \approx
\Cost(\param)\), so the gradient cost quotient \(\Theta\) is infinite and the
step sizes are therefore given by \(\stepsize^*(\infty)\) (also listed
in its own column in Table~\ref{table: optimal step size}). As we start
minimizing, the difference \(\mu-\Cost(\param)\) becomes positive. Towards the
end of minimization this difference no longer changes as the cost no longer
decreases. I.e. towards the end the gradient cost quotient \(\Theta\) is roughly
linear in the gradient \(\|\nabla\Cost(\param)\|\). The derivative
\(\frac{d}{d\Theta}\stepsize^*(0)\) of \(\stepsize^*(\Theta)\) at zero then effectively
results in a constant asymptotic learning rate.

\subsection{Asymptotic learning rate}
\label{subsec: asymptotic RFD step sizes}

To explain the claim above, note that the gradient cost quotient \(\Theta\)
converges to zero towards the end of optimization, because the gradient norm
converges to zero. A first order Taylor expansion of \(\stepsize^*\) would
therefore imply
\[
	\stepsize^*(\Theta)
	\approx \stepsize^*(0) + \tfrac{d}{d\Theta}\stepsize^*(0) \Theta
	= \underbrace{\tfrac{\tfrac{d}{d\Theta}\stepsize^*(0)}{\mu-\Cost(\param)}}_{\text{asymptotic learning rate}} \|\nabla\Cost(\param)\|
\]
assuming \(\stepsize^*(0) = 0\) and differentiability of \(\stepsize^*\), which
is a reasonable educated guess based on the the examples in Figure~\ref{fig: rfd
step sizes}.
But since the RFD step sizes \(\stepsize^*\) are abstractly defined as an
\(\argmin\), it is necessary to formalize this intuition for general covariance
models. First, we define asymptotic step sizes as an object towards which we can
prove convergence. Then we prove convergence, proving they are well defined. In
addition, we obtain a more explicit formula for the asymptotic learning rate.

\begin{definition}[A-RFD]
	\label{def: a-rfd}
	We define the step sizes of ``asymptotic RFD'' (A-RFD) to be the minimizer of the
	second order Taylor approximation \(T_2q_\Theta\) of \(q_\Theta\)
	around zero
	\[
		\hat{\stepsize}(\Theta)
		:= \argmin_\stepsize T_2q_\Theta(\stepsize)
		= \tfrac{\ikernel(0)}{-\ikernel'(0)}\Theta
		= \underbrace{\tfrac{\ikernel(0)}{\ikernel'(0)(\Cost(\param)-\mu)}}_{\text{asymptotic learning rate}}\|\nabla\Cost(\param)\|.
	\]
\end{definition}
In the following we prove that these are truly asymptotically equal to the step
sizes \(\stepsize^*\) of RFD.
\begin{prop}[A-RFD is well defined]
	\label{prop: a-rfd well defined}
	Let \(\Cost\sim\normal(\mu, \ikernel)\) and assume there exists \(\stepsize_0>0\) such
	that the correlation for larger distances \(\stepsize\ge \stepsize_0\) are
	bounded smaller than \(1\), i.e. \(\frac{\ikernel(\stepsize^2/2)}{\ikernel(0)} < \rho \in (0,1)\).
	Then the step sizes of RFD are asymptotically equal to the step sizes of
	A-RFD, i.e.
	\[
		\hat{\stepsize}(\Theta)\sim \stepsize^*(\Theta)
		\quad \text{as}\quad \Theta\to 0.
	\]
\end{prop}

Note that the assumption is essentially always satisfied, since the Cauchy-Schwarz inequality implies
\[
	\ikernel\bigl(\tfrac{\|\param-\tilde{\param}\|^2}2\bigr) = \Cov(\Cost(\param),\Cost(\tilde{\param}))
	\le \sqrt{\Var(\Cost(\param))\Var(\Cost(\tilde{\param}))} = \ikernel(0),
\]
where equality requires the random variables to be almost surely equal
\parencite{klenkeProbabilityTheoryComprehensive2014}. If the random function is
not periodic or constant, this will generally be strict.
In the proof, this requirement is only needed to ensure that \(\stepsize^*\) is not
very large. The smallest local minimum of \(q_\Theta\) is always close to
\(\hat{\stepsize}\) even without this assumption (which ensures it is a global minimum).

Figure~\ref{fig: rfd step sizes} illustrates that \(\stepsize^*\to 0\) should
imply \(\Theta\to 0\), resulting in a weak convergence guarantee.

\begin{restatable}{corollary}{convergence}\label{prop: convergence}
	Assume \(\stepsize^*\to 0\) implies \(\Theta\to 0\), the cost \(\Cost\) is
	bounded, has continuous gradients and RFD converges to some point
	\(\param_\infty\). Then \(\param_\infty\) is a critical point and
	the RFD step sizes \(\stepsize^*\) are asymptotically equal to
	\(\hat{\stepsize}\).
\end{restatable}

For the squared exponential covariance model we formally prove that \(\stepsize^*\)
is strictly monotonously increasing in \(\Theta\) and thus \(\stepsize^*\to 0\)
implies \(\Theta\to 0\) (Prop.~\ref{prop: sq exp is
strictly monotonous in xi}). The `bounded' and `continuous gradients'
assumptions are almost surely satisfied for all sufficiently smooth covariance
functions \parencite[cf.][]{adlerRandomFieldsGeometry2007}, where three times
differentiable is more than enough smoothness.

\subsection{RFD step sizes explain common step size heuristics}
\label{subsec: step size heuristics}

Asymptotically, RFD suggests constant learning rates, similar to the classical
\(L\)-smooth setting. We thus define these asymptotic learning rates (as the limit
of the learning rates \(\lr_n\) of iteration \(n\)) to be
\begin{equation}
	\label{eq: asymptotic learning rate}
	\lr_\infty := \frac{\ikernel(0)}{\ikernel'(0)(\Cost(\param_\infty) - \mu)},
\end{equation}
where \(\Cost(\param_\infty)\) is the cost we reach in the limit. If we used
these asymptotic learning rates from the start, step sizes would become too
large for large gradients, as RFD step sizes exhibit a plateau
(cf.~Figure~\ref{fig: rfd step sizes}). To emulate the behavior of RFD with a
piecewise linear function, we could introduce a cutoff whenever our step size
exceeds the initial step size
\(\stepsize^*(\infty)\), i.e.
\[
	\tag{gradient clipping}
	\param_{n+1}
	= \param_n - \min\Bigl\{
		\lr_\infty, \frac{\stepsize^*(\infty)}{\|\nabla\Cost(\param_n)\|}
	\Bigr\}\nabla\Cost(\param_n).
\]
At this point we have rediscovered `gradient clipping'
\parencite{pascanuDifficultyTrainingRecurrent2013}. Since the rational quadratic
covariance has the same asymptotic learning rate \(\lr_\infty\) for every
\(\beta\), its parameter \(\beta\) controls the step size bound
\(\stepsize^*(\infty)\) of gradient clipping (cf.~Table~\ref{table: optimal step
size}, Figure~\ref{fig: rfd step sizes}).

\citet{pascanuDifficultyTrainingRecurrent2013} motivated gradient clipping with the
geometric interpretation of movement towards a `wall' placed behind the minimum.
This suggests that clipping should happen towards the end of training. This
stands in contrast to a more recent
step size heuristic, ``(linear) warmup''
\parencite{goyalAccurateLargeMinibatch2018}, which suggests smaller learning
rates at the start (i.e.
\(\lr_0 = \frac{\stepsize^*(\infty)}{\|\nabla\Cost(\param_0)\|}\)) and gradual ramp-up
to the asymptotic learning rate \(\lr_\infty\). In other words, gradients are not clipped
due to some wall next to the minimum, but because the step sizes would be too
large at the start otherwise. \citet{goyalAccurateLargeMinibatch2018} further observe that
`constant warmup' (i.e. a step change of learning rates akin to gradient clipping) performs
worse than gradual warmup. Since RFD step sizes suggest this gradual increase,
we argue that they may have discovered RFD step sizes empirically (also cf.~Figure~\ref{fig:
summary figure}).
	\section{Mini-batch loss and covariance estimation}
\label{sec: stochastic loss and covariance estimation}

Since we do not have access to evaluations of the cost \(\Cost\) in practice,
we need to prove some results about stochastic losses \(\loss_i\)
before we can apply RFD in practice. For this, assume that we have
independent identically distributed (iid) data \(X_i\) independent of the true
relationship \(\rf\) drawn from \(\Pr_\rf\) resulting in labels
\(Y_i=\rf(X_i)+\noise_i\), where we have added independent iid noise
\(\noise_i\), resulting in loss and cost
\[
    \loss_i(\param) := \loss\bigl(\param, (X_i, Y_i)\bigr)
    \quad\text{and}\quad \Cost(\param) := \E[\loss_i(\param)\mid \rf].
\]
In this setting we confirm (cf.~Lemma~\ref{lem: stoch approx
errors cond. independent}), that the stochastic approximation errors
\[
    \epsilon_i(\param) := \loss_i(\param) - \Cost(\param)
\]
are independent conditional on the true relationship \(\rf\). In particular they
(and all their derivatives) are uncorrelated and also uncorrelated from
\(\Cost\). It follows that mini-batch
losses
\begin{equation}
    \label{eq: mini batch loss}
    \Loss_b(\param)
    := \frac1\batchsize\sum_{i=1}^\batchsize \loss_i(\param)   
    = \Cost(\param) + \frac1\batchsize\sum_{i=1}^\batchsize\epsilon_i(\param)
\end{equation}
have variance
\begin{equation}
    \label{eq: variance of mini batch}
    \Var(\Loss_\batchsize(\param)) = \Var(\Cost(\param)) + \tfrac1\batchsize\Var(\epsilon_1(\param))
    \overset{\text{isotropy}}=
    \ikernel(0) + \tfrac1\batchsize\ikernel_\epsilon(0),
\end{equation}
where we assume \(\Cost\sim\normal(\mu, \ikernel)\) and
\(\epsilon_i\sim\normal(0, \ikernel_\epsilon)\) in the last equation for
simplicity. But this step did not yet require the distributional Gaussian assumption
beyond the mean and variance.

\subsection{Variance estimation}
\label{subsec: non-parametric covariance estimation}

Recall that the asymptotic learning rate \(\lr_\infty\) in Equation~\eqref{eq:
asymptotic learning rate} only depends on \(\ikernel(0)\)
and \(\ikernel'(0)\). So if we estimate these values, we are certain to get
the right RFD step sizes asymptotically without knowing the entire covariance
kernel \(\ikernel\).

Equation~\eqref{eq: variance of mini batch} reveals that
for \(Z_\batchsize := (\Loss_\batchsize(\param) - \mu)^2\) we have
\[
    \E[Z_\batchsize] = \beta_0 + \tfrac1\batchsize \beta_1
    \qquad\text{i.e.}\qquad
    Z_\batchsize = \beta_0 + \tfrac1\batchsize \beta_1 + \text{noise}
\]
with bias \(\beta_0=\ikernel(0)\) and slope
\(\beta_1=\ikernel_\epsilon(0)\). So a linear regression on samples
\((\tfrac1{\batchsize_k}, Z_{\batchsize_k})_{k\le n}\) allows for the estimation of
\(\beta_0\) and \(\beta_1\). Using the Gaussian assumption from \eqref{eq:
variance of mini batch}, the variance of \(Z_\batchsize\) is the (centered)
fourth moment of \(\Loss_\batchsize\), which is
given by
\[
    \sigma_\batchsize^2 := \Var(Z_\batchsize)
    = \E[Z_\batchsize^4] - \E[Z_\batchsize^2]^2
    =  2 \Var(\Loss_\batchsize(\param))^2
    = 2(\beta_0 + \tfrac1\batchsize \beta_1)^2.
\]
In particular the variance of \(Z_\batchsize\) depends on the batch size \(\batchsize\).
The linear regression is therefore heteroskedastic. Weighted least squares (WLS) 
\citep[e.g.][Theorem~4.2]{kayFundamentalsStatisticalSignal1993} is designed to handle
this case, but for its application the variance of \(Z_\batchsize\) is needed. Since
\(\beta_0,\beta_1\) are the parameters we wish to estimate, we find ourselves in
the paradoxical situation that we need \(\beta\) to obtain \(\beta\). Our solution to
this problem is to start with a guess of \(\beta_0, \beta_1\), apply WLS to
obtain a better estimate and repeat this bootstrapping procedure until
convergence. Since all \(Z_\batchsize\) have the same underlying cost
\(\Cost\), we sample the parameters \(\param\) randomly to reduce their
covariance (details in Sec.~\ref{sec: variance estimation}).

The same procedure can be applied to obtain \(\ikernel'(0)\), where the counterpart of
Equation~\eqref{eq: variance of mini batch} is given by
\[
    \Var(\partial_i \Loss_\batchsize(\param))
    = \Var(\partial_i\Cost(\param)) + \tfrac1\batchsize\Var(\partial_i\epsilon_1(\param))
    \overset{\text{isotropy}}=
    -(\ikernel'(0) + \tfrac1\batchsize\ikernel_\epsilon'(0)).
\]
\begin{remark}
Under the isotropy assumption the partial derivatives are iid, so the expectation of
\(\|\nabla\Loss_\batchsize(\param)\|^2 = \sum_{i=1}^\dims (\partial_i\Loss_\batchsize(\param))^2\)
is this variance scaled by \(\dims\). In particular the variance needs to
scale with \(\frac1\dims\) to keep the gradient norms (and thus the Lipschitz
constant of \(\Cost\)) stable. This observation is closely related to
``isoperimetry'' \parencite[e.g.][]{bubeckUniversalLawRobustness2021}, for details
see \cite{benningGradientSpanAlgorithms2024}. Removing the isotropy
assumption and estimating the variance component-wise is most likely how
``adaptive'' step sizes
\parencite[e.g.][]{duchiAdaptiveSubgradientMethods2011,kingmaAdamMethodStochastic2015},
like the ones used by Adam, work (cf.~Sec.~\ref{sec: geometric anisotropy}).
\end{remark}

\paragraph*{Batch size distribution}

Before we can apply linear regression to the samples  \((\tfrac1{\batchsize_k},
Z_{\batchsize_k})_{k\le n}\), it is necessary to choose the batch sizes
\(\batchsize_k\). As this choice is left to us, we
calculate the variance of our estimator \(\hat{\beta}_0\) of \(\beta_0\)
explicitly (Lemma~\ref{lem: variance of beta_0}), in order to minimize this
variance subject to a sample budget \(\alpha\) over the selection of batch sizes
\begin{equation}
    \label{eq: bsize optimization problem}
    \min_{n, \batchsize_1,\dots, \batchsize_n}\Var(\hat{\beta_0}) 
    \quad\text{s.t.}\quad \underbrace{\sum_{k=1}^n \batchsize_k}_{\text{samples used}}\le \alpha.
\end{equation}
Since this optimization problem is very difficult to solve, we rephrase it in terms of
the empirical distribution of batch sizes \(\nu_n = \frac1n\sum_{i=1}^n
\delta_{\batchsize_i}\). Optimizing over distributions is still difficult, but we
explain in Section~\ref{sec: batch size distribution} how to heuristically arrive at the parametrization
\[
    \nu(b) \propto
    \exp\bigl(\lambda_1 \tfrac1{\sigma_\batchsize^2} - \lambda_2 \batchsize\bigr),
    \qquad b\in \nat
\]
where the parameters \(\lambda_1,\lambda_2\ge 0\) can then be used to optimize
\eqref{eq: bsize optimization problem}. Due to our usage of
\(\sigma_\batchsize^2\) this has to be bootstrapped.

\paragraph*{Covariance estimation} While the variance estimates above ensure correct
asymptotic learning rates, we motivated in Section~\ref{subsec: step size heuristics} that
asymptotic learning rates alone would result in too large step sizes at the
beginning. We therefore use the estimates of \(\ikernel(0)\) and \(\ikernel'(0)\)
to fit a covariance model, effectively acting as a gradient clipper while
retaining the asymptotic guarantees. Note that covariance models with less than
two parameters are generally fully determined by these values.

\subsection{Stochastic RFD (S-RFD)}

It is reasonable to ask whether there is a `stochastic gradient descent'-like
counterpart to the `gradient descent'-like RFD. The answer is yes, and we already
have all the required machinery. 

\begin{restatable}[S-RFD]{extension}{srfd}
    \label{ext: s-rfd}
    For loss \(\Cost\sim\normal(\mu, \ikernel)\) and stochastic errors \(\epsilon_i\overset{\iid}\sim\normal(0, \ikernel_\epsilon)\)
    we have 
    \[
        \argmin_{\step}\E[\Cost(\param - \step) \mid \Loss_\batchsize(\param), \nabla\Loss_\batchsize(\param)]
        = \stepsize^*(\Theta)\tfrac{\nabla\Loss_\batchsize(\param)}{\|\nabla\Loss_\batchsize(\param)\|}
    \]
    with the same step size function \(\stepsize^*\) as for RFD, but modified \(\Theta\)
    \[
        \Theta =
        \frac{
            \ikernel'(0)
        }{
            \ikernel'(0)+ \tfrac1\batchsize\ikernel_\epsilon'(0)
        }
        \frac{
            \ikernel(0)+\tfrac1\batchsize\ikernel_\epsilon(0)
        }{\ikernel(0)}
        \frac{\|\nabla\Loss_\batchsize(\param)\|}{\mu-\Loss_\batchsize(\param)}.
    \]
\end{restatable}
Note, that our non-parametric covariance estimation already provides us with estimates
of \(\ikernel_\epsilon(0)\) and \(\ikernel_\epsilon'(0)\), so no further adaptions are
needed. The resulting asymptotic learning rate is given by
\begin{equation}
    \label{eq: asymptotic lr s-rfd}
    \lr_\infty = \frac{\ikernel(0)+\tfrac1\batchsize\ikernel_\epsilon(0)}{
        (\ikernel'(0)+ \tfrac1\batchsize\ikernel_\epsilon'(0))(\Loss_\batchsize(\param_\infty) - \mu)
    }.
\end{equation}

	\section{MNIST case study}\label{sec: mnist case study}

For our case study we use the negative log likelihood loss to train a neural network
{
	\DeclareFieldFormat{postnote}{#1}
	\parencite[M7]{anEnsembleSimpleConvolutional2020}
}
on the MNIST dataset
\citep{lecunMNISTDATABASEHandwritten2010}. We choose this model as one of the
simplest state-of-the-art models at the time of selection, consisting only of
convolutional layers with ReLU activation interspersed by batch
normalization layers and a single dense layers at the end with softmax
activation.
Assuming isotropy, we estimate \(\mu\), \(\ikernel(0)\) and \(\ikernel'(0)\) as
described in Section~\ref{subsec: non-parametric covariance estimation} and
deduce the parameters \(\sigma^2\) and \(\scale\) of the respective covariance
model (more details in Section~\ref{sec: variance estimation}). We then use the
step sizes listed in Table~\ref{table: optimal step size} for the `squared
exponential' and `rational quadratic' covariance in our RFD algorithm.

\begin{figure}
	\includegraphics*[width=\linewidth]{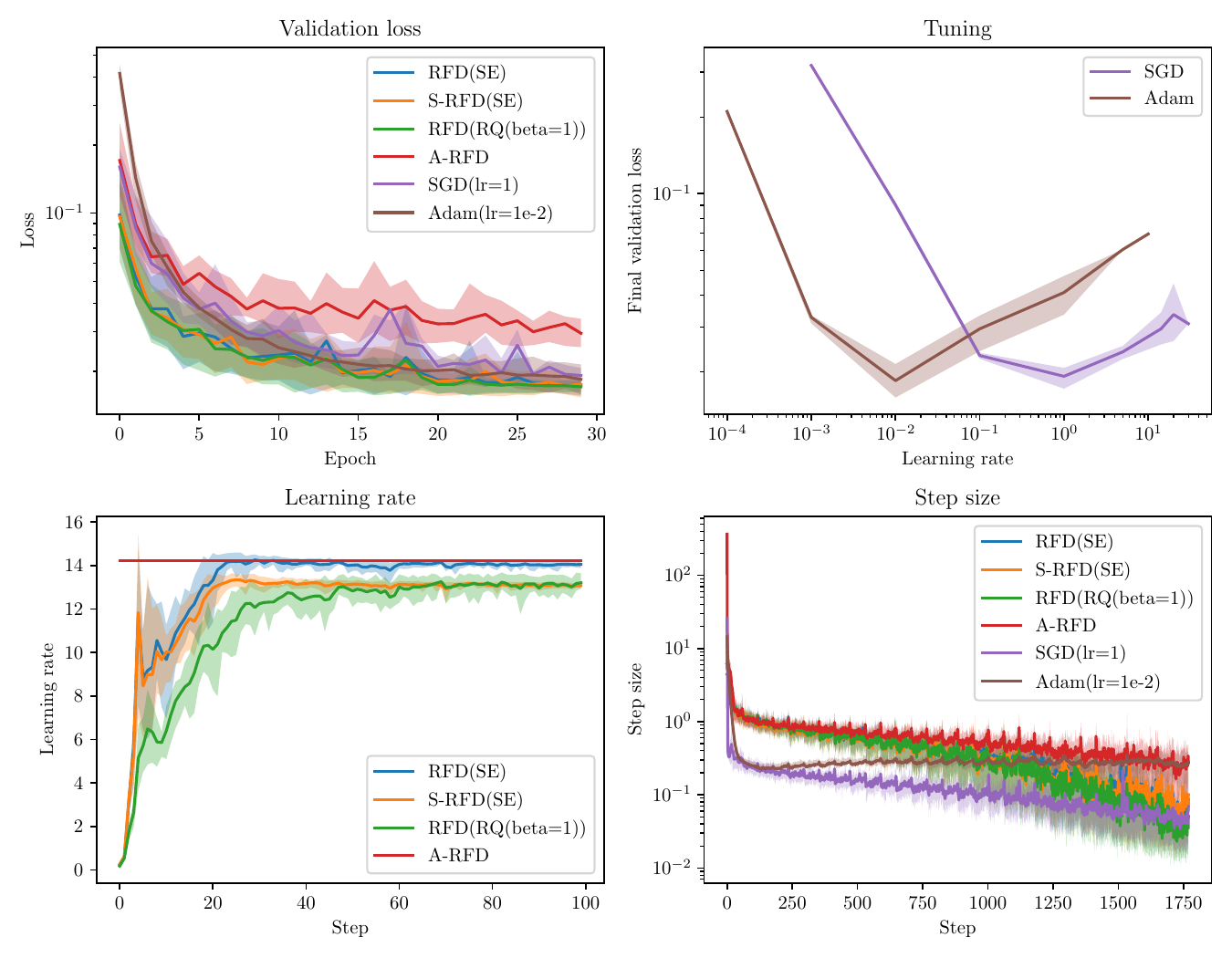}
	\caption{
		Training on the MNIST dataset (batch size \(1024\)). Ribbons describe the
		range between the \(10\%\) and \(90\%\) quantile of \(20\) repeated
		experiments while lines represent their mean. SE stands for the squared exponential \eqref{eq: sqExp
		covariance model} and RQ for the rational quadratic \eqref{eq:
		rational quadratic} covariance. The validation loss uses the test data set, 
		which provides a small advantage to Adam and SGD, as we also use it for tuning.
	}
	\label{fig: summary figure}
\end{figure}

In Figure~\ref{fig: summary figure}, RFD is benchmarked against step size tuned Adam
\citep{kingmaAdamMethodStochastic2015} and stochastic gradient descent (SGD). 
Even with early stopping, their tuning would typically require more than 1 epoch
worth of samples, \emph{in contrast to RFD} (Section~\ref{sec: sampling
efficiency and stability}). We
highlight that A-RFD performs significantly worse than either of the RFD versions
which effectively implement some form of learning rate warmup. This is despite the
RFD learning rates converging to the asymptotic one within one epoch (ca. \(30\) out
of \(60\) steps per epoch). The step sizes on the other hand are (up to noise)
monotonously decreasing. This stands in contrast to the ``wall next to the
minimum'' motivation of gradient clipping.

\textbf{Code availability:}
Our implementation of RFD can be found at \url{https://github.com/FelixBenning/pyrfd} and
the package can also be installed from \href{https://pypi.org/project/pyrfd/}{PyPI} via
`\verb|pip install pyrfd|'.
	\section{Limitations and extensions}
\label{sec: extensions and limitations}

To cover the vast amount of ground that lays between the `formulation of a
general average case optimization problem' and the `prototype of a working
optimizer with theoretical backing',
\begin{enumerate}[leftmargin=1cm]
	\item we used the common 
	\parencite[][]{frazierBayesianOptimization2018,wangBayesianOptimizationBillion2016,wuBayesianOptimizationGradients2017,dauphinIdentifyingAttackingSaddle2014}
	\emph{isotropic} and \emph{Gaussian} distributional assumption for \(\Cost\),
	\item we used very \emph{simple covariance models} for the actual implementation,
	\item we used WLS in our variance estimation procedure despite the \emph{violation of independence}.
\end{enumerate}
Since RFD is defined as the minimizer of an average instead of an upper bound -- 
making it more risk affine --
it naturally loses the improvement guarantee driving classical convergence
proofs. It is therefore impossible to extend classical optimization proofs and new
mathematical theory must be developed. This risk-affinity can also be observed
in its comparatively large step sizes (cf.~Fig.~\ref{fig: summary
figure}
and Sec.~\ref{sec: experiments}). On CIFAR-100
\parencite{krizhevskyLearningMultipleLayers2009}, the step sizes were \emph{too}
large and it is an open question whether assumptions were violated or whether
RFD is simply too risk-affine. But since the variance of random functions
vanishes asymptotic with high dimension
\parencite{benningGradientSpanAlgorithms2024} we highly suspect the former
(cf.~Remark~\ref{rem: high dimension}).

Future work will therefore have to target these assumptions. Some of the assumptions
were already simplifications for the sake of exposition, and we deferred their
relaxation to the appendix. The Gaussian assumption can be relaxed with a
generalization to the `BLUE' (Sec.~\ref{sec: BlUE}),
isotropy can be generalized to `geometric anisotropies' (Sec.~\ref{sec:
geometric anisotropy}) and the risk-affinity of RFD can be reduced with
confidence intervals (Sec.~\ref{sec: conservative rfd}). Since simple random
linear models already violate stationary isotropy (Sec.~\ref{sec: random regression}),
we believe that stationarity is the most important assumption to attack in future
work.


	\section{Conclusion}

In this paper we have demonstrated the \textbf{viability} (computability and
scalable complexity) and \textbf{advantages} (scale invariance, explainable step
size schedule which does not require expensive tuning) of replacing the
classical ``convex function'' framework with the ``random function'' framework.
Along the way we
bridged the gap between Bayesian optimization (not scalable so far) and
classical optimization methods (scalable). This theoretical framework not only
sheds light on existing step size heuristics, but can also be used to develop
future heuristics.

We envision the following improvements to RFD in the future:

\begin{enumerate}
	\item The \emph{reliability} of RFD can be improved by generalizing the
	distributional assumptions to cover more real world scenarios. In particular
	we are interested in the generalization to non-stationary isotropy because we
	suspect that regularization such as weight and batch normalization
	\parencite{salimansWeightNormalizationSimple2016,ioffeBatchNormalizationAccelerating2015}
	are used to patch violations of stationarity (cf. Section~\ref{sec: appendix: input invariance}).

	\item The \emph{performance} of RFD can also be improved. Since RFD is forgetful
	while momentum methods retains some information it is likely fruitful to
	relax the full forgetfulness. Furthermore, we suspect that adaptive learning
	rates \parencite[e.g.][]{duchiAdaptiveSubgradientMethods2011,kingmaAdamMethodStochastic2015}, such as those used by Adam, can be incorporated with
	geometric anisotropies (cf. Sec. E.1). Performance could also be further
	improved by estimating the covariance (locally) online instead of globally at
	the start. Finally, the implementation itself can be made more performant.
\end{enumerate}




	\section*{Acknowledgement}
	We extend our sincere gratitude to our colleagues at the University of
	Mannheim, with special thanks to Rainer Gemulla and Julie Naegelen for
	insightful discussions and invaluable feedback. 
	The Experiments in this work were partially carried out on the compute cluster of
	the state of Baden-Würtemberg (bwHPC).

	\printbibliography[heading=bibintoc]	

@book{adlerRandomFieldsGeometry2007,
  title = {Random {{Fields}} and {{Geometry}}},
  author = {Adler, Robert J. and Taylor, Jonathan E.},
  date = {2007},
  series = {Springer {{Monographs}} in {{Mathematics}}},
  publisher = {Springer New York},
  location = {New York, NY},
  issn = {1439-7382},
  doi = {10.1007/978-0-387-48116-6},
  isbn = {978-0-387-48112-8},
  langid = {english}
}

@article{agnihotriExploringBayesianOptimization2020,
  title = {Exploring {{Bayesian Optimization}}},
  author = {Agnihotri, Apoorv and Batra, Nipun},
  date = {2020-05-05},
  journaltitle = {Distill},
  shortjournal = {Distill},
  volume = {5},
  number = {5},
  pages = {e26},
  issn = {2476-0757},
  doi = {10.23915/distill.00026},
  langid = {english}
}

@online{anEnsembleSimpleConvolutional2020,
  title = {An {{Ensemble}} of {{Simple Convolutional Neural Network Models}} for {{MNIST Digit Recognition}}},
  author = {An, Sanghyeon and Lee, Minjun and Park, Sanglee and Yang, Heerin and So, Jungmin},
  date = {2020-10-04},
  doi = {10.48550/arXiv.2008.10400},
  pubstate = {prepublished}
}

@article{auffingerComplexityGaussianRandom2023,
  title = {Complexity of {{Gaussian Random Fields}} with {{Isotropic Increments}}},
  author = {Auffinger, Antonio and Zeng, Qiang},
  date = {2023-08-01},
  journaltitle = {Communications in Mathematical Physics},
  shortjournal = {Commun. Math. Phys.},
  volume = {402},
  number = {1},
  pages = {951--993},
  issn = {1432-0916},
  doi = {10.1007/s00220-023-04739-0},
  langid = {english}
}

@online{benningGradientSpanAlgorithms2024,
  title = {Gradient {{Span Algorithms Make Predictable Progress}} in {{High Dimension}}},
  author = {Benning, Felix and Döring, Leif},
  date = {2024-10-13},
  doi = {10.48550/arXiv.2410.09973},
  pubstate = {prepublished}
}

@book{borgwardtSimplexMethodProbabilistic1986,
  title = {The {{Simplex Method}}: {{A Probabilistic Analysis}}},
  shorttitle = {The {{Simplex Method}}},
  author = {Borgwardt, Karl Heinz},
  date = {1986-11-01},
  edition = {Softcover reprint of the original 1st ed. 1987 edition},
  publisher = {Springer},
  location = {Berlin Heidelberg},
  isbn = {978-3-540-17096-9},
  langid = {english},
  pagetotal = {282}
}

@inproceedings{bubeckUniversalLawRobustness2021,
  title = {A {{Universal Law}} of {{Robustness}} via {{Isoperimetry}}},
  booktitle = {Advances in {{Neural Information Processing Systems}}},
  author = {Bubeck, Sebastien and Sellke, Mark},
  date = {2021},
  volume = {34},
  eprint = {2105.12806},
  eprinttype = {arXiv},
  eprintclass = {cs, stat},
  pages = {28811--28822},
  publisher = {Curran Associates, Inc.},
  location = {Virtual Event},
  url = {https://proceedings.neurips.cc/paper/2021/hash/f197002b9a0853eca5e046d9ca4663d5-Abstract.html},
  urldate = {2023-09-22}
}

@inproceedings{choKernelMethodsDeep2009,
  title = {Kernel {{Methods}} for {{Deep Learning}}},
  booktitle = {Advances in {{Neural Information Processing Systems}}},
  author = {Cho, Youngmin and Saul, Lawrence},
  date = {2009},
  volume = {22},
  publisher = {Curran Associates, Inc.},
  url = {https://proceedings.neurips.cc/paper/2009/hash/5751ec3e9a4feab575962e78e006250d-Abstract.html},
  urldate = {2023-04-03}
}

@inproceedings{cunhaOnlyTailsMatter2022,
  title = {Only Tails Matter: {{Average-Case Universality}} and {{Robustness}} in the {{Convex Regime}}},
  shorttitle = {Only Tails Matter},
  booktitle = {Proceedings of the 39th {{International Conference}} on {{Machine Learning}}},
  author = {Cunha, Leonardo and Gidel, Gauthier and Pedregosa, Fabian and Scieur, Damien and Paquette, Courtney},
  date = {2022-06-28},
  pages = {4474--4491},
  publisher = {PMLR},
  issn = {2640-3498},
  url = {https://proceedings.mlr.press/v162/cunha22a.html},
  urldate = {2023-11-09},
  eventtitle = {International {{Conference}} on {{Machine Learning}}},
  langid = {english}
}

@inproceedings{dauphinIdentifyingAttackingSaddle2014,
  title = {Identifying and Attacking the Saddle Point Problem in High-Dimensional Non-Convex Optimization},
  booktitle = {Advances in {{Neural Information Processing Systems}}},
  author = {Dauphin, Yann N and Pascanu, Razvan and Gulcehre, Caglar and Cho, Kyunghyun and Ganguli, Surya and Bengio, Yoshua},
  date = {2014},
  volume = {27},
  publisher = {Curran Associates, Inc.},
  location = {Montréal, Canada},
  url = {https://proceedings.neurips.cc/paper/2014/hash/17e23e50bedc63b4095e3d8204ce063b-Abstract.html},
  urldate = {2022-06-10}
}

@inproceedings{defazioLearningRateFreeLearningDAdaptation2023,
  title = {Learning-{{Rate-Free Learning}} by {{D-Adaptation}}},
  booktitle = {Proceedings of the 40th {{International Conference}} on {{Machine Learning}}},
  author = {Defazio, Aaron and Mishchenko, Konstantin},
  date = {2023-07-03},
  eprint = {2301.07733},
  eprinttype = {arXiv},
  eprintclass = {cs.LG},
  pages = {7449--7479},
  publisher = {PMLR},
  issn = {2640-3498},
  url = {https://proceedings.mlr.press/v202/defazio23a.html},
  urldate = {2024-03-31},
  eventtitle = {International {{Conference}} on {{Machine Learning}}},
  langid = {english}
}

@article{deiftConjugateGradientAlgorithm2021,
  title = {The Conjugate Gradient Algorithm on Well-Conditioned {{Wishart}} Matrices Is Almost Deterministic},
  author = {Deift, Percy and Trogdon, Thomas},
  date = {2021-03},
  journaltitle = {Quarterly of Applied Mathematics},
  shortjournal = {Quart. Appl. Math.},
  volume = {79},
  number = {1},
  eprint = {1901.09007},
  eprinttype = {arXiv},
  eprintclass = {cs, math},
  pages = {125--161},
  issn = {0033-569X, 1552-4485},
  doi = {10.1090/qam/1574},
  langid = {english}
}

@article{deuflhardAffineInvariantConvergence1979,
  title = {Affine {{Invariant Convergence Theorems}} for {{Newton}}’s {{Method}} and {{Extensions}} to {{Related Methods}}},
  author = {Deuflhard, P. and Heindl, G.},
  date = {1979-02},
  journaltitle = {SIAM Journal on Numerical Analysis},
  shortjournal = {SIAM J. Numer. Anal.},
  volume = {16},
  number = {1},
  pages = {1--10},
  publisher = {{Society for Industrial and Applied Mathematics}},
  issn = {0036-1429},
  doi = {10.1137/0716001}
}

@article{duchiAdaptiveSubgradientMethods2011,
  title = {Adaptive {{Subgradient Methods}} for {{Online Learning}} and {{Stochastic Optimization}}},
  author = {Duchi, John and Hazan, Elad and Singer, Yoram},
  date = {2011-07-01},
  journaltitle = {The Journal of Machine Learning Research},
  shortjournal = {J. Mach. Learn. Res.},
  volume = {12},
  pages = {2121--2159},
  issn = {1532-4435}
}

@article{elalaouiOptimizationMeanfieldSpin2021,
  title = {Optimization of Mean-Field Spin Glasses},
  author = {El Alaoui, Ahmed and Montanari, Andrea and Sellke, Mark},
  date = {2021-11},
  journaltitle = {The Annals of Probability},
  volume = {49},
  number = {6},
  pages = {2922--2960},
  publisher = {Institute of Mathematical Statistics},
  doi = {10.1214/21-AOP1519}
}

@incollection{frazierBayesianOptimization2018,
  title = {Bayesian {{Optimization}}},
  booktitle = {Recent {{Advances}} in {{Optimization}} and {{Modeling}} of {{Contemporary Problems}}},
  author = {Frazier, Peter I.},
  date = {2018-10},
  series = {{{INFORMS TutORials}} in {{Operations Research}}},
  eprint = {1807.02811},
  eprinttype = {arXiv},
  eprintclass = {cs, math, stat},
  pages = {255--278},
  publisher = {INFORMS},
  location = {Phoenix, Arizona, USA},
  doi = {10.1287/educ.2018.0188},
  isbn = {978-0-9906153-2-3}
}

@article{gaoImplementingNelderMeadSimplex2012,
  title = {Implementing the {{Nelder-Mead}} Simplex Algorithm with~Adaptive Parameters},
  author = {Gao, Fuchang and Han, Lixing},
  date = {2012-01-01},
  journaltitle = {Computational Optimization and Applications},
  shortjournal = {Comput Optim Appl},
  volume = {51},
  number = {1},
  pages = {259--277},
  issn = {1573-2894},
  doi = {10.1007/s10589-010-9329-3},
  langid = {english}
}

@inproceedings{glorotUnderstandingDifficultyTraining2010,
  title = {Understanding the Difficulty of Training Deep Feedforward Neural Networks},
  booktitle = {Proceedings of the {{Thirteenth International Conference}} on {{Artificial Intelligence}} and {{Statistics}}},
  author = {Glorot, Xavier and Bengio, Yoshua},
  date = {2010-03-31},
  pages = {249--256},
  publisher = {{JMLR Workshop and Conference Proceedings}},
  location = {Sardinia, Italy},
  issn = {1938-7228},
  url = {https://proceedings.mlr.press/v9/glorot10a.html},
  urldate = {2023-04-11},
  eventtitle = {Proceedings of the {{Thirteenth International Conference}} on {{Artificial Intelligence}} and {{Statistics}}},
  langid = {english}
}

@book{goodfellowDeepLearning2016,
  title = {Deep {{Learning}}},
  author = {Goodfellow, Ian and Bengio, Yoshua and Courville, Aaron},
  date = {2016-11-10},
  eprint = {omivDQAAQBAJ},
  eprinttype = {googlebooks},
  publisher = {MIT Press},
  isbn = {978-0-262-33737-3},
  langid = {english},
  pagetotal = {801}
}

@report{goyalAccurateLargeMinibatch2018,
  title = {Accurate, {{Large Minibatch SGD}}: {{Training ImageNet}} in 1 {{Hour}}},
  shorttitle = {Accurate, {{Large Minibatch SGD}}},
  author = {Goyal, Priya and Dollár, Piotr and Girshick, Ross and Noordhuis, Pieter and Wesolowski, Lukasz and Kyrola, Aapo and Tulloch, Andrew and Jia, Yangqing and He, Kaiming},
  date = {2018-04-30},
  eprint = {1706.02677},
  eprinttype = {arXiv},
  eprintclass = {cs},
  institution = {arXiv},
  url = {http://arxiv.org/abs/1706.02677},
  urldate = {2024-04-02}
}

@book{hanssonOptimizationLearningControl2023,
  title = {Optimization for Learning and Control},
  author = {Hansson, Anders},
  namea = {Andersen, Martin},
  nameatype = {collaborator},
  date = {2023},
  publisher = {John Wiley \& Sons, Inc.},
  location = {Hoboken, New Jersey},
  isbn = {978-1-119-80914-2},
  langid = {english}
}

@unpublished{hintonNeuralNetworksMachine2012,
  type = {Massive Open Online Course},
  title = {Neural {{Networks}} for {{Machine Learning}}},
  author = {Hinton, Geoffrey},
  namea = {Sirvastava, Nitish and Swersky, Kevin},
  nameatype = {collaborator},
  date = {2012},
  url = {https://www.cs.toronto.edu/~hinton/coursera_lectures.html},
  urldate = {2021-11-16},
  venue = {Coursera}
}

@article{hoareQuicksort1962,
  title = {Quicksort},
  author = {Hoare, C. A. R.},
  date = {1962-01-01},
  journaltitle = {The Computer Journal},
  shortjournal = {The Computer Journal},
  volume = {5},
  number = {1},
  pages = {10--16},
  issn = {0010-4620},
  doi = {10.1093/comjnl/5.1.10}
}

@inproceedings{huangTightLipschitzHardness2022,
  title = {Tight {{Lipschitz Hardness}} for Optimizing {{Mean Field Spin Glasses}}},
  booktitle = {2022 {{IEEE}} 63rd {{Annual Symposium}} on {{Foundations}} of {{Computer Science}} ({{FOCS}})},
  author = {Huang, Brice and Sellke, Mark},
  date = {2022-10},
  pages = {312--322},
  issn = {2575-8454},
  doi = {10.1109/FOCS54457.2022.00037},
  eventtitle = {2022 {{IEEE}} 63rd {{Annual Symposium}} on {{Foundations}} of {{Computer Science}} ({{FOCS}})}
}

@inproceedings{ioffeBatchNormalizationAccelerating2015,
  title = {Batch {{Normalization}}: {{Accelerating Deep Network Training}} by {{Reducing Internal Covariate Shift}}},
  shorttitle = {Batch {{Normalization}}},
  booktitle = {Proceedings of the 32nd {{International Conference}} on {{Machine Learning}}},
  author = {Ioffe, Sergey and Szegedy, Christian},
  date = {2015-06-01},
  eprint = {1502.03167},
  eprinttype = {arXiv},
  eprintclass = {cs},
  pages = {448--456},
  publisher = {PMLR},
  issn = {1938-7228},
  url = {https://proceedings.mlr.press/v37/ioffe15.html},
  urldate = {2021-10-06},
  eventtitle = {International {{Conference}} on {{Machine Learning}}},
  langid = {english}
}

@article{jaynesInformationTheoryStatistical1957,
  title = {Information {{Theory}} and {{Statistical Mechanics}}},
  author = {Jaynes, E. T.},
  date = {1957-05-15},
  journaltitle = {Physical Review},
  shortjournal = {Phys. Rev.},
  volume = {106},
  number = {4},
  pages = {620--630},
  publisher = {American Physical Society},
  doi = {10.1103/PhysRev.106.620}
}

@book{johnsonAppliedMultivariateStatistical2007,
  title = {Applied {{Multivariate Statistical Analysis}}},
  author = {Johnson, Richard Arnold and Wichern, Dean W.},
  date = {2007},
  edition = {6},
  publisher = {Pearson College Div},
  location = {Upper Saddle River, N.J},
  isbn = {978-0-13-187715-3},
  langid = {english},
  pagetotal = {767}
}

@book{kayFundamentalsStatisticalSignal1993,
  title = {Fundamentals of Statistical Signal Processing: Estimation Theory},
  shorttitle = {Fundamentals of Statistical Signal Processing},
  author = {Kay, Steven M.},
  date = {1993-02},
  publisher = {Prentice-Hall, Inc.},
  location = {USA},
  isbn = {978-0-13-345711-7},
  pagetotal = {595}
}

@inproceedings{kingmaAdamMethodStochastic2015,
  title = {Adam: {{A Method}} for {{Stochastic Optimization}}},
  shorttitle = {Adam},
  booktitle = {Proceedings of the 3rd {{International Conference}} on {{Learning Representations}}},
  author = {Kingma, Diederik P. and Ba, Jimmy},
  date = {2015},
  eprint = {1412.6980},
  eprinttype = {arXiv},
  location = {San Diego},
  eventtitle = {{{ICLR}}}
}

@book{klenkeProbabilityTheoryComprehensive2014,
  title = {Probability {{Theory}}: {{A Comprehensive Course}}},
  shorttitle = {Probability {{Theory}}},
  author = {Klenke, Achim},
  date = {2014},
  series = {Universitext},
  publisher = {Springer},
  location = {London},
  doi = {10.1007/978-1-4471-5361-0},
  isbn = {978-1-4471-5360-3 978-1-4471-5361-0},
  langid = {english}
}

@report{krizhevskyLearningMultipleLayers2009,
  title = {Learning Multiple Layers of Features from Tiny Images},
  author = {Krizhevsky, Alex},
  date = {2009},
  url = {https://www.cs.toronto.edu/ kriz/learning-features-2009-TR.pdf},
  urldate = {2024-05-21}
}

@article{kushnerNewMethodLocating1964,
  title = {A {{New Method}} of {{Locating}} the {{Maximum Point}} of an {{Arbitrary Multipeak Curve}} in the {{Presence}} of {{Noise}}},
  author = {Kushner, H. J.},
  date = {1964-03-01},
  journaltitle = {Journal of Basic Engineering},
  shortjournal = {Journal of Basic Engineering},
  volume = {86},
  number = {1},
  pages = {97--106},
  issn = {0021-9223},
  doi = {10.1115/1.3653121}
}

@inproceedings{lacotteOptimalRandomizedFirstOrder2020,
  title = {Optimal {{Randomized First-Order Methods}} for {{Least-Squares Problems}}},
  booktitle = {Proceedings of the 37th {{International Conference}} on {{Machine Learning}}},
  author = {Lacotte, Jonathan and Pilanci, Mert},
  date = {2020-11-21},
  pages = {5587--5597},
  publisher = {PMLR},
  issn = {2640-3498},
  url = {https://proceedings.mlr.press/v119/lacotte20a.html},
  urldate = {2023-11-09},
  eventtitle = {International {{Conference}} on {{Machine Learning}}},
  langid = {english}
}

@dataset{lecunMNISTDATABASEHandwritten2010,
  title = {{{THE MNIST DATABASE}} of Handwritten Digits},
  author = {LeCun, Yann and Cortes, Corinna and Burges, Christopher J.C.},
  date = {2010},
  url = {http://yann.lecun.com/exdb/mnist/},
  urldate = {2024-01-24}
}

@thesis{lizottePracticalBayesianOptimization2008,
  type = {phdthesis},
  title = {Practical Bayesian Optimization},
  author = {Lizotte, Daniel James},
  date = {2008},
  institution = {University of Alberta},
  location = {Alberta, Canada},
  pagetotal = {168}
}

@article{matheronIntrinsicRandomFunctions1973,
  title = {The Intrinsic Random Functions and Their Applications},
  author = {Matheron, G.},
  date = {1973-12},
  journaltitle = {Advances in Applied Probability},
  volume = {5},
  number = {3},
  pages = {439--468},
  publisher = {Cambridge University Press},
  issn = {0001-8678, 1475-6064},
  doi = {10.2307/1425829},
  langid = {english}
}

@article{montanariOptimizationSherringtonKirkpatrick2021,
  title = {Optimization of the {{Sherrington--Kirkpatrick Hamiltonian}}},
  author = {Montanari, Andrea},
  date = {2021-01-07},
  journaltitle = {SIAM Journal on Computing},
  shortjournal = {SIAM J. Comput.},
  pages = {FOCS19-1},
  publisher = {{Society for Industrial and Applied Mathematics}},
  issn = {0097-5397},
  doi = {10.1137/20M132016X},
  pagetotal = {FOCS19-38}
}

@book{nesterovLecturesConvexOptimization2018,
  title = {Lectures on {{Convex Optimization}}},
  author = {Nesterov, Yurii Evgen'evič},
  date = {2018},
  series = {Springer Optimization and {{Its}} Applications; Volume 137},
  edition = {Second edition},
  publisher = {Springer},
  location = {Cham},
  doi = {10.1007/978-3-319-91578-4},
  isbn = {978-3-319-91578-4},
  langid = {english}
}

@inproceedings{padidarScalingGaussianProcesses2021,
  title = {Scaling {{Gaussian Processes}} with {{Derivative Information Using Variational Inference}}},
  booktitle = {Advances in {{Neural Information Processing Systems}}},
  author = {Padidar, Misha and Zhu, Xinran and Huang, Leo and Gardner, Jacob and Bindel, David},
  date = {2021},
  volume = {34},
  pages = {6442--6453},
  publisher = {Curran Associates, Inc.},
  url = {https://proceedings.neurips.cc/paper/2021/hash/32bbf7b2bc4ed14eb1e9c2580056a989-Abstract.html},
  urldate = {2023-05-17}
}

@article{paquetteHaltingTimePredictable2022,
  title = {Halting {{Time}} Is {{Predictable}} for {{Large Models}}: {{A Universality Property}} and {{Average-Case Analysis}}},
  shorttitle = {Halting {{Time}} Is {{Predictable}} for {{Large Models}}},
  author = {Paquette, Courtney and family=Merriënboer, given=Bart, prefix=van, useprefix=true and Paquette, Elliot and Pedregosa, Fabian},
  date = {2022-02},
  journaltitle = {Foundations of Computational Mathematics},
  shortjournal = {Found Comput Math},
  volume = {23},
  number = {2},
  eprint = {2006.04299},
  eprinttype = {arXiv},
  eprintclass = {math, stat},
  pages = {597--673},
  issn = {1615-3383},
  doi = {10.1007/s10208-022-09554-y},
  langid = {english}
}

@article{paquetteUniversalityConjugateGradient2022,
  title = {Universality for the {{Conjugate Gradient}} and {{MINRES Algorithms}} on {{Sample Covariance Matrices}}},
  author = {Paquette, Elliot and Trogdon, Thomas},
  date = {2022-09-01},
  journaltitle = {Communications on Pure and Applied Mathematics},
  volume = {76},
  number = {5},
  pages = {1085--1136},
  issn = {1097-0312},
  doi = {10.1002/cpa.22081},
  langid = {english}
}

@inproceedings{pascanuDifficultyTrainingRecurrent2013,
  title = {On the Difficulty of Training Recurrent Neural Networks},
  booktitle = {Proceedings of the 30th {{International Conference}} on {{Machine Learning}}},
  author = {Pascanu, Razvan and Mikolov, Tomas and Bengio, Yoshua},
  date = {2013-05-26},
  pages = {1310--1318},
  publisher = {PMLR},
  location = {Atlanta},
  issn = {1938-7228},
  url = {https://proceedings.mlr.press/v28/pascanu13.html},
  urldate = {2024-04-02},
  eventtitle = {International {{Conference}} on {{Machine Learning}}},
  langid = {english}
}

@inproceedings{pedregosaAccelerationSpectralDensity2020,
  title = {Acceleration through Spectral Density Estimation},
  booktitle = {Proceedings of the 37th {{International Conference}} on {{Machine Learning}}},
  author = {Pedregosa, Fabian and Scieur, Damien},
  date = {2020-11-21},
  pages = {7553--7562},
  publisher = {PMLR},
  location = {Virtual Event (formerly Vienna)},
  issn = {2640-3498},
  url = {https://proceedings.mlr.press/v119/pedregosa20a.html},
  urldate = {2023-11-09},
  eventtitle = {International {{Conference}} on {{Machine Learning}}},
  langid = {english}
}

@book{rasmussenGaussianProcessesMachine2006,
  title = {Gaussian Processes for Machine Learning},
  author = {Rasmussen, Carl Edward and Williams, Christopher K.I.},
  date = {2006},
  series = {Adaptive Computation and Machine Learning},
  edition = {2},
  number = {3},
  publisher = {MIT Press},
  location = {Cambridge, Massachusetts},
  url = {http://gaussianprocess.org/gpml/chapters/RW.pdf},
  urldate = {2023-04-06},
  isbn = {0-262-18253-X},
  langid = {english},
  pagetotal = {248}
}

@inproceedings{roosHighDimensionalGaussianProcess2021,
  title = {High-{{Dimensional Gaussian Process Inference}} with {{Derivatives}}},
  booktitle = {Proceedings of the 38th {{International Conference}} on {{Machine Learning}}},
  author = {family=Roos, given=Filip, prefix=de, useprefix=false and Gessner, Alexandra and Hennig, Philipp},
  date = {2021-07-01},
  pages = {2535--2545},
  publisher = {PMLR},
  issn = {2640-3498},
  url = {https://proceedings.mlr.press/v139/de-roos21a.html},
  urldate = {2023-05-15},
  eventtitle = {International {{Conference}} on {{Machine Learning}}},
  langid = {english}
}

@inproceedings{salimansWeightNormalizationSimple2016,
  title = {Weight {{Normalization}}: {{A Simple Reparameterization}} to {{Accelerate Training}} of {{Deep Neural Networks}}},
  shorttitle = {Weight {{Normalization}}},
  booktitle = {Advances in {{Neural Information Processing Systems}}},
  author = {Salimans, Tim and Kingma, Durk P},
  date = {2016},
  volume = {29},
  publisher = {Curran Associates, Inc.},
  location = {Barcelona, Spain},
  url = {https://proceedings.neurips.cc/paper/2016/hash/ed265bc903a5a097f61d3ec064d96d2e-Abstract.html},
  urldate = {2023-10-16}
}

@article{shahriariTakingHumanOut2016,
  title = {Taking the {{Human Out}} of the {{Loop}}: {{A Review}} of {{Bayesian Optimization}}},
  shorttitle = {Taking the {{Human Out}} of the {{Loop}}},
  author = {Shahriari, Bobak and Swersky, Kevin and Wang, Ziyu and Adams, Ryan P. and family=Freitas, given=Nando, prefix=de, useprefix=true},
  date = {2016-01},
  journaltitle = {Proceedings of the IEEE},
  volume = {104},
  number = {1},
  pages = {148--175},
  issn = {1558-2256},
  doi = {10.1109/JPROC.2015.2494218},
  eventtitle = {Proceedings of the {{IEEE}}}
}

@inproceedings{smithCyclicalLearningRates2017,
  title = {Cyclical {{Learning Rates}} for {{Training Neural Networks}}},
  booktitle = {2017 {{IEEE Winter Conference}} on {{Applications}} of {{Computer Vision}} ({{WACV}})},
  author = {Smith, Leslie N.},
  date = {2017-03},
  pages = {464--472},
  doi = {10.1109/WACV.2017.58},
  eventtitle = {2017 {{IEEE Winter Conference}} on {{Applications}} of {{Computer Vision}} ({{WACV}})}
}

@online{smithDisciplinedApproachNeural2018,
  title = {A Disciplined Approach to Neural Network Hyper-Parameters: {{Part}} 1 -- Learning Rate, Batch Size, Momentum, and Weight Decay},
  shorttitle = {A Disciplined Approach to Neural Network Hyper-Parameters},
  author = {Smith, Leslie N.},
  date = {2018-04-24},
  eprint = {1803.09820},
  eprinttype = {arXiv},
  eprintclass = {cs, stat},
  doi = {10.48550/arXiv.1803.09820},
  pubstate = {prepublished}
}

@book{steinInterpolationSpatialData1999,
  title = {Interpolation of {{Spatial Data}}},
  author = {Stein, Michael L.},
  date = {1999},
  series = {Springer {{Series}} in {{Statistics}}},
  publisher = {Springer},
  location = {New York, NY},
  doi = {10.1007/978-1-4612-1494-6},
  isbn = {978-1-4612-7166-6 978-1-4612-1494-6}
}

@article{subagFollowingGroundStates2021,
  title = {Following the {{Ground States}} of {{Full-RSB Spherical Spin Glasses}}},
  author = {Subag, Eliran},
  date = {2021},
  journaltitle = {Communications on Pure and Applied Mathematics},
  volume = {74},
  number = {5},
  pages = {1021--1044},
  issn = {1097-0312},
  doi = {10.1002/cpa.21922},
  langid = {english}
}

@article{wangBayesianOptimizationBillion2016,
  title = {Bayesian {{Optimization}} in a {{Billion Dimensions}} via {{Random Embeddings}}},
  author = {Wang, Ziyu and Hutter, Frank and Zoghi, Masrour and Matheson, David and family=Feitas, given=Nando, prefix=de, useprefix=false},
  date = {2016-02-19},
  journaltitle = {Journal of Artificial Intelligence Research},
  volume = {55},
  pages = {361--387},
  issn = {1076-9757},
  doi = {10.1613/jair.4806},
  langid = {english}
}

@article{williamsBayesianClassificationGaussian1998,
  title = {Bayesian Classification with {{Gaussian}} Processes},
  author = {Williams, C.K.I. and Barber, D.},
  date = {1998-12},
  journaltitle = {IEEE Transactions on Pattern Analysis and Machine Intelligence},
  volume = {20},
  number = {12},
  pages = {1342--1351},
  issn = {1939-3539},
  doi = {10.1109/34.735807},
  eventtitle = {{{IEEE Transactions}} on {{Pattern Analysis}} and {{Machine Intelligence}}}
}

@article{williamsComputationInfiniteNeural1998,
  title = {Computation with {{Infinite Neural Networks}}},
  author = {Williams, Christopher K. I.},
  date = {1998-07-01},
  journaltitle = {Neural Computation},
  shortjournal = {Neural Computation},
  volume = {10},
  number = {5},
  pages = {1203--1216},
  issn = {0899-7667},
  doi = {10.1162/089976698300017412}
}

@inproceedings{wuBayesianOptimizationGradients2017,
  title = {Bayesian {{Optimization}} with {{Gradients}}},
  booktitle = {Advances in {{Neural Information Processing Systems}}},
  author = {Wu, Jian and Poloczek, Matthias and Wilson, Andrew G and Frazier, Peter},
  date = {2017},
  volume = {30},
  publisher = {Curran Associates, Inc.},
  url = {https://proceedings.neurips.cc/paper/2017/hash/64a08e5f1e6c39faeb90108c430eb120-Abstract.html},
  urldate = {2022-06-02}
}

@online{xiaoFashionMNISTNovelImage2017,
  title = {Fashion-{{MNIST}}: A {{Novel Image Dataset}} for {{Benchmarking Machine Learning Algorithms}}},
  shorttitle = {Fashion-{{MNIST}}},
  author = {Xiao, Han and Rasul, Kashif and Vollgraf, Roland},
  date = {2017-09-15},
  eprint = {1708.07747},
  eprinttype = {arXiv},
  eprintclass = {cs, stat},
  doi = {10.48550/arXiv.1708.07747},
  pubstate = {prepublished}
}

@unpublished{zeilerADADELTAAdaptiveLearning2012,
  title = {{{ADADELTA}}: {{An Adaptive Learning Rate Method}}},
  shorttitle = {{{ADADELTA}}},
  author = {Zeiler, Matthew D.},
  date = {2012-12-22},
  eprint = {1212.5701},
  eprinttype = {arXiv},
  eprintclass = {cs}
}

@inproceedings{zhangWhichAlgorithmicChoices2019,
  title = {Which {{Algorithmic Choices Matter}} at {{Which Batch Sizes}}? {{Insights From}} a {{Noisy Quadratic Model}}},
  shorttitle = {Which {{Algorithmic Choices Matter}} at {{Which Batch Sizes}}?},
  booktitle = {Advances in {{Neural Information Processing Systems}}},
  author = {Zhang, Guodong and Li, Lala and Nado, Zachary and Martens, James and Sachdeva, Sushant and Dahl, George and Shallue, Chris and Grosse, Roger B},
  date = {2019},
  volume = {32},
  eprint = {1907.04164},
  eprinttype = {arXiv},
  eprintclass = {cs, stat},
  publisher = {Curran Associates, Inc.},
  url = {https://proceedings.neurips.cc/paper/2019/hash/e0eacd983971634327ae1819ea8b6214-Abstract.html},
  urldate = {2023-11-09}
}

	\clearpage

	\appendix

	\begin{center}
		\textbf{\Large Appendix: Random Function Descent}
	\end{center}

	\section{Experiments}
\label{sec: experiments}

\subsection{Covariance estimation}

\begin{figure}
	\includegraphics*[width=\linewidth]{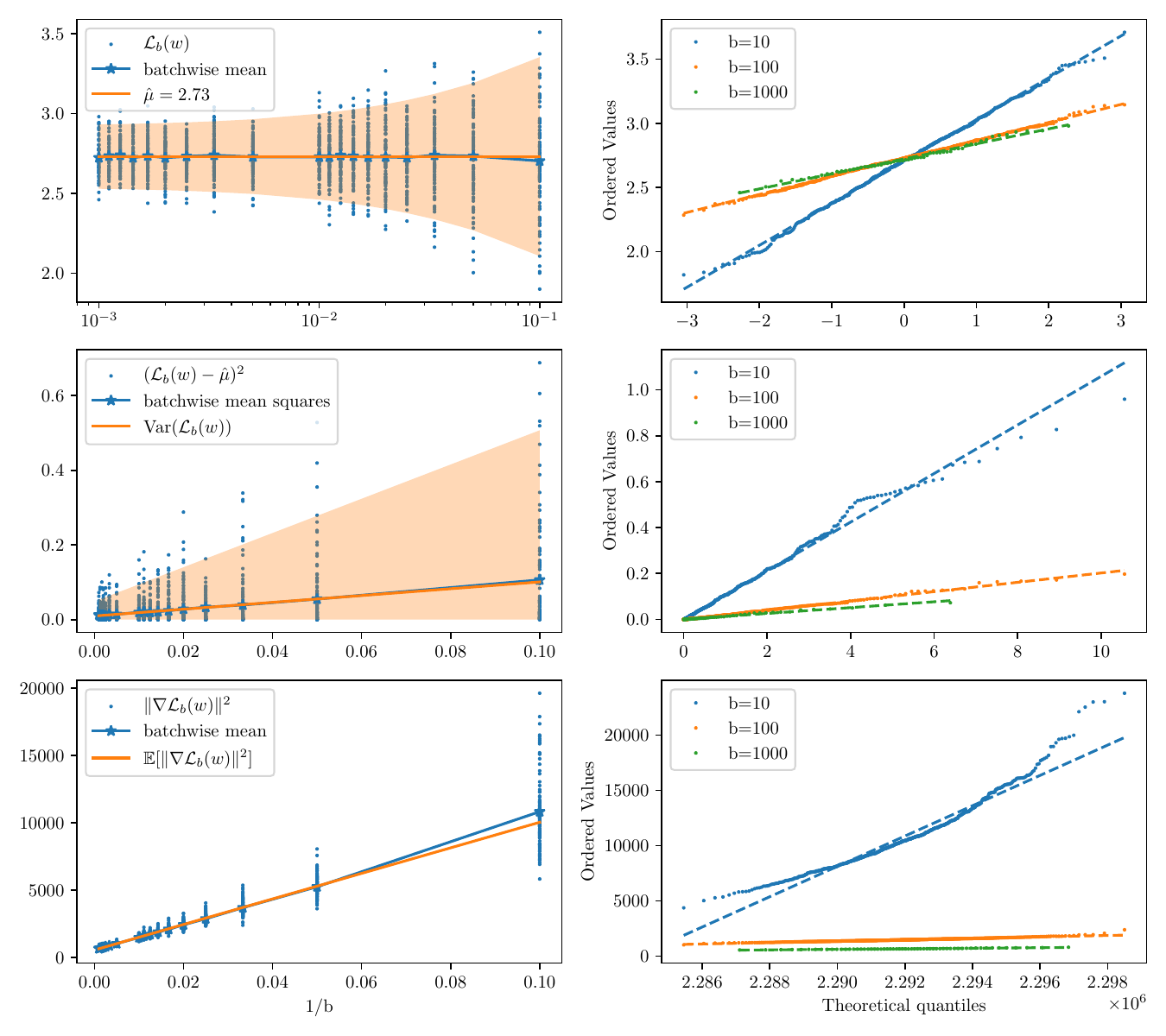}
	\caption{
		Visualization of the variance estimation (Section~\ref{subsec:
		non-parametric covariance estimation}) with \(95\%\)-confidence intervals
		based on the assumed distribution.  Quantile-quantile (QQ) plots of the
		losses (against a normal distribution), squared losses (against a
		\(\chi^2(1)\) distribution) and squared gradient norms (against a
		\(\chi^2(\dims)\)-distribution) are displayed on the right for a selection
		of batch sizes.
	}
	\label{fig: covariance estimation}
\end{figure}

In Figure~\ref{fig: covariance estimation} we visualize weighted
least squares (WLS) regression of the covariance estimation from
Section~\ref{subsec: non-parametric covariance estimation}.
Note, that we sampled much more samples per batch size for these plots than RFD
would typically require by itself in order to be able to plot batch-wise means and
batch-wise QQ-plots. The batch size distribution we described in
Section~\ref{sec: batch size distribution} would avoid sampling the
same batch size multiple times to ensure better stability of the regression and
generally requires much fewer samples than were used for this visualization
(cf.~\ref{sec: sampling efficiency and stability})

We can observe from the QQ-plots on the right, that the Gaussian assumption is
essentially justified for the losses, resulting in a \(\chi^2(1)\) distribution
for the squared losses and a \(\chi^2(\dims)\) distribution for the gradient
norms squared. The confidence interval estimate for the squared norms appears to
be much too small (it is plotted, but too small to be visible). Perhaps this
signifies a violation of the isotropy assumption as the variance of
\[
	\|\nabla\Loss_\batchsize(\param)\|^2
	= \sum_{i=1}^\dims (\partial_i \Loss_\batchsize(\param))^2
\]
does not appear to be the variance of independent \(\chi^2(\dims)\) Gaussian
random variables, and the independence only follows from the isotropy
assumption.

\subsubsection{Sampling efficiency and stability}
\label{sec: sampling efficiency and stability}

\begin{figure}
	\includegraphics*[width=0.5\linewidth]{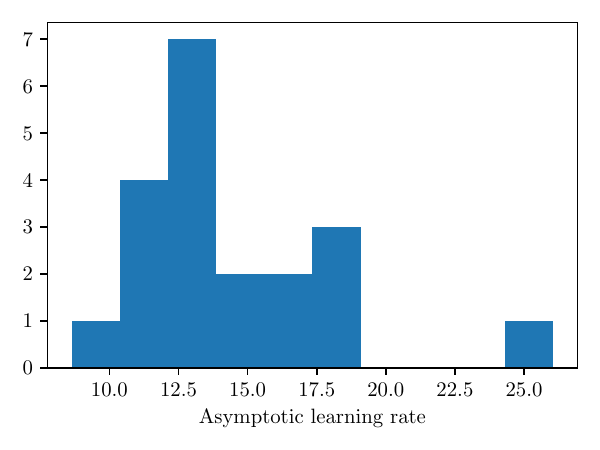}
	\includegraphics*[width=0.5\linewidth]{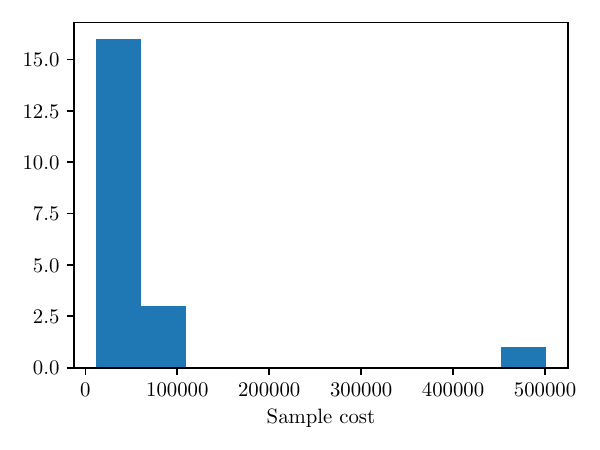}
	\caption{
		\(20\) repeated covariance estimations of model M7
		\parencite{anEnsembleSimpleConvolutional2020} applied to the
		MNIST dataset. On the left are the resulting asymptotic learning rates
		(assuming a final loss of zero) and on the right are the
		samples used until the stopping criterion interrupted sampling. 
	}
	\label{fig: sampling efficiency and stability}
\end{figure}

To evaluate the sampling efficiency and stability of our variance estimation process,
we repeated the covariance estimation of the model model M7
\parencite{anEnsembleSimpleConvolutional2020} applied to the
MNIST dataset \(20\) times (Figure~\ref{fig: sampling efficiency and stability}).
We used a tolerance of \(\mathrm{tol}=0.3\) as a stopping criterion for
the estimated relative standard deviation \eqref{eq: relative std}.

At this tolerance, the asymptotic learning rate already seems relatively stable (in the
same order of magnitude) and the sample cost is quite cheap. The majority of
runs (\(16/20\) runs or \(80\%\)) required less than \(60\,000\) samples (1
epoch). There was one large outlier which used \(500\,589\) samples. A closer
inspection revealed, that after the initial sample to estimate the optimal batch size
distribution, it sampled almost exclusively at batch sizes \(20\) (which was the minimal
cutoff to avoid instabilities caused by batch normalization) and batch sizes between
\(1700\) and \(1900\). It therefore seems like the initial batch of samples caused
a very unfavorable batch size distribution which then required a lot of samples to
recover from. Our selection of an initial sample size of \(6000\) might
therefore have been too small.

A more extensive empirical study is needed to tune this estimation process, but the
process promises to be very sample efficient. Classical step size tuning would
train models for a short duration in order to evaluate the performance of a particular
learning rate \parencite[e.g.][]{smithDisciplinedApproachNeural2018}, but a single
epoch worth of samples is very hard to beat.

Our implementation of this process on the other hand is very inefficient as of writing.
Piping data of differing batch sizes into a model is not a standard use case. We
implement this by repeatedly initializing data loaders, which is anything but
performance friendly.

\subsection{Other models and datasets}

\begin{figure}
	\includegraphics*[width=\linewidth]{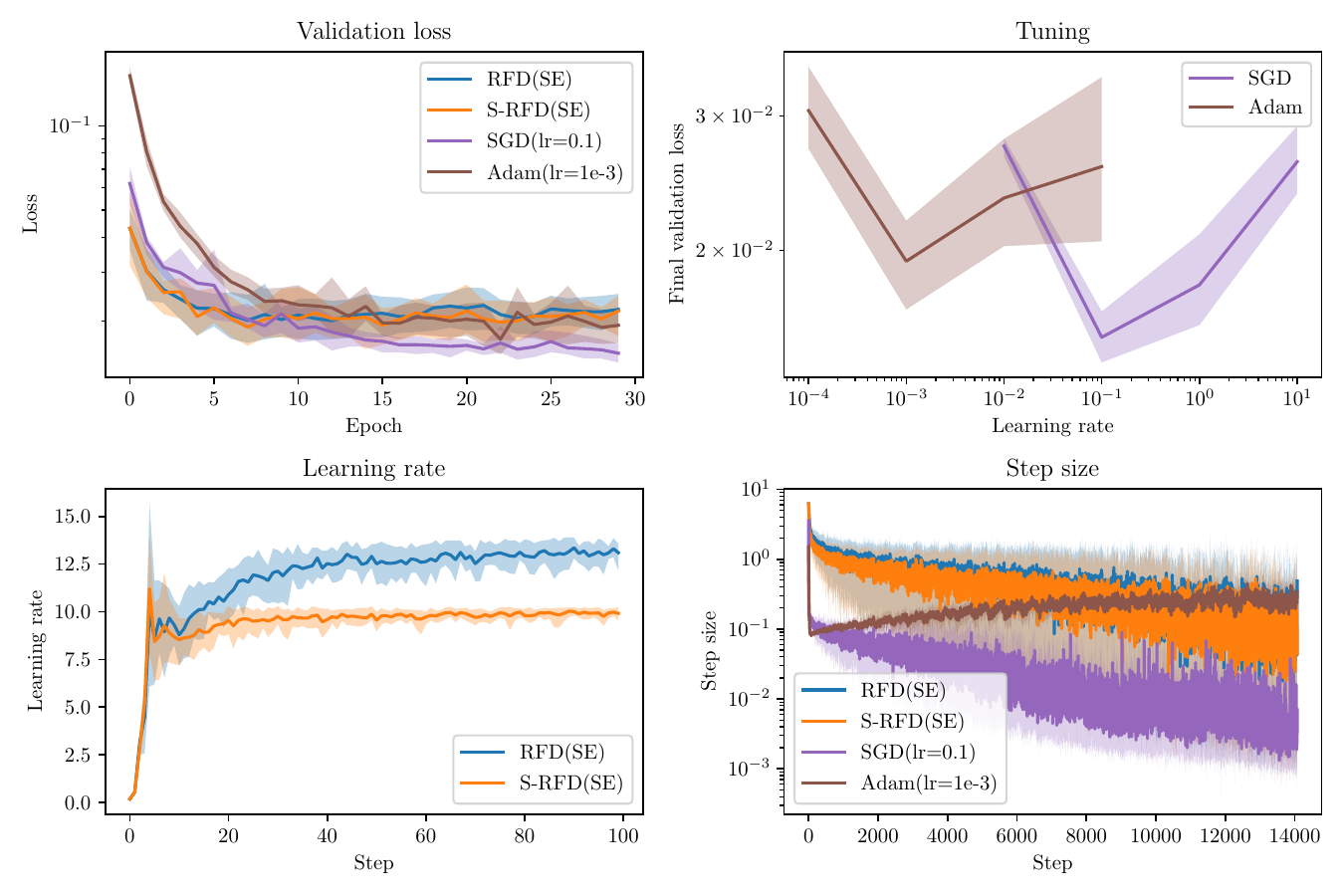}
	\caption{
		Training model M7 \parencite{anEnsembleSimpleConvolutional2020} with batch
		size \(128\) on MNIST \parencite{lecunMNISTDATABASEHandwritten2010}.
	}
	\label{fig: mnist cnn7 bsize=128}
\end{figure}

To estimate the effect of the batch size on RFD, we trained the same model (M7
\parencite{anEnsembleSimpleConvolutional2020}) on MNIST with batch size \(128\)
(Figure~\ref{fig: mnist cnn7 bsize=128}). We can see that the asymptotic learning
rate of S-RFD is reduced at a smaller batch size (cf.~Equation~\ref{eq:
asymptotic lr s-rfd}) but the performance is barely different. Overall, RFD
seems to be slightly too risk-affine, selecting larger step sizes than
the tuned SGD models.

\begin{figure}
	\includegraphics*[width=\linewidth]{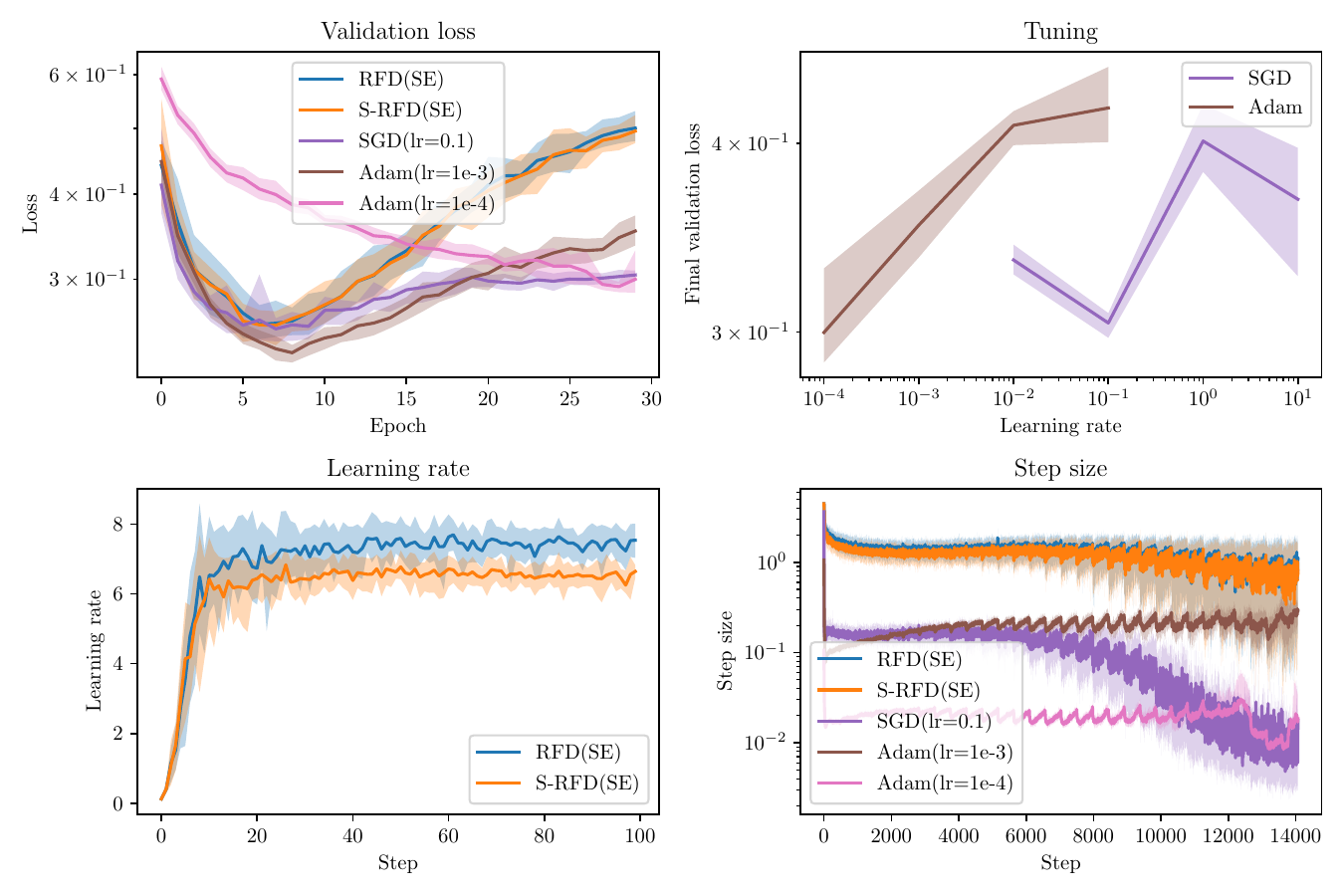}
	\caption{
		Model M5 \parencite{anEnsembleSimpleConvolutional2020} trained on Fashion
		MNIST \parencite{xiaoFashionMNISTNovelImage2017} with batch size \(128\).
	}	
	\label{fig: fashion mnist}
\end{figure}

We also trained a different model (M5 \parencite{anEnsembleSimpleConvolutional2020})
on the Fashion MNIST dataset \parencite{xiaoFashionMNISTNovelImage2017}
with batch size \(128\) (Figure~\ref{fig: fashion mnist}). Since the validation loss increases after
epoch \(5\), early stopping would have been appropriate. We therefore include
Adam with learning rate \(10^{-3}\), despite Adam with learning rate \(10^{-4}\)
technically performing better at the end of training. We can generally see, that
RFD comes very close to tuned performance at the time early stopping would have been
appropriate. Again, learning rates seem to be slightly too large (risk-affine)
in comparison to tuned SGD.

	\section{Variance estimation in detail}
\label{sec: variance estimation}

Recall that we are interested in the regression
\[
   Z_\batchsize(\param)
	= (\Loss_\batchsize(\param) - \mu)^2
	\sim \beta_0 + \frac1\batchsize \beta_1
\]
where the variance of \(Z_\batchsize\) is given by
\[
    \sigma_{b}^2 = 2(\beta_0 + \tfrac1\batchsize \beta_1)^2.
\]
under the Gaussian assumption on \(\Loss_\batchsize\).

More specifically we for minibatch sizes \((b_k)_{k\le n}\) and
parameter vectors \((\param_k)_{k\le n}\) we want to sample mini batch losses
\[
	\Loss^{(k)} := \Loss_{\batchsize_k}(\param_k)
	= \Cost(\param_k) + \frac1{\batchsize_k} \sum_{i=1}^{\batchsize_k} \epsilon_{k,i}(\param_k)
\]
As the \(\epsilon_{k,i}\) are all conditionally independent and therefore uncorrelated, we
have
\[
	\Cov(\Loss^{(k)}, \Loss^{(l)})
	= \Cov(\Cost(\param_k), \Cost(\param_l))
	= \ikernel\bigl(\tfrac{\|\param_k - \param_l\|^2}2\bigr)
\]
Since the covariance kernel \(\ikernel\) is typically monotonously falling in the
distance of parameters \(\|\param_k - \param_l\|^2\), we want to select them as spaced
out as possible to minimize the covariance of \(\Loss^{(k)}\) (which is the next best
thing to iid samples). Randomly selecting \(\param_i\) with Glorot initialization
\parencite{glorotUnderstandingDifficultyTraining2010} will ensure a good spread.

Note that Glorot initialization places all parameters approximately on the same
sphere. This is because Glorot initialization initializes all parameters independently,
therefore their norm is the sum of independent squares, which converges by the
law of large numbers due to the normalization Glorot uses.
Since stationary isotropy and non-stationary isotropy coincides on the sphere,
this is an important effect to consider (cf.~Section~\ref{sec: appendix: input
invariance}).

What is left, is the selection of the batch sizes \(\batchsize_k\).

\subsection{Batch size distribution}
\label{sec: batch size distribution}

Since we plan to use the data set \((\frac1{\batchsize_k}, \Loss^{(k)})_{k\le n}\)
for weighted least squares (WLS) regression and do not have a selection process for
the batch sizes \(\batchsize_k\) yet, it might be appropriate to select the batch
sizes \(\batchsize_k\) in such a way, that the variance of our estimator \(\hat{\beta}_0\)
of \(\beta_0\) is minimized. Here we choose \(\Var(\hat{\beta}_0)\) and not
\(\Var(\hat{\beta}_1)\) as our optimization target, since \(\beta_0 =
\ikernel(0)\) is used to fit the covariance model, while \(\beta_1 =
\ikernel_\epsilon(0)\) is only required for S-RFD. Without deeper analysis
\(\beta_0\) therefore seems to be more important.

Optimization over \(n\) parameters \(\batchsize_k\) is quite difficult, but we can
simplify this optimization problem by considering the empirical batch size distribution
\[
	\nu_n = \frac1n\sum_{k=1}^n \delta_{\batchsize_k}.
\]
Using a random variable \(\Batchsize\) distributed according to \(\nu_n\), the
total number of sample losses can then be expressed as 
\[
	\sum_{k=1}^n \batchsize_k	= n \E[B] = \text{samples used}.
\]
Under an (unrealistic) independence assumption, the variance
\(\Var(\hat{\beta}_0)\) also has a simple representation in terms of \(\nu_n\)
(Lemma~\ref{lem: variance of beta_0}).  We now want to minimize this variance
subject to compute constraint \(\alpha\)
limiting the number of sample losses we can use resulting in the optimization problem
\begin{equation}
    \label{eq: batch size distribution opt problem}
    \Var(\hat{\beta_0}) = \frac1n 
    \underbrace{\frac{1}{\E[\frac1{\sigma_\Batchsize^2}]}}_{\text{`variance of \(Z_\Batchsize\)'}}
    \underbrace{
        \frac{\E[\frac1{\sigma_\Batchsize^2 \Batchsize^2}]}{\E\Bigl[\frac1{\sigma_\Batchsize^2}
        \bigl(\frac1\Batchsize - \E[\frac1{\Batchsize\sigma_\Batchsize^2 \E[1/\sigma_\Batchsize^2]}]\bigr)^2\Bigr]}
    }_{\text{inverse of the `spread'  of \(\frac1\Batchsize\)}}
    \quad\text{s.t.}\quad n\E[\Batchsize]\le \alpha.
\end{equation}
where we recall that \(\sigma_\Batchsize^2\) is the variance of \(Z_\Batchsize\). So the
inverse of the expectation of its inverse is roughly the average variance of \(Z_\Batchsize\).
The second half is the fraction of a weighted second moment divided by the weighted
variance. Unless the mean is at zero, the former will be larger. In particular
we want a spread of data otherwise the variance would be zero. This is in some
conflict with the variance of \(Z_\Batchsize\).

But first, let us get rid of \(n\). Note that we would always increase \(n\)
until our compute budget is used up, since this always reduces variance. So we
approximately have \(n\E[\Batchsize] = \alpha\). Thus
\[
    \Var(\hat{\beta_0})
	 = \frac{\E[\Batchsize]}{\alpha}
    \frac{1}{\E[\frac1{\sigma_\Batchsize^2}]}
        \frac{\E[\frac1{\sigma_\Batchsize^2 \Batchsize^2}]}{\E\Bigl[\frac1{\sigma_\Batchsize^2}
        \bigl(\frac1\Batchsize - \E[\frac1{\Batchsize\sigma_\Batchsize^2 \E[1/\sigma_\Batchsize^2]}]\bigr)^2\Bigr]}
\]
Since \(\alpha\) is now just resulting in a constant factor, it can be assumed
to be \(1\) without loss of generality. Over batch size distributions \(\nu\) we
therefore want to solve the minimization problem
\begin{equation}
	\label{eq: the optimization target}
	\min_{\nu}	
	\underbrace{
   \frac{\E[\Batchsize]}{\E[\frac1{\sigma_\Batchsize^2}]}
	}_{\text{moments}}
	\underbrace{
		\frac{\E[\frac1{\sigma_\Batchsize^2 \Batchsize^2}]}{\E\Bigl[\frac1{\sigma_\Batchsize^2}
		\bigl(\frac1\Batchsize - \E[\frac1{\Batchsize\sigma_\Batchsize^2 \E[1/\sigma_\Batchsize^2]}]\bigr)^2\Bigr]}
	}_{\text{spread}}
\end{equation}

\begin{example}[If we did not require spread]
	If we were not concerned with the variance of batch sizes, we could select a
	constant \(\Batchsize=\batchsize\).
	Then it is straightforward to minimize the moments factor manually
	\[
		\min_{\batchsize} \frac{\E[\batchsize]}{\E[\frac1{\sigma_{\batchsize}^2}]}
		= \batchsize \sigma_\batchsize^2
		= 2\batchsize(\beta_0 + \frac1\batchsize \beta_1)^2,
	\]
	resulting in \(\frac1\batchsize = \frac{\beta_0}{\beta_1}\).
	In other words: If we did not have to be concerned with the spread of \(\Batchsize\)
	there is one optimal selection to minimize the first factor. But in reality we
	have to trade-off this target with the spread of \(\Batchsize\).
\end{example}

To ensure a good spread of data, we use the maximum entropy distribution for \(B\),
with the moment constraints
\begin{align*}
    \E[-B] &\ge -\frac{\alpha}{n} \quad &&\text{average sample usage}\\
    \E[\frac1{\sigma_\Batchsize^2}] &\ge \theta \quad &&\text{\(Z_\Batchsize\) variance}
\end{align*}
which capture the first factor. Maximizing entropy under moment
constraints is known \parencite{jaynesInformationTheoryStatistical1957}
to result in the Boltzmann (a.k.a. Gibbs) distribution
\[
    \nu(\batchsize)
	 = \Pr(\Batchsize = \batchsize) \propto
    \exp\Bigl(\lambda_1 \frac1{\sigma_\batchsize^2} - \lambda_2 \batchsize\Bigr),
\]
where \(\lambda_1,\lambda_2\) depend on the momentum constraints. We can
now forget the origin of this distribution and use \(\lambda_1,\lambda_2\) as
parameters for the distribution \(\nu\) in Equation~\eqref{eq: the
optimization target} to get close to its minimum. In practice we use a zero
order black box optimizer (Nelder-Mead
\parencite{gaoImplementingNelderMeadSimplex2012}). One could calculate the expectations
of \eqref{eq: the optimization target} under this distribution explicitly and take
manual derivatives with respect to \(\lambda_i\) to investigate this further,
but we wanted to avoid getting too distracted by this tangent. 

We also use the estimated relative standard deviation
\begin{equation}
	\label{eq: relative std}
	\mathrm{rel\_std} = \frac{\sqrt{\widehat{\Var}(\hat{\beta}_0)}}{\hat{\beta}_0}
\end{equation}
as a stopping criterion for sampling. Without extensive testing we found a
tolerance of \(\mathrm{rel\_std}<\mathrm{tol}=0.3\) to be reasonable, cf.~Section~\ref{sec:
sampling efficiency and stability}.

\begin{lemma}[Variance of \(\hat{\beta}_0\) in terms of the empirical batch size distribution]
	\label{lem: variance of beta_0}
	Assuming independence of the samples \(((\frac1{\batchsize_k}),
	Z_{\batchsize_k})_{k\le n}\), the variance of \(\hat{\beta}_0\) is given by
	\[
		\Var(\hat{\beta_0}) = \frac1n 
		\frac{1}{\E[\frac1{\sigma_\Batchsize^2}]}
			\frac{\E[\frac1{\sigma_\Batchsize^2 B^2}]}{\E\Bigl[\frac1{\sigma_\Batchsize^2}
			\bigl(\frac1B - \E[\frac1{B\sigma_\Batchsize^2 \E[1/\sigma_\Batchsize^2]}]\bigr)^2\Bigr]}
	\]
	where \(B\) is distributed according to the empirical batch size distribution
	\(\nu_n = \frac1n\sum_{k=1}^n \delta_{\batchsize_k}\).
\end{lemma}
\begin{proof}
	With the notation \(\sigma_{k}^2 = \sigma^2_{\batchsize_k}\) to describe the
	variance of \(Z_{\batchsize_k}\) it follows from
	\parencite[cf.][Thm.~4.2]{kayFundamentalsStatisticalSignal1993}
	that the variance of the estimator \(\hat{\beta}\) of \(\beta\) using \(n\) samples
	is given by
	\begin{align*}
		\Var(\hat{\beta})
		&= (H^T C^{-1} H)^{-1}
		\\
		&= \frac{1}{
			\bigl(\sum_k \frac1{\sigma_k^2}\bigr)\bigl(\sum_k \frac1{(\sigma_k\batchsize_k)^2}\bigr)
			- (\sum_k \frac1{\sigma_k^2\batchsize_k})^2
		}
		\begin{pmatrix}
			\sum_k \frac1{(\sigma_k\batchsize_k)^2}
			& - \sum_k \frac1{\sigma_k^2\batchsize_k}
			\\
			- \sum_k \frac1{\sigma_k^2\batchsize_k}
			&  \sum_k \frac1{\sigma_k^2}
		\end{pmatrix}
	\end{align*}
	where
	\[
		C := \begin{pmatrix}
			\sigma_{1}^2
			\\
			& \ddots
			\\
			&& \sigma_{n}^2
		\end{pmatrix}
		\qquad
		H := \begin{pmatrix}
			1 & \frac1{\batchsize_1}
			\\
			\vdots
			\\
			1 & \frac1{\batchsize_n}
		\end{pmatrix}.
	\]
	In particular we have
	\[
		\Var(\hat{\beta}_0)= \frac{
			\sum_k \frac1{\sigma_k^2\batchsize_k^2}
		}{
			\bigl(\sum_k \frac1{\sigma_k^2}\bigr)\bigl(\sum_k \frac1{\sigma_k^2\batchsize_k^2}\bigr)
			- (\sum_k \frac1{\sigma_k^2\batchsize_k})^2
		}.
	\]
	With the help of \(\theta :=  \sum_j \frac1{\sigma_j^2}\) and \(\lambda_k :=
	\frac1{\sigma_k^2\theta}\),
	we can reorder the divisor. For this note that since the \(\lambda_k\) sum to \(1\)
	we have
	\begin{align*}
			\sum_k \lambda_k \Bigl(\frac1{\batchsize_k}
			- \sum_j \lambda_j\frac{1}{\batchsize_j}\Bigr)^2
			&=
			\sum_k \lambda_k \Bigl(\frac1{\batchsize_k^2}
			- 2 \frac1{\batchsize_k}\sum_j \lambda_j\frac{1}{\batchsize_j}
			+ \Bigl(\sum_j \lambda_j\frac{1}{\batchsize_j}\Bigr)^2
			\Bigr)
			\\
			&=
			\sum_k \lambda_k \frac1{\batchsize_k^2}
			- 2  \Bigl(
				\sum_k \lambda_k\frac1{\batchsize_k}
			\Bigr)
			+ \Bigl(\sum_k \lambda_k\frac{1}{\batchsize_j}\Bigr)^2
			\\
			&=
			\sum_k \lambda_k \frac1{\batchsize_k^2}
			-  \Bigl(
				\sum_k \lambda_k\frac1{\batchsize_k}
			\Bigr)^2
	\end{align*}
	Where the above is essentially the well known statement \(\E[(Y-\E[Y])^2] =
	\E[Y^2] - \E[Y]^2\) for an appropriate selection of \(Y\). This implies that
	our divisor is given by a weighted variance
	\[
			\theta^2\sum_k \lambda_k \Bigl(\frac1{\batchsize_k}
			- \sum_j \lambda_j\frac{1}{\batchsize_j}\Bigr)^2
			= 
			\theta \sum_k \frac1{\sigma_k^2\batchsize_k^2}
			-  \Bigl(
				\sum_k \frac1{\sigma_k^2\batchsize_k}
			\Bigr)^2,
	\]
	where it is only necessary to plug in the definition of \(\theta\) to see the
	right term is exactly our divisor. Expanding both the enumerator as well as the divisor
	by \(\frac1n\), we obtain
	\[
		\Var(\hat{\beta}_0)
		= \frac1{\theta} \frac{
			\frac1n \sum_{k} \frac1{\sigma_k^2 \batchsize_k^2}
		}{
			\frac1n \sum_{k} \frac1{\sigma_k^2} \bigl(
				\frac1{\batchsize_k} - \sum_{j} \lambda_j \frac1{\batchsize_j}
			\bigr)^2
		}
	\]
	Since \(\theta = n \E[1/\sigma_\Batchsize^2]\) for \(\Batchsize\sim \frac1n\sum_{k=1}^n \delta_{\batchsize_k}\)
	and \(\lambda_k = \frac{1}{n\sigma_k^2 \E[1/\sigma_\Batchsize^2]}\), the above can thus
	be written as
	\[
		\Var(\hat{\beta_0}) = \frac1n 
		\frac{1}{\E[1/\sigma_\Batchsize^2]}
			\frac{\E[\frac1{\sigma_\Batchsize^2 \Batchsize^2}]}{\E\Bigl[\frac1{\sigma_\Batchsize^2}
			\bigl(\frac1\Batchsize - \E[\frac1{\Batchsize\sigma_\Batchsize^2 \E[1/\sigma_\Batchsize^2]}]\bigr)^2\Bigr]},
	\]
	which proves our claim.
\end{proof}

	\section{Covariance models}\label{appendix: covariance models}

In this section we calculate the step sizes of the covariance models listed
in Table~\ref{table: optimal step size} and plotted in Figure~\ref{fig: rfd step sizes}.
Additionally we calculate the asymptotic learning rate of A-RFD and prove an
Assumption of Corollary~\ref{prop: convergence} for the squared exponential covariance
(Prop.~\ref{prop: sq exp is strictly monotonous in xi}).

\subsection{Squared exponential}

The squared exponential covariance function is given by
\begin{equation}
	\label{eq: sqExp covariance model}
	\ikernel\bigl(\tfrac{\|x-y\|^2}2\bigr)
	= \sigma^2 \exp\bigl(-\tfrac{\|x-y\|^2}{2\scale^2}\bigr).
\end{equation}
Note that \(\sigma^2\) will play no role in the step sizes of RFD due to its
scale invariance (cf. Advantage~\ref{advant: scale invariance}).

\begin{theorem}
	Let \(\Cost\sim \normal(\mu, \ikernel)\) where \(\ikernel\) is the
	\textbf{squared exponential} covariance function \eqref{eq: sqExp covariance
	model}, then we have
	\begin{equation*}
		\stepsize^* \frac{\nabla\Cost(\param)}{\|\nabla\Cost(\param)\|}
		= \argmin_{\step} \E[\Cost(\param - \step)\mid \Cost(\param),\grad\Cost(\param)]
	\end{equation*}	
	with RFD step size
	\[
		\stepsize^*
		= \frac{
			\scale^2\|\nabla\Cost(\param)\|
		}{
			\sqrt{\bigl(\frac{\mu-\Cost(\param)}{2}\bigr)^2+\scale^2\|\nabla\Cost(\param)\|^2}
			+ \frac{\mu-\Cost(\param)}{2}
		}.
	\]
\end{theorem}
\begin{proof}
	The covariance function \(\ikernel\) is of the form
	\[
		\ikernel(h)
		= \sigma^2 e^{-\frac{h}{\scale^2}}.
	\]
	By Equation~\eqref{eq: simplified step size opt}
	\[
		\stepsize^*
		= -\argmin_{\stepsize} \tfrac{\ikernel(\frac{\stepsize^2}2)}{\ikernel(0)}
		- \stepsize\tfrac{\ikernel'(\frac{\stepsize^2}2)}{\ikernel'(0)}\Theta.
	\]
	where \(\Theta = \frac{\|\nabla\Cost(\param)\|}{\mu-\Cost(\param)}\). We calculate
	\[
		 -\frac{\ikernel\bigl(\frac{\stepsize^2}2\bigr)}{\sqC(0)}
		- \stepsize \frac{\ikernel'\bigl(\frac{\stepsize^2}2\bigr)}{\ikernel'(0)}\Theta
		= -e^{-\frac{\stepsize^2}{2\scale^2}}
		(1 + \stepsize\Theta).
	\]
	This results in the first order condition
	\[
		0 \overset{!}{=} \frac{\stepsize}{\scale^2} e^{-\frac{\stepsize^2}{2\scale^2}}
		(1 + \stepsize\Theta)
		- e^{-\frac{\stepsize^2}{2\scale^2}}\Theta
		= \frac{e^{-\frac{\stepsize^2}{2\scale^2}}}{\scale^2}
		(\stepsize^2\Theta + \stepsize - \scale^2\Theta).
	\]
	Since the exponential can never be zero, we have to solve a quadratic equation.
	Its solution results in
	\begin{equation}
		\label{eq: unstable formula}
		\stepsize^*(\Theta)
		= \sqrt{
			\bigl(\tfrac1{2\Theta}\bigr)^2 + \scale^2 
		}
		- \tfrac1{2\Theta}.
	\end{equation}
	At this point we could stop, but the result is numerically unstable as it suffers
	from catastrophic cancellation. To solve this issue we set \(x=\frac1{2\Theta}\)
	and reorder
	\[
		\stepsize^*
		= \sqrt{x^2 + \scale^2} - x
		= (\sqrt{x^2 + \scale^2} - x) \frac{\sqrt{x^2 + \scale^2} + x}{\sqrt{x^2 + \scale^2} + x}
		= \frac{\cancel{x^2} + \scale^2 - \cancel{x^2}}{\sqrt{x^2+\scale^2} + x}.
	\]
	Re-substituting \(x=\frac1{2\Theta} = \frac{\mu-\Cost(\param)}{2\|\nabla\Cost(\param)\|}\),
	we finally get
	\[
		\stepsize^*
		= \frac{\scale^2}{\sqrt{\bigl(\frac{\mu-\Cost(\param)}{2\|\nabla\Cost(\param)\|}\bigr)^2+\scale^2} + \frac{\mu-\Cost(\param)}{2\|\nabla\Cost(\param)\|}}
		= \frac{\scale^2\|\nabla\Cost(\param)\|}{\sqrt{\bigl(\frac{\mu-\Cost(\param)}{2}\bigr)^2+\scale^2\|\nabla\Cost(\param)\|^2} + \frac{\mu-\Cost(\param)}{2}}.
		\qedhere
	\]
\end{proof}

\begin{prop}[A-RFD for the Squared Exponential Covariance]
	\label{prop: A-RFD for squared exponential covariance}
	If \(\Cost\) is isotropic with squared exponential covariance \eqref{eq:
	sqExp covariance model}, then the step size of A-RFD is given by
	\[
		\hat{\stepsize}
		= \frac{\scale^2}{\mu-\Cost(\param)} \|\nabla\Cost(\param)\|,
	\]
\end{prop}
\begin{proof}
	By Definition~\ref{def: a-rfd} of A-RFD and \(\Theta=\frac{\|\nabla\Cost(\param)\|}{\mu-\Cost(\param)}\)
	we have
	\[
		\hat{\stepsize}(\Theta)
		= \frac{\ikernel(0)}{-\ikernel'(0)}\Theta
		= \frac{\sigma^2\exp(0)}{\frac{\sigma^2}{\scale^2}\exp(0)} \frac{\|\nabla\Cost(\param)\|}{\mu-\Cost(\param)}
		= \scale^2\frac{\|\nabla\Cost(\param)\|}{\mu-\Cost(\param)}.
		\qedhere
	\]
\end{proof}
\begin{prop}
	\label{prop: sq exp is strictly monotonous in xi}
	If \(\Cost\) is isotropic with squared exponential covariance \eqref{eq:
	sqExp covariance model}, then the RFD step sizes are strictly monotonously increasing
	in \(\Theta\).
\end{prop}
\begin{proof}
	Since we know that \(\Theta \to 0\) implies \(\stepsize^* \sim \hat{\stepsize}\to 0\)
	strict monotonicity of \(\stepsize^*\) in \(\Theta\) is sufficient to show that
	\(\stepsize^*\to 0\) also implies \(\Theta\to 0\). So we take the derivative
	of \eqref{eq: unstable formula} resulting in
	\[
		\frac{d}{d\Theta}\stepsize^*
		= 	\frac{1-\frac{1}{\sqrt{1+\scale^2(2\Theta)^2}}}{2\Theta^2},
	\]
	which is greater zero for all \(\Theta>0\).
\end{proof}

\subsection{Rational quadratic}

The \textbf{rational quadratic} covariance function is given by
\begin{equation}
	\label{eq: rational quadratic}
	\ikernel\bigl(\tfrac{\|x-y\|}{2}\bigr)
	= \sigma^2\left(1+\frac{\|x-y\|^2}{\beta\scale^2}\right)^{-\beta/2}
	\quad \beta > 0.
\end{equation}
It can be viewed as a scale mixture of the squared exponential and converges
to the squared exponential in the limit \(\beta\to\infty\)
\citep[87]{rasmussenGaussianProcessesMachine2006}.

\begin{theorem}[Rational Quadratic]
	For \(\Cost\sim\normal(\mu, \ikernel)\) where \(\ikernel\) is the \textbf{rational
	quadratic covariance} we have for \(\Theta=\frac{\|\nabla\Cost(\param)\|}{\mu-\Cost(\param)}\ge 0\)
	that the RFD step size is given by
	\begin{align*}
		\stepsize^*
		&= \scale \sqrt{\beta} \Root_\stepsize\left(
			-1 + \tfrac{\sqrt{\beta}}{\scale\Theta}\stepsize
			+ (1+\beta)\stepsize^2 + \tfrac{\sqrt{\beta}}{\scale\Theta}\stepsize^3
		\right).
	\end{align*}
	The unique root of the polynomial in \(\stepsize\) can be found either directly 
	with a formula for polynomials of third degree (e.g. using Cardano's method)
	or by bisection as it is contained in \([0, 1/\sqrt{1+\beta}]\). 
\end{theorem}

\begin{proof}
	By Theorem~\ref{thm: explicit rfd}	we have
	\begin{equation*}
		\stepsize^*
		= \argmin_{\stepsize}-\frac{\ikernel\bigl(\frac{\stepsize^2}2\bigr)}{\ikernel(0)}
		-  \stepsize\frac{\ikernel'\bigl(\frac{\stepsize^2}2\bigr)}{\ikernel'(0)}\Theta
	\end{equation*}
	for \(\ikernel(x) = \sigma^2(1+\frac{2x}{\beta\scale^2})^{-\beta/2}\).
	We therefore
	need to minimize
	\begin{equation*}
		f\bigl(\frac\stepsize{\sqrt{\beta}\scale}\bigr) := -\left(1+\frac{\stepsize^2}{\beta\scale^2}\right)^{-\beta/2}
		- \stepsize \left(1+\frac{\stepsize^2}{\beta\scale^2}\right)^{-\beta/2-1}\Theta.
	\end{equation*}
	Substitute in \(\tilde{\stepsize}:=\frac\stepsize{\sqrt{\beta}\scale}\), then the first
	order condition is
	\begin{equation*}
		0\overset!=f'(\tilde{\stepsize})
		= -\frac{d}{d\tilde{\stepsize}}
		\left(1+\tilde{\stepsize}^2\right)^{-\beta/2}
		+ \sqrt{\beta}\scale\tilde{\stepsize}\left(1+\tilde{\stepsize}^2\right)^{-\beta/2-1}\Theta
	\end{equation*}
	Dividing both sides by \(\sqrt{\beta}\scale\Theta\) we get
	\begin{align*}
		0 = \frac{f'(\tilde{\stepsize})}{\sqrt{\beta}\scale\Theta}&= \tfrac\beta2(1+\tilde{\stepsize}^2)^{-\frac\beta2-1}
			2\tilde{\stepsize}
			\tfrac1{\sqrt{\beta}\scale\Theta}
			+ (1+\tilde{\stepsize}^2)^{-\frac\beta2 -2}
			\left[1 + \tilde{\stepsize}^2 - (\tfrac\beta2+1)2\tilde{\stepsize}^2 \right]
		\\
		&= (1+\tilde{\stepsize}^2)^{-\frac\beta2-2}\underbrace{\left[
			\beta\tilde{\stepsize}\tfrac{1}{\sqrt{\beta}\scale\Theta}
			(1+\tilde{\stepsize}^2)
			- [1-\tilde{\stepsize}^2(1+\beta)]
		\right]}_{
			= -1 + \frac{\sqrt{\beta}}{\scale\Theta}\tilde{\stepsize}
			+ (1+\beta)\tilde{\stepsize}^2
			+ \frac{\sqrt{\beta}}{\scale\Theta}\tilde{\stepsize}^3
		}
	\end{align*}
	Since \(\Theta\ge 0\) and \(\beta > 0\) all coefficients of the polynomial are
	positive except for the shift. The polynomial thus starts out at \(-1\) in
	zero and only increases from there. Therefore there exists a unique positive
	critical point which is a minimum.

	At the point \(\tilde{\stepsize} = \sqrt{1+\beta}\) the quadratic
	term is already larger than \(1\) so the polynomial is positive and we have
	passed the root. The minimum is
	therefore contained in the interval \([0, \sqrt{1+\beta}]\).
	
	After finding the minimum in \(\tilde{\stepsize}\) we return to \(\stepsize\)
	by multiplication with \(\sqrt{\beta}s\).
\end{proof}

\begin{prop}[A-RFD for the Rational Quadratic Covariance]
	If \(\Cost\) is isotropic with rational quadratic covariance \eqref{eq:
	rational quadratic}, then the step size of A-RFD is given by
	\[
		\hat{\stepsize} = \frac{\scale^2}{\mu-\Cost(\param)} \|\nabla\Cost(\param)\|.
	\]
\end{prop}
\begin{proof}
	\(\ikernel(x) = \sigma^2(1+\frac{2x}{\beta\scale^2})^{-\beta/2}\)
	implies by Definition~\ref{def: a-rfd} of A-RFD and \(\Theta=\frac{\|\nabla\Cost(\param)\|}{\mu-\Cost(\param)}\)
	\[
		\hat{\stepsize}(\Theta)
		= \frac{\ikernel_\Cost(0)}{-\ikernel_\Cost'(0)}\Theta 
		= \frac{\sigma^2(1+ 0)^{-\beta/2}}{\frac{\sigma^2}{\scale^2}(1+0)^{-\beta/2-1}} \frac{\|\nabla\Cost(x)\|}{\mu-\Cost(x)}
		= \scale^2\frac{\|\nabla\Cost(x)\|}{\mu-\Cost(x)}.
		\qedhere
	\]
\end{proof}

\subsection{Matérn}

\begin{definition}
	The Matérn model parametrized by \(\scale >0, \nu\ge 0, \sigma^2\ge 0\) is given by
	\begin{equation}
		\label{eq: matern model}
		\ikernel\bigl(\tfrac{\|x-y\|^2}{2}\bigr)
		= \sigma^2 \frac{2^{1-\nu}}{\Gamma(\nu)}
		\left(\tfrac{\sqrt{2\nu}\|x-y\|}{\scale}\right)^\nu
		\modifiedBessel\left(\tfrac{\sqrt{2\nu}\|x-y\|}{\scale}\right)
	\end{equation}
	where \(\modifiedBessel\) is the modified Bessel function. 
	
	For \(\nu=p+\frac12\) with \(p\in\nat_0\), it can be simplified
	\citep[cf.][sec.~4.2.1]{rasmussenGaussianProcessesMachine2006} to
	\[
		\ikernel\bigl(\tfrac{\|x-y\|^2}{2}\bigr)
		= \sigma^2 e^{- \tfrac{\sqrt{2\nu} \|x-y\|}{s}}\tfrac{p!}{(2p)!}
		\sum_{k=0}^p \tfrac{(2p-k)!}{(p-k)!k!}\left(\tfrac{2\sqrt{2\nu}}{s}\|x-y\|\right)^k
	\]
\end{definition}

The Matérn model encompasses \citet{rasmussenGaussianProcessesMachine2006}
\begin{itemize}
	\item the \textbf{nugget effect} for \(\nu=0\) (independent randomness)
	\item the \textbf{exponential model} for \(\nu=\tfrac12\) (Ornstein-Uhlenbeck process)
	\item the \textbf{squared exponential model} for \(\nu\to \infty\) with the same
	scale \(\scale\) and variance \(\sigma^2\).
\end{itemize}
The random functions induced by the Matérn model are a.s.
\(\lfloor \nu\rfloor\)-times differentiable \citet{rasmussenGaussianProcessesMachine2006}, i.e. the smoothness of the model
increases with increasing \(\nu\). While the exponential covariance model with
\(\nu=\tfrac12\) results in a random function
which is not yet differentiable, larger \(\nu\) result in increasing
differentiability. As differentiability starts with \(\nu=\frac32\) and we have
a more explicit formula for \(\nu=p+\tfrac12\) the cases \(\nu=\frac32\) and
\(\nu=\frac52\) are
of particular interest.
\begin{quote}
	``[F]or \(\nu \ge 7/2\), in the absence of explicit prior knowledge about the existence
	of higher order derivatives, it is probably very hard from finite noisy
	training examples to distinguish between values of \(\nu \ge 7/2\) (or even to
	distinguish between finite values of \(\nu\) and \(\nu \to\infty\), the smooth squared
	exponential, in this case)'' \citep[85]{rasmussenGaussianProcessesMachine2006}.
\end{quote}

\begin{theorem}
	Assuming \(\Cost\sim\normal(\mu, \ikernel)\) is a random function where
	\(\ikernel\) is the Matérn covariance such that \(\nu = p+\tfrac12\) with
	\(p\in \{1,2\}\). Then the RFD step
	is given for \(\Theta := \frac{\|\nabla\Cost(\param)\|}{\mu-\Cost(\param)}\ge 0\) by
	\begin{itemize}
		\item \(p=1\)
		\begin{equation*}
			\stepsize^*
			= \frac{s}{\sqrt{3}}\frac{1}{\left(
				1 + \frac{\sqrt{3}}{s\Theta}
			\right)}
		\end{equation*}

		\item \(p=2\)
		\begin{equation*}
			\stepsize^*
			= \frac{\scale}{\sqrt{5}}\frac{
				(1-\zeta)+\sqrt{4 + (1+\zeta)^2}
			}{2(1+\zeta)} \qquad \zeta := \frac{\sqrt{5}}{3\scale\Theta}.
		\end{equation*}
	\end{itemize}
\end{theorem}

\begin{proof}
	We define \(\C(\stepsize) := \ikernel(\frac{\stepsize^2}2)\), which implies
	\[
		\C'(\stepsize) = \ikernel'\bigl(\tfrac{\stepsize^2}2\bigr)\stepsize 
	\]
	or conversely
	\begin{equation}
		\label{eq: quadratic representation from absolute}
		\ikernel'\bigl(\tfrac{\stepsize^2}2\bigr) = \frac{1}\stepsize \C'(\stepsize).
	\end{equation}
	By Theorem~\ref{thm: explicit rfd}, we need to calculate
	\begin{equation}
		\label{eq: reminder minimization theorem}	
        \stepsize^*
        = \argmin_{\stepsize} -\tfrac{\ikernel(\frac{\stepsize^2}2)}{\ikernel(0)}
        - \stepsize\tfrac{\ikernel'(\frac{\stepsize^2}2)}{\ikernel'(0)}\Theta.
	\end{equation} Discarding \(\sigma\) w.l.o.g. due to scale invariance
	(Advantage~\ref{advant: scale invariance}), we have in the case \(p=1\) 
	\[
		\C(\stepsize) = \Bigl(1+\frac{\sqrt{3}}{\scale}\stepsize\Bigr)
		\exp\Bigl(-\frac{\sqrt{3}}{\scale}\stepsize\Bigr).
	\]
	The derivative is then given by
	\[
		\C'(\stepsize) = -\bigl(\tfrac{\sqrt{3}}\scale\bigr)^2\stepsize \exp\bigl(-\tfrac{\sqrt{3}}{\scale}\stepsize\bigr)
	\]
	which implies using \eqref{eq: quadratic representation from absolute}
	\begin{equation}
		\label{eq: derivative matern p=1}
		\ikernel'\bigl(\tfrac{\stepsize^2}2\bigr) = -\bigl(\tfrac{\sqrt{3}}\scale\bigr)^2 \exp\bigl(-\tfrac{\sqrt{3}}{\scale}\stepsize\bigr)
	\end{equation}
	We therefore need to minimize \eqref{eq: reminder minimization theorem} which is given by
	\[
		\argmin_\stepsize
		- \bigl(1+\tfrac{\sqrt{3}}{\scale}\stepsize\bigr)
		\exp\bigl(-\tfrac{\sqrt{3}}{\scale}\stepsize\bigr)
		- \stepsize \exp\bigl(-\tfrac{\sqrt{3}}{\scale}\stepsize\bigr)\Theta
		= \argmin_\stepsize
		- \bigl(1+(\tfrac{\sqrt{3}}{\scale} + \Theta)\stepsize\bigr)
		\exp\bigl(-\tfrac{\sqrt{3}}{\scale}\stepsize\bigr).
	\]
	The first order condition is
	\[
		0\overset!=
		\Bigl(
			\tfrac{\sqrt{3}}{\scale}\bigl(\cancel{1}+(\tfrac{\sqrt{3}}{\scale} + \Theta)\stepsize\bigr)
			- (\cancel{\tfrac{\sqrt{3}}{\scale}} + \Theta)
		\Bigr)
	\]
	which (divided by \(\Theta\) and noting that the exponential can never be zero) is equivalent to
	\[
			0\overset!=\tfrac{\sqrt{3}}{\scale}(\tfrac{\sqrt{3}}{\scale\Theta} + 1)\stepsize - 1
	\]
	reordering for \(\stepsize\) implies
	\[
		\stepsize \overset! = \frac{\scale}{\sqrt{3}} \frac1{\bigl(1+\tfrac{\sqrt{3}}{\scale\Theta}\bigr)}.
	\]
	It is also not difficult to see that this is the point where the derivative switches
	from negative to positive (i.e. a minimum).

	Let us now consider the case \(p=2\), i.e.
	\[
		\C(\stepsize)
		= \bigl(1 + \tfrac{\sqrt{5}}{\scale}\stepsize + \tfrac{5}{3\scale^2}\stepsize^2\bigr)
		\exp\bigl(-\tfrac{\sqrt{5}}{\scale}\stepsize\bigr),
	\]
	which results in
	\[
		\C'(\stepsize)
		= -\tfrac{5}{3\scale^2}\bigl(\stepsize+\tfrac{\sqrt{5}}{\scale}\stepsize^2\bigr)
		\exp\bigl(-\tfrac{\sqrt{5}}{\scale}\stepsize\bigr),
	\]
	i.e. by \eqref{eq: quadratic representation from absolute}
	\begin{equation}
		\label{eq: derivative matern p=2}
		\ikernel'\bigl(\tfrac{\stepsize^2}2\bigr)
		= -\tfrac{5}{3\scale^2}\bigl(1+\tfrac{\sqrt{5}}{\scale}\stepsize\bigr)
		\exp\bigl(-\tfrac{\sqrt{5}}{\scale}\stepsize\bigr).
	\end{equation}
	We therefore need to minimize \eqref{eq: reminder minimization theorem} which is given by
	\[
		\underbrace{
			\Bigl(-\bigl(1 + \tfrac{\sqrt{5}}{\scale}\stepsize + \tfrac{5}{3\scale^2}\stepsize^2\bigr)
			- \stepsize\bigl(1+\tfrac{\sqrt{5}}{\scale}\stepsize\bigr)\Theta\Bigr)
		}_{
			= -\bigl(1 + \bigl(\tfrac{\sqrt{5}}{\scale}+\Theta\bigr)\stepsize
			+ \bigl(\tfrac{5}{3\scale^2} +\tfrac{\sqrt{5}}{\scale}\Theta\bigr)\stepsize^2
			\bigr)
		}
		\exp\bigl(-\tfrac{\sqrt{5}}{\scale}\stepsize\bigr).
	\]
	The first order condition results in
	\begin{align*}
		0 \overset{!}&{=}	
		\tfrac{\sqrt{5}}\scale\Bigl(
			\cancel{1}
			+ \bigl(\tfrac{\sqrt{5}}{\scale}+\Theta\bigr)\stepsize
			+ \bigl(\tfrac{5}{3\scale^2} + \tfrac{\sqrt{5}}{\scale}\Theta\bigr)\stepsize^2
		\Bigr)
		- \Bigl(
			\bigl(\cancel{\tfrac{\sqrt{5}}{\scale}}+\Theta\bigr)
			+ 2\bigl(\tfrac{5}{3\scale^2} +\tfrac{\sqrt{5}}{\scale}\Theta\bigr)\stepsize
		\Bigr)
		\\
		&= -\Theta + \bigl(\tfrac{5}{3\scale^2} - \tfrac{\sqrt{5}}\scale\Theta\bigr)\stepsize
		+\tfrac{\sqrt{5}}\scale\bigl(\tfrac{5}{3\scale^2}+\tfrac{\sqrt{5}}{\scale}\Theta\bigr)\stepsize^2
	\end{align*}
	Dividing everything by \(\Theta\) and using \(\zeta := \frac{\sqrt{5}}{3\scale\Theta}\) we get
	\[
		0\overset!= -1 -\bigl(\zeta - 1\bigr)\bigl(\tfrac{\sqrt{5}}\scale\stepsize\bigr)
		+ \bigl(\zeta + 1\bigr)\bigl(\tfrac{\sqrt{5}}\scale\stepsize\bigr)^2
	\]
	Taking a closer look at the sign changes of the derivative it becomes
	obvious, that the positive root is the
	minimum, i.e.
	\[
		\frac{\sqrt{5}}{\scale} \stepsize
		\overset{!}= \frac{
			(1-\zeta)+\sqrt{(1-\zeta)^2 +4(1+\zeta)}
		}{2(1+\zeta)}
		= \frac{
			(1-\zeta)+\sqrt{4 + (1+\zeta)^2}
		}{2(1+\zeta)}.
		\qedhere
	\]
\end{proof}

\begin{prop}[A-RFD for the Matérn Covariance]
	If \(\Cost\) is isotropic with Matérn covariance \eqref{eq: matern model}
	such that \(\nu = p+\tfrac12\), then the step size of A-RFD for \(p\in
	\{1,2\}\) is given by
	\begin{itemize}
		\item \(p=1\)
		\[
			\hat{\stepsize}
			= \frac{\scale^2}{3} \frac{\|\nabla\Cost(x)\|}{\mu - \Cost(x)}
		\]
		\item \(p=2\)
		\[
			\hat{\stepsize}
			= \frac{3\scale^2}{5} \frac{\|\nabla\Cost(x)\|}{\mu - \Cost(x)}
		\]
	\end{itemize}
\end{prop}
\begin{proof}
	Noting \(\Theta = \frac{\|\nabla\Cost(x)\|}{\mu - \Cost(x)}\), we have by
	Definition~\ref{def: a-rfd} of A-RFD for \(p=1\)
	\[
		\hat{\stepsize}
		= \frac{\ikernel(0)}{-\ikernel'(0)} \Theta
		\overset{\eqref{eq: derivative matern p=1}}=
		\frac{\scale^2}{3}\Theta,
	\]
	and in the case \(p=2\)	
	\[
		\hat{\stepsize}
		= \frac{\ikernel(0)}{-\ikernel'(0)} \Theta
		\overset{\eqref{eq: derivative matern p=2}}=
		\frac{3\scale^2}{5} \Theta.
		\qedhere
	\]
\end{proof}

	\section{Proofs}
\label{sec: proofs}

In this section we prove all the claims made in the main body.

\subsection{Section~\ref{sec: rfd}: Random function descent}

\subsubsection{Formal RFD}\label{sec: formal rfd}

As we mentioned in a footnote at the definition of RFD, the fact that the parameters
become random variables as they are selected by random gradients poses some mathematical
challenges which would have been distracting to address in the main body.
In following paragraphs leading up to Definition~\ref{def: formal rfd} we introduce
and discuss the probability theory required to provide a mathematically
sound definition.

For a fixed cost distribution \(\Pr_\Cost\) and any weight vectors \(\param\)
and \(\tilde{\param}\) the conditional distribution
\[
	\E[\Cost(\tilde{\param}) \mid \Cost(\param), \nabla\Cost(\param)]
\]
is by its axiomatic definition a \((\Cost(\param),
\nabla\Cost(\param))\)-measurable random variable. By the factorization lemma
\parencite[Cor.~1.9.7]{klenkeProbabilityTheoryComprehensive2014},
there therefore exists a measurable function \((j,g)\mapsto \varphi_{\param, \tilde{\param}}(j, g)\) such that
the following equation holds almost surely
\begin{equation}
	\label{eq: factorization of conditional expec}
	\varphi_{\param, \theta}(\Cost(\param), \nabla\Cost(\param))
	= \E[\Cost(\tilde{\param}) \mid \Cost(\param), \nabla\Cost(\param)].
\end{equation}
Since it is possible to calculate \(\varphi_{\param, \tilde{\param}}\)
explicitly in the Gaussian case (cf.~\ref{thm: conditional gaussian
distribution}), the function  
\[
		\Phi_{\Pr_\Cost}:
		\begin{cases}
			(\real^\dims \times \real \times \real^\dims) \to \real\\
			(\param, j, g) \mapsto \argmin_{\tilde{\param}} \varphi_{\param, \tilde{\param}}(j, g),
		\end{cases}
\]
which implements some tie-breaker rules for set valued \(\argmin\) is measurable
when \(\Cost\) is Gaussian and its covariance function is sufficiently smooth.
To prove measurability in the general case is a difficult problem of its own,
which we do not attempt to solve here, since we would not utilize the conditional
expectation outside of the Gaussian case anyway (cf.~Section~\ref{sec: BlUE}). For deterministic \(\param\), we
therefore have
\begin{align*}
	\Phi_{\Pr_\Cost}(\param, \Cost(\param), \nabla\Cost(\param))
	&= \argmin_{\tilde{\param}} \varphi_{\param, \tilde{\param}}(
		\Cost(\param), \nabla\Cost(\param)
	)
	\\
	\overset{\eqref{eq: factorization of conditional expec}}&=
	\argmin_{\tilde{\param}} \E[
		\Cost(\tilde{\param}) \mid \Cost(\param), \nabla\Cost(\param)
	].
\end{align*}
So if the parameter vectors \(\param_n\) were deterministic, our formal definition of
RFD and our initial definition would coincide.
But for random weights \(\Param\) \eqref{eq: factorization of conditional
expec} stops to hold in general\footnote{
	E.g. consider the random variable
	\[
		\Param
		= \bigl(\argmin \Cost\bigr)\ind_{\Cost(\tilde{\param}) > 0}
		+  \bigl(\argmax \Cost \bigr)\ind_{\Cost(\tilde{\param}) < 0}.
	\]
	In this case, \(\Cost(\Param)\) is much more informative of \(\Cost(\tilde{\param})\) than
	\(\Cost(\param)\) at some deterministic \(\param\).
}, i.e.
\[
	\varphi_{\Param, \tilde{\param}}(\Cost(\Param), \nabla\Cost(\Param))
	\neq \E[\Cost(\tilde{\param})\mid \Cost(\Param), \nabla\Cost(\Param)].
\]
If this equation does not need to hold, we similarly have in general
\[
	\Phi_{\Pr_\Cost}(\Param, \Cost(\Param), \nabla\Cost(\Param))
	\neq \argmin_{\tilde{\param}} \E[\Cost(\tilde{\param}) \mid \Cost(\Param), \nabla\Cost(\Param)].
\]
So the following definition is not just a restatement of the original definition of RFD.

\begin{definition}[Formal RFD]
	\label{def: formal rfd}
	For a Gaussian random cost function \(\Cost\), we define the RFD algorithm
	with starting point \(\Param_0 = \param_0 \in \real^\dims\) by
	\[
		\Param_{n+1} := \Phi_{\Pr_\Cost}(\Param_n, \Cost(\Param_n), \nabla\Cost(\Param_n))
	\]
\end{definition}
This is what we effectively do in Theorem~\ref{thm: explicit rfd} under the
additional isotropy assumption, where we calculate the \(\argmin\) under the
assumption that \(\param\) is deterministic (i.e.  we determine
\(\Phi_{\Pr_\Cost}\)), before we plug-in the random variables
\(\Param_n\) to obtain \(\Param_{n+1}\). Similarly this is how the step size
prescriptions of RFD actually work. We first assume deterministic weights and later
plug the random variables into our formulas. For this reason, we avoided large letters
indicating random variables for parameters \(\param\) in the main body.

\subsubsection{Scale invariance}

\scaleInvariance*

Before we get to the proof, let us quickly formulate the statement in mathematical terms.
Let \(\param_n\) be the parameters selected optimizing \(\Cost\) starting in
\(\param_0\) and \(\tilde{\param}_n\) the parameters selected by the same optimizer
optimizing \(\tilde{\Cost}\) starting in \(\tilde{\param}_0\).

If we apply affine linear scaling to cost \(\Cost\) such that
\(\tilde{\Cost}(\param)= a\Cost(\param) + b\) and start optimization in the same
point, i.e. \(\param_0 = \tilde{\param}_0\), then we expect a scale invariant optimizer
to select
\[
	\param_n = \tilde{\param}_n.
\]
If we scale inputs on the other hand (or more generally map them with a bijection \(\phi\)),
then we expect for \(\tilde{\Cost} := \Cost \circ \phi\) and starting point
\(
	\tilde{\param}_0 = \phi^{-1}(\param_0)
\),
that this relationship is retained by a scale invariant optimizer, i.e.
\[
	\tilde{\param}_n = \phi^{-1}(\param_n).
\]
Why do we use a different starting point? As an illustrating example, assume that
\(\phi\) maps miles into kilometers. Then \(\tilde{\Cost}\) accepts miles, while
\(\Cost\) accepts kilometers. Then we have to map the initial starting point
\(\param_0\) of \(\Cost\) measured in kilometers into miles
\(\tilde{\param}_0\). \(\phi^{-1}\) is precisely this transformation from
kilometers into miles. A scale invariant optimizer should retain this relation, i.e.
no matter if the input is measured in miles or kilometers the same points are
selected.

\begin{proof}
	The following proof will be split into three parts. The first two parts of
	the proof will address a more general audience and ignore the mathematical
	subtleties we discussed in Section~\ref{sec: formal rfd}. In the third part
	we explain to the interested probabilists how to resolve these issues.

\begin{enumerate}[wide]
	\item \textbf{Invariance with regard to affine linear scaling}

	Let \(\tilde{\Cost}(\param):= a\Cost(\param) + b\)
	where \(a>0\) and \(b\in \real\) and assume
	\(\tilde{\param}_0 = \param_0\). With the induction start given, we
	only require the induction step to prove \(\tilde{\param}_n = \param_n\).

	For the induction step, we assume this equation holds up to \(n\).
	Since \(\phi(x) = ax +b\) is a measurable bijection, the sigma algebra\footnote{
		if you are unfamiliar with sigma algebras read them as ``information''.
	}
	generated by
	\[
		(\tilde{\Cost}(\param_n),\nabla\tilde{\Cost}(\param_n))
		= (\phi\circ\Cost(\param_n),a\nabla\Cost(\param_n))
	\]
	is therefore equal to the sigma algebra generated by
	\((\Cost(\param_n),\nabla\Cost(\param_n))\). This implies
	\begin{equation}
	\label{eq: output scale invariance}	
	\begin{aligned}
		\tilde{\param}_{n+1}
		&= \argmin_{\param}\E[\tilde{\Cost}(\param) \mid \tilde{\Cost}(\tilde{\param}_n), \nabla\tilde{\Cost}(\tilde{\param}_n)]
		\\
		\overset{\text{induction}}&= \argmin_\param \E[\tilde{\Cost}(\param) \mid \tilde{\Cost}(\param_n), \nabla\tilde{\Cost}(\param_n)]
		\\
		\overset{\text{sigma alg.}}&= \argmin_\param \E[\tilde{\Cost}(\param) \mid \Cost(\param_n), \nabla\Cost(\param_n)]
		\\
		\overset{\text{linearity}}&= \argmin_\param a\E[\Cost(\param) \mid \Cost(\param_n), \nabla\Cost(\param_n)] + b
		\\
		\overset{\text{monotonicity}}&= \argmin_\param \E[\Cost(\param) \mid \Cost(\param_n), \nabla\Cost(\param_n)]
		\\
		\overset{\text{def.}}&= \param_{n+1}
	\end{aligned}
	\end{equation}
	
	Where we have used the linearity of the conditional expectation and the
	strict monotonicity of \(\phi(x) = ax + b\).

	\item \textbf{Invariance with regard to certain input bijections}
		
	Let \(\phi\) be a differentiable bijection whose jacobian is invertible everywhere
	and assume \(\tilde{\Cost} := \Cost \circ \phi\). Since \(\phi\) is a
	bijection, \(\phi(M)\) is the domain of \(\Cost\) whenever \(M\) is the domain
	of \(\tilde{\Cost}\).

	For a starting point \(\param_0\in \phi(M)\) we now assume \(\tilde{\param}_0
	= \phi^{-1}(\param_0) \in M\) and are again going to prove the claim 
	\[
		\tilde{\param}_n = \phi^{-1}(\param_n).
	\]
	by induction. Assume that we have this claim up to \(n\). Then we have by
	induction
	\begin{equation}
		\tilde{\Cost}(\tilde{\param}_n)
		= \Cost \circ \phi(\phi^{-1}(\param_n))
		= \Cost(\param_n)
	\end{equation}
	and
	\[
		\nabla\tilde{\Cost}(\tilde{\param}_n)
		= \nabla_{\tilde{\param}_n} (\Cost \circ \phi(\tilde{\param}_n))
		= \phi'(\tilde{\param}_n)(\nabla\Cost)(\phi(\tilde{\param}_n))
		= \phi'(\tilde{\param}_n) \nabla\Cost(\param_n).
	\]
	Since \(\phi'(\tilde{\param}_n)\) is invertible by assumption, the sigma
	algebras generated by \((\tilde{\Cost}(\tilde{\param}_n),
	\nabla\tilde{\Cost}(\tilde{\param}_n))\) and \(\Cost(\param_n), \nabla\Cost(\param_n)\)
	are identical. But this results in the induction step
	\begin{equation}
	\label{eq: input scale invariance}	
	\begin{aligned}
		\tilde{\param}_{n+1}
		&= \argmin_{\param\in M}\E[\tilde{\Cost}(\param) \mid \tilde{\Cost}(\tilde{\param}_n), \nabla\tilde{\Cost}(\tilde{\param}_n)]
		\\
		\overset{\text{sigma alg.}}&=
		\argmin_{\param\in M}\E[\tilde{\Cost}(\param) \mid \Cost(\param_n), \nabla\Cost(\param_n)]
		\\
		\overset{\text{def.}}&= \argmin_{\param\in M}\E[\Cost\circ \phi(\param) \mid \Cost(\param_n), \nabla\Cost(\param_n)]
		\\
		&= \phi^{-1}\Bigl(
			\underbrace{
				\argmin_{\theta\in \phi(M)}
				\E[\Cost(\theta) \mid \Cost(\param_n), \nabla\Cost(\param_n)]
			}_{\overset{\text{def.}}= \param_{n+1}}
		\Bigr).
	\end{aligned}
	\end{equation}
	where we simply optimize over \(\theta = \phi(\param)\) instead of \(\param\) and
	correct the \(\argmin\) at the end.

	\item \textbf{Addressing the subtleties}

	In equation \eqref{eq: output scale invariance} we have really proven for deterministic \(\param\)
	\[
		\Phi_{\Pr_{\tilde{\Cost}}}(\param, \tilde{\Cost}(\param), \nabla\tilde{\Cost}(\param))
		= \Phi_{\Pr_\Cost}(\param, \Cost(\param), \nabla\Cost(\param)).
	\]
	But this implies with the induction assumption \(\Param_n = \tilde{\Param}_n\)
	\[
		\tilde{\Param}_{n+1} = 
		\Phi_{\Pr_{\tilde{\Cost}}}(
			\tilde{\Param}_n,
			\tilde{\Cost}(\tilde{\Param}_n),
			\nabla\tilde{\Cost}(\tilde{\Param}_n)
		)
		\overset{\text{ind.}}= \Phi_{\Pr_\Cost}(\Param_n, \Cost(\Param_n), \nabla\Cost(\Param_n))
		= \Param_{n+1}.
	\]
	Similarly we have proven in \eqref{eq: input scale invariance} that
	\[
		\Phi_{\Pr_{\tilde{\Cost}}}\bigl(
			\phi^{-1}(\param),
			\tilde{\Cost}(\phi^{-1}(\param)),
			\nabla\tilde{\Cost}(\phi^{-1}(\param))
		\bigr)
		= \phi^{-1}\bigl(
			\Phi_{\Pr_\Cost}(\param, \Cost(\param), \nabla\Cost(\param))
		\bigr).
	\]
	By the induction assumption \(\tilde{\Param}= \phi^{-1}(\Param_n)\), this
	implies
	\begin{align*}
		\tilde{\Param}_{n+1}
		&= \Phi_{\Pr_{\tilde{\Cost}}}(
				\tilde{\Param}_n,
				\tilde{\Cost}(\tilde{\Param}_n),
				\nabla\tilde{\Cost}(\tilde{\Param}_n))
		\\
		\overset{\text{ind.}}&=\Phi_{\Pr_{\tilde{\Cost}}}\bigl(
			\phi^{-1}(\Param_n),
			\tilde{\Cost}(\phi^{-1}(\Param_n)),
			\nabla\tilde{\Cost}(\phi^{-1}(\Param_n))
		\bigr)
		\\
		&= \phi^{-1}\bigl(
			\Phi_{\Pr_\Cost}(\Param_n, \Cost(\Param_n), \nabla\Cost(\Param_n))
		\bigr)
		\\
		&= \phi^{-1}(\Param_{n+1}).
		\qedhere
	\end{align*}
\end{enumerate}

\end{proof}

\subsection{Section~\ref{subsec: explicit rfd}: Relation to gradient descent}

\firstStochTaylor*
\begin{proof}
	\((\Cost(\param),\nabla\Cost(\param),\Cost(\param-\step))\) is a Gaussian vector for
	which the conditional distribution is well known. It is only necessary to calculate
	the covariance matrix. The key ingredient here is to observe that \(\Cost(\param), \partial_1\Cost(\param), \dots, \partial_\dims\Cost(\param)\)
	are all independent, trivializing matrix inversion.

	More formally, by Lemma~\ref{lem: cov of derivatives} we have
	\[
		\Cov\Bigl(\begin{pmatrix}
			\Cost(\param)\\
			\nabla\Cost(\param)
		\end{pmatrix}\Bigr)
		= \begin{pmatrix}
			\ikernel(0) & \\
			& -\ikernel'(0)\identity_{\dims\times\dims}
		\end{pmatrix}
	\]
	and
	\[
		\Cov\Bigl(\Cost(\param-\step), \begin{pmatrix}
			\Cost(\param)\\
			\nabla\Cost(\param)
		\end{pmatrix}\Bigr)
		= \begin{pmatrix}
				\ikernel(\frac{\|\step\|^2}{2})\\
				\ikernel'(\frac{\|\step\|^2}{2})\step
		\end{pmatrix}.
	\]
	By Theorem~\ref{thm: conditional gaussian distribution} we therefore know that
	\[
		\E[\Cost(\param-\step) \mid \Cost(\param),\nabla\Cost(\param)]	
		=  \mu +
		\begin{pmatrix}
				\ikernel(\frac{\|\step\|^2}{2})\\
				\ikernel'(\frac{\|\step\|^2}{2})\step
		\end{pmatrix}^T
		\begin{pmatrix}
			\ikernel(0) & \\
			& -\ikernel'(0)\identity_{\dims\times\dims}
		\end{pmatrix}^{-1}
		\begin{pmatrix}
			\Cost(\param) - \mu\\
			\nabla\Cost(\param)
		\end{pmatrix},
	\]
	which immediately yields the claim.
\end{proof}

\explicitRFD*
\begin{proof}
	The explicit version of RFD follows essentially by fixing the step size
	\(\stepsize = \|\step\|\) and optimizing over the direction first.
	With Lemma~\ref{lem: first stoch Taylor} we have
	\begin{align*}
		&\min_\step 
		\E[\Cost(\param-\step)\mid \Cost(\param),\nabla\Cost(\param)]
		\\
		&= \min_{\stepsize \ge 0} \min_{\step: \|\step\|=\stepsize}
		\mu + \frac{\ikernel\bigl(\frac{\stepsize^2}2\bigr)}{\ikernel(0)}
		(\Cost(\param)-\mu) - \frac{\ikernel'\bigl(\frac{\stepsize^2}2\bigr)}{\ikernel'(0)}
		\langle \step, \nabla\Cost(\param)\rangle
		\\
		&= \min_{\stepsize \ge 0}
		\mu + \frac{\ikernel\bigl(\frac{\stepsize^2}2\bigr)}{\ikernel(0)}
		(\Cost(\param)-\mu) - \frac{\ikernel'\bigl(\frac{\stepsize^2}2\bigr)}{\ikernel'(0)}
		\begin{cases}
			\displaystyle \max_{\step: \|\step\|=\stepsize}
			\langle \step, \nabla\Cost(\param)\rangle
			& \frac{\ikernel'(\frac{\stepsize^2}2)}{\ikernel'(0)} \ge 0
			\\
			\displaystyle \min_{\step: \|\step\|=\stepsize}
			\langle \step, \nabla\Cost(\param)\rangle
			& \frac{\ikernel'(\frac{\stepsize^2}2)}{\ikernel'(0)} < 0.
		\end{cases}
	\end{align*}
	By Lemma~\ref{lem: constrained maximiziation of scalar products} and
	Corollary~\ref{cor: constrained minimization of scalar products} the maximizing or
	minimizing step direction is then given by
	\[
		\step(\stepsize) = \pm \stepsize \frac{\nabla\Cost(\param)}{\|\nabla\Cost(\param)\|}.
	\]
	Where it is typically to be expected, that we have a positive sign. Since that depends
	on the covariance though, we avoid this problem with the following argument:
	Since \(\stepsize\) only appears as \(\stepsize^2\) in the remaining equation, we can
	optimize over \(\stepsize\in \real\) in the outer minimization instead of
	over \(\stepsize \ge 0\) to move the sign into the step size \(\stepsize\) and set
	without loss of generality
	\[
		\step(\stepsize) = \stepsize \frac{\nabla\Cost(\param)}{\|\nabla\Cost(\param)\|}.
	\]
	Since \(\langle \step(\stepsize), \nabla\Cost(\param)\rangle = \stepsize \|\nabla\Cost(\param)\|\)
	the remaining outer minimization problem over the step size is then given by
	\[
		\min_{\stepsize\in\real} \frac{\ikernel\bigl(\frac{\stepsize^2}2\bigr)}{\ikernel(0)}(\Cost(\param)-\mu)
		- \stepsize \frac{\ikernel'\bigl(\frac{\stepsize^2}2\bigr)}{\ikernel'(0)}\|\nabla\Cost(\param)\|,
	\]
	Its minimizer is by definition the RFD step size as given in the Theorem.
\end{proof}

\subsection{Section~\ref{subsec: rfd step sizes}: RFD-step sizes}

\begin{prop}[Tayloring the step size optimization problem]
	The second order Taylor approximation of the step size optimization problem
	\[
		q_\Theta(\stepsize)
		= -\frac{\ikernel\bigl(\tfrac{\stepsize^2}2\bigr)}{\ikernel(0)} 
		- \stepsize\frac{\ikernel'\bigl(\tfrac{\stepsize^2}2\bigr)}{\ikernel'(0)}\Theta
	\]
	around zero is given by
	\[
		T_2q_\Theta(\stepsize)
		= -1 
		- \stepsize \Theta
		+ \stepsize^2\frac{-\ikernel'(0)}{2\ikernel(0)}
		\quad
		\text{minimized by}
		\quad
		\hat{\stepsize}
		:= \argmin_\stepsize T_2q_\Theta(\stepsize)
		= \tfrac{\ikernel(0)}{-\ikernel'(0)}\Theta.
	\]
	Furthermore, the Taylor residual is bounded by
	\[
		\bigl|
			q(\stepsize)
			- T_2q(\stepsize)
		\bigr|
		\le
		\stepsize^3 c_0
		\bigl(\tfrac{\stepsize}{4} + \Theta\bigr)
	\]
	with  \(c_0=\frac12\max\{\sup_{\theta\in[0,1]}|\ikernel''(\theta)|, |\ikernel'(0)|\}(\tfrac1{\ikernel(0)} + \tfrac1{|\ikernel'(0)|})<\infty\).
\end{prop}

\begin{proof}
	Using the Taylor approximation with the mean value reminder for \(\ikernel\), we get
	\begin{align*}
		\ikernel\bigl(\tfrac{\stepsize^2}2\bigr)
		&= \ikernel(0) + \ikernel'(0)\tfrac{\stepsize^2}2
		+ \ikernel''(\theta_2)\frac{\bigl(\frac{\stepsize^2}2\bigr)^2}{2!}
		\\
		\ikernel'\bigl(\tfrac{\stepsize^2}2\bigr)
		&= \ikernel'(0) + \ikernel''(\theta_1)\tfrac{\stepsize^2}2
	\end{align*}
	for some \(\theta_1,\theta_2\in [0, \frac{\stepsize^2}2]\). This implies
	\[
		q(\stepsize)
		- \underbrace{\Bigl(
			-\bigl(1+\tfrac{\ikernel'(0)}{\ikernel(0)}\tfrac{\stepsize^2}2\bigr)
			- \stepsize \Theta 
		\Bigr)}_{=:T_2q_\Theta(\stepsize)}
		= -\frac{\ikernel''(\theta_2)}{\ikernel(0)}\frac{\stepsize^4}{2^3}
		- \frac{\ikernel''(\theta_1)}{\ikernel'(0)}\frac{\stepsize^3}2\Theta
	\]
	By the following error the optimistically defined \(T_2q_\Theta(\stepsize)\) is
	really the second Taylor approximation (which can be confirmed manually, but we
	deduce it by arguing that its residual is in \(\bigO(\stepsize^3)\)). More
	specifically,
	\[
		\bigl|
			q(\stepsize)
			- T_2q(\stepsize)
		\bigr|
		\le
		\stepsize^3
		\Bigl(
			\tfrac{\sup_{\theta\in [0, \frac{\stepsize^2}2]}|\ikernel''(\theta)|}{2\ikernel(0)}
			\frac{\stepsize}{4}
			+ \tfrac{\sup_{\theta\in [0, \frac{\stepsize^2}2]}|\ikernel''(\theta)|}{2|\ikernel'(0)|}
			\Theta
		\Bigr)
		\overset{\text{Lem.~\ref{lem: bound on the second derivative of the covariance}}}\le
		\stepsize^3 c_0
		\bigl(\tfrac{\stepsize}{4} + \Theta\bigr)
	\]
	It is easy to see for \(\Cost(\param)<\mu\) that \(T_2q(\stepsize)\) is a convex
	parabola due to \(\ikernel'(0)<0\). We thus have
	\[
		\hat{\stepsize}
		:= \argmin_\stepsize T_2q_\Theta(\stepsize)
		= \tfrac{\ikernel(0)}{-\ikernel'(0)}\Theta.
		\qedhere
	\]
\end{proof}

\begin{theorem}[Details of Proposition~\ref{prop: a-rfd well defined}]
	Let \(\Cost\sim\normal(\mu, \ikernel)\) and assume there exists \(\stepsize_0>0\) such
	that the correlation for larger distances \(\stepsize\ge \stepsize_0\) are
	bounded smaller than \(1\), i.e. \(\frac{\ikernel(\stepsize^2/2)}{\ikernel(0)} < \rho \in (0,1)\).
	Then there exists \(K, \Theta_0>0\) such that for all \(\Theta < \Theta_0\)
	\[
		1-K\Theta \le \frac{\stepsize^*(\Theta)}{\hat{\stepsize}(\Theta)}\le 1+ K\Theta.
	\]
	In particular we have \(\stepsize^*(\Theta)\sim \hat{\stepsize}(\Theta)\) as \(\Theta\to 0\) or
	equivalently as \(\hat{\stepsize}\to 0\).
	
\end{theorem}
\begin{proof}
	This follows immediately from Lemma~\ref{lem: small step sizes},
	Lemma~\ref{lem: medium step sizes} and Lemma~\ref{lem: large step sizes}.
\end{proof}

\convergence*
\begin{proof}
	Assuming RFD converges, its step sizes \(\stepsize^*\) converge to zero. But this implies
	\(\Theta\to 0\) by assumption, i.e.	
	\[
		\Theta = \frac{\|\nabla\Cost(\param)\|}{\mu - \Cost(\param)} \to 0
	\]
	Since \(\Cost(\param)\) is bounded, this implies \(\|\nabla\Cost(\param)\|\to 0\)
	and by continuity the of the gradient, it is zero in its limit. Thus we converge
	to a stationary point. The asymptotic equality follows by Lemma~\ref{lem: small step sizes}
	and Lemma~\ref{lem: medium step sizes}, as we know \(\stepsize^*\) converges so we do
	not require the assumptions of Lemma~\ref{lem: large step sizes}.
\end{proof}

\subsubsection{Locating the Minimizer}

In the following we want to rule out locations for the RFD step size \(\stepsize^*\) by
proving \(q_\Theta(\stepsize) > q_\Theta(\hat{\stepsize})\) for a wide range of \(\stepsize\).
For this endeavour the relative position of the step size \(\stepsize\) relative to \(\hat{\stepsize}\)
is a useful re-parametrization \[
	\stepsize := \stepsize(\lambda) = \lambda \hat{\stepsize}.
\]
Due to \(\hat{\stepsize}=\frac{\ikernel(0)}{-\ikernel'(0)}\Theta\)
we obtain
\[
	T_2q_\Theta(\stepsize)
	= -1 - \stepsize \Theta + \tfrac{\stepsize^2}2\tfrac{-\ikernel'(0)}{\ikernel(0)}
	= -1 + \lambda(\tfrac{\lambda}2 - 1)\hat{\stepsize}\Theta
\]
On the other hand we have for the bound
\begin{align*}
	|q_\Theta(\stepsize)- T_2q_\Theta(\stepsize)|
	\le 
	\lambda^3\hat{\stepsize}^3 c_0\Bigl(\lambda\tfrac{\ikernel(0)}{4|\ikernel'(0)|} + 1\Bigr)\Theta
\end{align*}
Since \(\hat{\stepsize} = \stepsize(1)\) we thus obtain
\begin{align}
	\nonumber
	\frac{q_\Theta(\stepsize) - q_\Theta(\hat{\stepsize})}{\hat{\stepsize}\Theta}
	&\ge \frac{
		\magenta{T_2q_\Theta(\stepsize)}
		- |q_\Theta(\stepsize) - T_2q_\Theta(\stepsize)| 
		- \teal{T_2q_\Theta(\hat{\stepsize})}
		- \blue{|q_\Theta(\hat{\stepsize}) - T_2q_\Theta(\hat{\stepsize})|}
	}{\hat{\stepsize}\Theta}
	\\
	\nonumber
	&\ge
	\underbrace{\bigl(
		\magenta{\lambda(\tfrac{\lambda}2 - 1)} - (\teal{-\tfrac12})
	\bigr)}_{
		= \tfrac12 - \lambda + \tfrac{\lambda^2}2
	}
	- \hat{\stepsize}^2 c_0\Bigl[
		\lambda^3\Bigl(\lambda\tfrac{\ikernel(0)}{4|\ikernel'(0)|} + 1\Bigr)
		+ \blue{\Bigl(\tfrac{\ikernel(0)}{4|\ikernel'(0)|} + 1\Bigr)}
	\Bigr]
	\\
	\label{eq: minimizer rule-out equation}
	&=
	\tfrac12 (1-\lambda)^2
	- \hat{\stepsize}^2 c_0\Bigl[
		\lambda^3\Bigl(\lambda\tfrac{\ikernel(0)}{4|\ikernel'(0)|} + 1\Bigr)
		+ \Bigl(\tfrac{\ikernel(0)}{4|\ikernel'(0)|} + 1\Bigr)
	\Bigr].
\end{align}
This equation will be the basis of a number of lemmas ruling out various step sizes as minimizers.

\begin{lemma}[Ruling out small step sizes]
	\label{lem: small step sizes}
	If the step size is (much) smaller than the asymptotic step size
	\(\hat{\stepsize}=\hat{\stepsize}(\Theta)\), then it can not be a minimizer. More specifically
	\[
		\frac{\stepsize}{\hat{\stepsize}} \in [0, 1-c_1\Theta) \implies q_\Theta(\stepsize) > q_\Theta(\hat{\stepsize})
	\]	
	where \(c_1 := 2 \tfrac{\ikernel(0)}{|\ikernel'(0)|} \sqrt{c_0\bigl(\tfrac{\ikernel(0)}{4|\ikernel'(0)|} + 1\bigr)}<\infty\).
\end{lemma}
\begin{proof}
	Here we consider the case \(\stepsize\le \hat{\stepsize}\), i.e. \(\lambda \in [0,1]\).
	By \eqref{eq: minimizer rule-out equation} we have
	\begin{align*}
		\frac{q_\Theta(\stepsize) - q_\Theta(\hat{\stepsize})}{\hat{\stepsize}\Theta}
		&\ge
		\tfrac12 (1-\lambda)^2
		- \hat{\stepsize}^2 c_0\Bigl[
			\lambda^3\Bigl(\lambda\tfrac{\ikernel(0)}{4|\ikernel'(0)|} + 1\Bigr)
			+ \Bigl(\tfrac{\ikernel(0)}{4|\ikernel'(0)|} + 1\Bigr)
		\Bigr]
		\\
		&\ge
		\tfrac12 (1-\lambda)^2
		- 2\hat{\stepsize}^2 c_0 
			\Bigl(\tfrac{\ikernel(0)}{4|\ikernel'(0)|} + 1\Bigr)
		\\
		\overset{!}&> 0
	\end{align*}
	for which
	\[
		(1-\lambda)^2
		> 4 \hat{\stepsize}^2 c_0\Bigl(\tfrac{\ikernel(0)}{4|\ikernel'(0)|} + 1\Bigr)
	\]
	is sufficient or equivalently
	\[
		\lambda < 1-2 \hat{\stepsize} \sqrt{c_0\Bigl(\tfrac{\ikernel(0)}{4|\ikernel'(0)|} + 1\Bigr)}
		= 1- \Theta \underbrace{
			2 \tfrac{\ikernel(0)}{|\ikernel'(0)|} \sqrt{c_0\Bigl(\tfrac{\ikernel(0)}{4|\ikernel'(0)|} + 1\Bigr)}
		}_{=:c_1}
	\]
	So for \(\lambda \in [0, 1-\Theta c_1)\) we have \(q_\Theta(\stepsize) > q_\Theta(\hat{\stepsize})\).
\end{proof}

\begin{lemma}[Ruling out medium sized step sizes as minimizer]
	\label{lem: medium step sizes}
	For \(c_2 = 2c_1\) and \(\Theta \le \Theta_0 := \frac1{5c_1}\), we have
	\[
		\frac{\stepsize}{\hat{\stepsize}} \in \bigl(1+ c_2\Theta, \tfrac1{c_2\Theta}\bigr)
		\implies q_\Theta(\stepsize) > q_\Theta(\hat{\stepsize})
	\]
\end{lemma}
\begin{proof}
	Here we consider the case \(\lambda \ge 1\), i.e. \(\stepsize>\hat{\stepsize}\).
	Again starting with \eqref{eq: minimizer rule-out equation} we get
	\begin{align*}
		\frac{q(\stepsize) - q(\hat{\stepsize})}{\hat{\stepsize}\Theta}
		&\ge
		\tfrac12 (1-\lambda)^2
		- \hat{\stepsize}^2 c_0\Bigl[
			\lambda^3\Bigl(\lambda\tfrac{\ikernel(0)}{4|\ikernel'(0)|} + 1\Bigr)
			+ \Bigl(\tfrac{\ikernel(0)}{4|\ikernel'(0)|} + 1\Bigr)
		\Bigr]
		\\
		&\ge
		\tfrac12 (\lambda-1)^2
		- 2\lambda^4\hat{\stepsize}^2 c_0\Bigl(\tfrac{\ikernel(0)}{4|\ikernel'(0)|} + 1\Bigr)
		\\
		\overset{!}&> 0,
	\end{align*}
	for which
	\[
		\lambda -1 > 2 \lambda^2\hat{\stepsize}\sqrt{c_0\Bigl(\tfrac{\ikernel(0)}{4|\ikernel'(0)|} + 1\Bigr)}
		= c_1\Theta \lambda^2
	\]
	or equivalently
	\[
		\lambda -1 - c_1\Theta\lambda^2
		> 0
	\]
	is sufficient. Note that this is a concave parabola in \(\lambda\). So it is positive
	between its zeros which are characterized by
	\[
		c_1\Theta\lambda^2 - \lambda + 1 = 0.
	\]
	They are thus given by
	\[
		\lambda_{1/2} = \frac{1 \pm \sqrt{1 - 4c_1\Theta}}{2c_1\Theta}.
	\]
	So whenever \(\lambda \in (\lambda_1, \lambda_2)\) we have that
	\(q_\Theta(\stepsize) > q_\Theta(\hat{\stepsize})\). In particular for \(4c_1\Theta\le 1\)
	or equivalently \(\Theta \le \frac1{4c_1}\) we have
	\[
		\lambda_2 \ge \frac1{2c_1\Theta}	 = \frac1{c_2\Theta}
	\]
	To get a bound on \(\lambda_1\) note that the original equation was essentially
	\[
		\lambda \ge 1+c_1\Theta\lambda^2
	\]
	with equality for \(\lambda=\lambda_1\), if \(\Theta\) is reduced, the inequality remains,
	which implies that \(\lambda_1\) is decreasing with \(\Theta\). So
	assuming the inequality is satisfied for a particular \(\lambda\) e.g.
	\(\lambda=\sqrt{2}\) which requires
	\[
		\sqrt{2} \ge 1+ 2c_1\Theta  \iff \Theta \le \tfrac{\sqrt{2} -1}{2c_1},
	\]
	then we know that \(\lambda_1\le \sqrt{2}\) for all smaller \(\Theta\). This
	implies for \(\Theta \le \Theta_0 = \frac1{5c_1} \le \frac{\sqrt{2} -1}{2c_1}\)
	\[
		\lambda_1 = 1+ c_1\Theta\lambda_1^2 \le 1+ \underbrace{2c_1}_{c_2}\Theta.
		\qedhere
	\]
\end{proof}

\begin{lemma}[Ruling out large step sizes as minimizer]
	\label{lem: large step sizes}
	If there exists step size \(\stepsize_0>0\) such that the correlation is
	bounded by some \(\rho<1\), i.e.
	\[
		\frac{\ikernel\bigl(\frac{\stepsize^2}2\bigr)}{\ikernel(0)} \le \rho\in (0,1),
	\]
	for larger step sizes \(\stepsize\ge \stepsize_0\), then there exist
	\(\Theta_0>0\) such that for all \(\Theta<\Theta_0\)
	\[
		\frac{\stepsize}{\hat{\stepsize}} \in
		\bigl(1+c_2\Theta,\infty\bigr)
		\implies q(\stepsize) > q(\hat{\stepsize}),
	\]
	where \(c_2\) is the constant from Lemma~\ref{lem: medium step sizes}.
\end{lemma}
\begin{proof}
	The upper bound \(\frac{1}{c_2\Theta}\) in Lemma~\ref{lem: medium step sizes} is
	only due to the loss of precision of the Taylor approximation. To remove it,
	we take a closer look at the actual \(q_\Theta\) itself. We have the following
	bound for our asymptotic minimum
	\begin{align*}
		\frac{q_\Theta(\hat{\stepsize})}{\Theta}
		&\le \frac{T_2 q_\Theta(\hat{\stepsize}) + |q_\Theta(\hat{\stepsize}) - T_2q_\Theta(\hat{\stepsize})|}{\Theta}
		= -\frac1\Theta - \frac12 \hat{\stepsize}
		+ \hat{\stepsize}^3 \underbrace{c_0\Bigl(\tfrac{\ikernel(0)}{4|\ikernel'(0)|} + 1\Bigr)}_{=:c_3}
		\\
		&\le - \frac1\Theta + \hat{\stepsize}^3 c_3
	\end{align*}
	Which means we have for
	\begin{align*}
		\frac{q_\Theta(\stepsize) - q_\Theta(\hat{\stepsize})}{\Theta}
		&\ge \Bigl(1- \frac{\ikernel\bigl(\tfrac{\stepsize^2}2\bigr)}{\ikernel(0)}\Bigr)\frac1\Theta
		- \stepsize\frac{\ikernel'\bigl(\tfrac{\stepsize^2}2\bigr)}{\ikernel'(0)}
		- \hat{\stepsize}^3 c_3
		\\
		\overset{\text{Lemma~\ref{lem: bound on first derivative of covariance}}}&\ge
		\Bigl(1- \frac{\ikernel\bigl(\tfrac{\stepsize^2}2\bigr)}{\ikernel(0)}\Bigr)\frac1\Theta
		- \tfrac{\sqrt{\ikernel(0)}}{\sqrt{-\ikernel'(0)}}
		- \hat{\stepsize}^3 c_3
		\\
		&\ge
		\bigl(1- M\bigr)\frac1\Theta
		- \tfrac{\sqrt{\ikernel(0)}}{\sqrt{-\ikernel'(0)}}
		- \hat{\stepsize}^3 c_3
		\\
		\overset{!}&> 0,
	\end{align*}
	where we use the assumption that there exists \(\rho\in (0,1)\) such that
	\(\rho \ge \frac{\ikernel\bigl(\tfrac{\stepsize^2}2\bigr)}{\ikernel(0)}\) for
	all \(\stepsize\ge \stepsize_0\) and the fact that we only need
	to consider \(\stepsize \ge \frac{1}{c_2\Theta}\) (due to Lemma~\ref{lem: medium
	step sizes}) which allows a translation of \(\stepsize_0\) into some maximal
	\(\Theta_0\).  Note that \(\hat{\stepsize}\sim \Theta\) vanishes as \(\Theta\to 0\),
	so eventually the term \((1-M)\frac1\Theta\) dominates. Selecting
	\(\Theta_0\) small small enough is thus sufficient to cover everything that is not
	already covered by Lemma~\ref{lem: medium step sizes}.
\end{proof}

\subsubsection{Technical bounds}

\begin{lemma}[Bound on the first derivative of the covariance]
	\label{lem: bound on first derivative of covariance}	
	\[
		\sup_{\stepsize\ge 0}|\ikernel'\bigl(\tfrac{\stepsize^2}2\bigr)\stepsize|
		\le \sqrt{-\ikernel'(0)\ikernel(0)}
	\]
\end{lemma}
\begin{proof}
	Since we have	
	\[
		\Cov(D_v \Cost(x), \Cost(y))
		= \ikernel'\bigl(\tfrac{\|x-y\|^2}2\bigr)\langle x-y, v\rangle
	\]
	we have for a standardized vector \(\|v\|=1\) and \(x-y = \stepsize v\)
	by Cauchy-Schwarz
	\[
		|\ikernel'\bigl(\tfrac{\stepsize^2}2\bigr)\stepsize|
		= |\Cov(D_v \Cost(x), \Cost(y))|
		\overset{\text{C.S.}}\le \sqrt{\Var(D_v \Cost(x)) \Var(\Cost(y))}
		= \sqrt{-\ikernel'(0)\ikernel(0)}.
	\]
	As the bound is independent of \(\stepsize\) this yields the claim.
\end{proof}

\begin{lemma}[Bound on the second derivative of the covariance]
	\label{lem: bound on the second derivative of the covariance}
	\[
		\sup_{\theta\ge 0}|\ikernel''(\theta)|
		\le \max\Bigl\{\sup_{\theta\in[0,1]}|\ikernel''(\theta)|, |\ikernel'(0)|\Bigr\}
	\]
\end{lemma}
\begin{proof}
	Note that
	\[
		\Cov(D_v \Cost(x), D_w\Cost(y))
		= -\ikernel''\bigl(\tfrac{\|x-y\|^2}{2}\bigr)\langle x-y, v\rangle\langle x-y, w\rangle
		- \ikernel'\bigl(\tfrac{\|x-y\|^2}{2}\bigr)\langle v, w\rangle
	\]
	Selecting \(v,w\) as orthonormal vectors (e.g. \(v=e_1, w=e_2\)) and \(x-y := \stepsize (v+w)\)
	for some \(\stepsize > 0\) results in \(\|x-y\|^2 = 2\stepsize^2\) and thus
	by the Cauchy-Schwarz inequality
	\[
		\bigl|-\ikernel''(\stepsize^2)\stepsize^2\bigr|
		= \bigl| \Cov(D_v \Cost(x), D_w\Cost(y))\bigr|
		\overset{\text{C.S.}}\le\sqrt{\Var(D_v\Cost(x))\Var(D_w\Cost(y))}
		= \sqrt{(-\ikernel'(0))^2}
	\]
	This implies the claim.
\end{proof}

\subsection{Section~\ref{sec: stochastic loss and covariance estimation}: Stochastic loss}

\begin{lemma}
    \label{lem: stoch approx errors cond. independent}
    The stochastic approximation errors
    \[
        \epsilon_i(\param) := \loss_i(\param) - \Cost(\param)
    \]
    are identically distributed, centered random functions, which are
    independent conditional on \(\rf\). In particular,
    \[
        \E[\epsilon_i(\param) \epsilon_j(\tilde{\param})]
        = \E[\epsilon_i(\param) \epsilon_j(\tilde{\param})\mid \rf] = 0
        \quad \forall j\neq i.
    \]
\end{lemma}
\begin{proof}
	The \(\epsilon_i\) are independent random functions conditional on \(\rf\), since
	for any \(n\in \nat\), any bounded measurable functions \(h\) and \(g\)
	\begin{align*}
		&\E\Bigl[
			h\bigl(\epsilon_i(\param_1),\dots,\epsilon_i(\param_n)\bigr)
			g\bigl(\epsilon_j(\param_1),\dots,\epsilon_j(\param_n)\bigr) \mid \rf
		\Bigr]
		\\
		&= 
		\Bigl[
			h\bigl(\epsilon_i(\param_1),\dots,\epsilon_i(\param_n)\bigr)
		\underbrace{
				\E\Bigl[g\bigl(\epsilon_j(\param_1),\dots,\epsilon_j(\param_n)\bigr) \mid \rf, X_i, \noise_i\Bigr]
			}_{\overset{(*)}=\E\bigl[g(\epsilon_j(\param_1),\dots,\epsilon_j(\param_n)) \mid \rf\bigr]}
			\Bigm| \rf
		\Bigr]
		\\
		&= \E\Bigl[
			h\bigl(\epsilon_i(\param_1),\dots,\epsilon_i(\param_n)\bigr) \mid\rf
		\Bigr]
		\E\Bigl[
			g\bigl(\epsilon_j(\param_1),\dots,\epsilon_j(\param_n)\bigr) \mid \rf
		\Bigr],
	\end{align*}
	where \((*)\) uses the fact that \(\epsilon_j\) does not depend on the
	independent \(X_i, \noise_i\). Since almost by definition
	\[
		\E[\epsilon_i \mid \rf] = \E[\loss(\cdot, X_i, Y_i) \mid \rf] - \Cost(\cdot) = 0,
	\]
	the stochastic approximation errors are thus uncorrelated
	\[
		\E[\epsilon_i \epsilon_j]
		= \E\Bigl[ \E[\epsilon_i \epsilon_j\mid \rf]\Bigr]
		= \E\Bigl[ \E[\epsilon_i \mid\rf] \E[\epsilon_j \mid \rf]\Bigr]
		= 0.
	\]
\end{proof}

\srfd*
\begin{proof}
	Since \(\epsilon_i\) are conditionally independent between each other and to \(\Cost\),
	as entire functions, the same holds true for \(\nabla\epsilon_i\). As all the mixed
	covariances disappear we have
	\begin{align*}
		\Cov\Bigl(\begin{pmatrix}
			\Loss_\batchsize(\param)\\
			\nabla\Loss_\batchsize(\param)
		\end{pmatrix}\Bigr)	
		&= 
		\Cov\Bigl(\begin{pmatrix}
			\Cost(\param)\\
			\nabla\Cost(\param)
		\end{pmatrix}\Bigr)	
		+ \frac1{\batchsize^2}
		\sum_{i=1}^\batchsize
		\Cov\Bigl(\begin{pmatrix}
			\epsilon_i(\param)\\
			\nabla\epsilon_i(\param)
		\end{pmatrix}\Bigr)	
		\\
		&= 
		\begin{pmatrix}
			\ikernel(0) &
			\\
			& -\ikernel'(0)\identity_{\dims\times\dims}
		\end{pmatrix}
		+ \frac1{\batchsize^2}
		\sum_{i=1}^\batchsize
		\begin{pmatrix}
			\ikernel_\epsilon(0) &
			\\
			& -\ikernel_\epsilon'(0)\identity_{\dims\times\dims}
		\end{pmatrix}
		\\
		&= 
		\begin{pmatrix}
			\ikernel(0) + \frac1\batchsize\ikernel_\epsilon(0) &
			\\
			& -\Bigl(\ikernel'(0) + \frac1\batchsize \ikernel_\epsilon'(0)\Bigr)
			\identity_{\dims\times\dims}.
		\end{pmatrix}
	\end{align*}
	by Lemma~\ref{lem: cov of derivatives}. If you want to break up the first
	step we recommend considering individual entries of the covariance matrix to
	convince yourself that all the mixed covariances disappear. Together with the
	fact
	\begin{align*}
		&\Cov\Bigl(
			\Cost(\param-\step),
			\begin{pmatrix}
				\Loss_\batchsize(\param)\\
				\nabla\Loss_\batchsize(\param)
			\end{pmatrix}
		\Bigr)	
		\\
		&= 
		\Cov\Bigl(
			\Cost(\param-\step),
			\begin{pmatrix}
				\Cost(\param)\\
				\nabla\Cost(\param)
			\end{pmatrix}
		\Bigr)
		+ \frac1{\batchsize^2}\sum_{i=1}^\batchsize
		\underbrace{
		\Cov\Bigl(
			\Cost(\param-\step),
			\begin{pmatrix}
				\epsilon_i(\param)\\
				\nabla\epsilon_i(\param)
			\end{pmatrix}
		\Bigr)
		}_{=0}
		\\
		&= \begin{pmatrix}
				\ikernel(\frac{\|\step\|^2}{2})\\
				\ikernel'(\frac{\|\step\|^2}{2})\step
		\end{pmatrix}.
	\end{align*}
	The rest is analogous to Lemma~\ref{lem: first stoch Taylor} and
	Theorem~\ref{thm: explicit rfd}, so we only sketch the remaining steps.
	
	Applying Theorem~\ref{thm: conditional gaussian distribution} as in Lemma~\ref{lem:
	first stoch Taylor} we obtain a
	stochastic version of the stochastic Taylor approximation
	(``stochastic\({}^2\) Taylor approximation'' perhaps?)
	\[
		\E[\Cost(\param-\step) \mid \Loss_\batchsize(\param), \nabla\Loss_\batchsize(\param)]	
		= \mu
		+ \frac{\ikernel\bigl(\frac{\|\step\|^2}2\bigr)}{\ikernel(0) + \frac1\batchsize\ikernel_\epsilon(0)}(\Loss_\batchsize(\param) - \mu)
		-  \frac{\ikernel'\bigl(\frac{\|\step\|^2}2\bigr)}{\ikernel'(0) + \frac1\batchsize\ikernel_\epsilon'(0)}\langle \step, \Loss_\batchsize(\param)\rangle.
	\]
	Minimizing this subject to a constant step size as in Theorem~\ref{thm: explicit rfd}
	results in 
	\begin{align*}
		\stepsize^*
		&= \argmin_{\stepsize\in \real}
		\frac{\ikernel\bigl(\frac{\|\step\|^2}2\bigr)}{\ikernel(0) + \frac1\batchsize\ikernel_\epsilon(0)}
		(\Loss_\batchsize(\param) - \mu)
		-  \stepsize\frac{\ikernel'\bigl(\frac{\|\step\|^2}2\bigr)}{\ikernel'(0) + \frac1\batchsize\ikernel_\epsilon'(0)}
		\|\Loss_\batchsize(\param)\|
		\\
		&= \argmin_{\stepsize\in \real}
		- \frac{\ikernel\bigl(\frac{\|\step\|^2}2\bigr)}{\ikernel(0)}
		-  \stepsize\frac{\ikernel'\bigl(\frac{\|\step\|^2}2\bigr)}{\ikernel'(0) + \frac1\batchsize\ikernel_\epsilon'(0)}
		\frac{\ikernel(0)}{\ikernel(0) + \frac1\batchsize\ikernel_\epsilon(0)}
		\frac{\|\Loss_\batchsize(\param)\|}{\mu -\Loss_\batchsize(\param)},
	\end{align*}
	where we divided the term by \(\frac{\ikernel(0)}{\ikernel(0) +
	\frac1\batchsize\ikernel_\epsilon(0)}\frac{1}{\mu
	-\Loss_\batchsize(\param)} \ge 0\)
	to obtain the last equation. The claim follows by definition of \(\stepsize^*(\Theta)\) and
	our redefinition of \(\Theta\).
\end{proof}

	\section{Extensions}
\label{sec: extensions}

In this section we present a few possible extensions to Theorem~\ref{thm:
explicit rfd}, which are all composable, i.e. it is possible to combine these
extensions without any major problems
(including S-RFD, i.e. Extension~\ref{ext: s-rfd}).

\subsection{Geometric anisotropy/Adaptive step sizes}
\label{sec: geometric anisotropy}

In this section, we discuss the generalization of isotropy to ``geometric anisotropies''
\parencite[17]{steinInterpolationSpatialData1999}, which provide good insights into
the inner workings of adaptive learning rates (e.g. AdaGrad \parencite{duchiAdaptiveSubgradientMethods2011} and
Adam \parencite{kingmaAdamMethodStochastic2015}).
\begin{definition}[Geometric Anisotropy]
	We say a random function \(\Cost\) exhibits a ``geometric anisotropy'', if there exists
	an invertible matrix \(A\) such that \(\Cost(x) = \rg(Ax)\) for some isotropic random
	function \(\rg\).
\end{definition}

This implies that the expectation of \(\Cost\) is still constant (\(\E[\Cost(x)] =
\E[\rg(Ax)] = \mu\)) and the covariance function of \(\Cost\) is given by
\begin{equation}
	\label{eq: geometric anisotropy characterization}
	\Cov(\Cost(x), \Cost(y))
	= \Cov(\rg(Ax), \rg(Ay))
	= \ikernel\Bigl(\frac{\|A(x-y)\|^2}2\Bigr)
	= \ikernel\Bigl(\frac{\|x-y\|_{A^TA}^2}2\Bigr)
\end{equation}
where \(\|\cdot\|_\Sigma\) is the norm induced by \(\langle x,y\rangle_\Sigma :=
\langle x, \Sigma y\rangle\) for some strictly positive definite matrix \(\Sigma= A^T A\).
Here \eqref{eq: geometric anisotropy characterization} characterizes the set of
random functions with a geometric anisotropy in the Gaussian case, because for
an \(\Cost\) with such a covariance we can always obtain an isotropic \(\rg\) by
\(\rg(x):= \Cost(A^{-1}x)\). This is the whitening transformation we suggest looking
for in order to ensure isotropy in the context of scale invariance (Section~\ref{sec: rfd}).

An important observation is, that Theorem~\ref{thm: characterization of weak input
invariances} implies that \(\Cost\) is still stationary, so the distribution of
\(\Cost\) is still invariant to translations. If stationarity is a problem, this is
therefore not the solution. But geometric anisotropies are a beautiful model to
explain preconditioning and adaptive step sizes. For this, we first determine the RFD steps.

\begin{restatable}[RFD steps under geometric anisotropy]{extension}{geometricAnisotropyRFD}
	\label{ext: geometric anisotropy rfd steps}
	Let \(\Cost\) be a Gaussian random function which exhibits a ``geometric
	anisotropy'' \(A\) and is based on an isotropic random function
	\(\rg\sim\normal(\mu, \ikernel)\). Then the RFD steps are given by
	\[
		\stepsize^* \frac{\Sigma^{-1}\nabla\Cost(\param)}{\|\Sigma^{-1}\nabla\Cost(\param)\|_\Sigma}
		= \argmin_{\step} \E[\Cost(\param- \step) \mid \Cost(\param), \nabla\Cost(\param)]
	\]
	with
	\[
		\stepsize^* = \argmin_\stepsize q_\Theta(\stepsize)	
		\quad\text{where}\quad
		\Theta = \frac{\|\Sigma^{-1}\nabla\Cost(\param)\|_\Sigma}{\mu - \Cost(\param)}.
	\]	
\end{restatable}
\begin{proof}[Proofsketch]
	There are two ways to see this. Either we apply scale invariance
	(Advantage~\ref{advant: scale invariance}) directly to translate the steps on \(\rg\) into
	steps on \(\Cost\). Alternatively one can manually retrace the steps of the proof.
	Details in Subsection~\ref{sec: proof of rfd steps under geometric anisotropy}
\end{proof}

The step direction is therefore
\[
	\Sigma^{-1}	\nabla\Cost(x)
\]
and \(\Sigma^{-1}\) acts as a preconditioner. So how would one obtain \(\Sigma\)? As it
turns out the following holds true (by Lemma~\ref{lem: cov of derivatives})
\[
	\E[\nabla\Cost(\param)\nabla\Cost(\param)^T]
	= A^T \E[\nabla\rg(\param)\nabla\rg(\param)^T] A
	= A^T (-\ikernel'(0)\identity)  A
	= -\ikernel'(0)\Sigma
\]
In their proposal of the first ``adaptive'' method, AdaGrad,
\textcite{duchiAdaptiveSubgradientMethods2011} suggest to use the matrix
\[
	G_t = \sum_{k=1}^t \nabla\Cost(\param_k)\nabla\Cost(\param_k)^T,
\]
which is basically already looking like an estimation method of \(\Sigma\).
They then restrict themselves to \(\diag(G_t)\) due to the computational costs of a
full matrix inversion. This results in entry-wise (``adaptive'') learning rates.
Later adaptive methods like RMSProp \parencite{hintonNeuralNetworksMachine2012},
AdaDelta
\parencite{zeilerADADELTAAdaptiveLearning2012} and in particular Adam
\parencite{kingmaAdamMethodStochastic2015} replace this sum with an exponential mean
estimate, i.e. in the case of Adam the decay rate \(\beta_2\) is used to get an
exponential moving average
\[
	v_t
	= \beta_2 v_{t-1} + (1-\beta_2)\diag(\nabla\Cost(\param_t)\nabla\Cost(\param_t)^T)
	= \beta_2 v_{t-1} + (1-\beta_2)(\nabla\Cost(\param_t))^2.
\]
They then take the expectation
\[
	\E[v_t]
	= \E\Bigl[(1-\beta_2)\sum_{k=1}^t \beta_2^{t-k}\nabla\Cost(\param_k)^2\Bigr]
	= \E\Bigl[(1-\beta_2)\sum_{k=1}^t \beta_2^{t-k}\nabla\Cost(\param_k)^2\Bigr]
	= (1-\beta_2^t)\underbrace{\E[\nabla\Cost(x_t)^2]}_{\propto \diag(\Sigma)}
\]
So \(\hat{v}_t = v_t/(1-\beta_2^t)\) in the Adam optimizer is essentially an
estimator for \(\diag(\Sigma)\). It is noteworthy, that
\textcite{kingmaAdamMethodStochastic2015} already used the expectation symbol.
This is despite the fact, that they did not yet model the optimization objective
\(\Cost\) as a random function.

We can not yet explain why they then use the square root of their estimate
\(\diag(\Sigma)^{-1/2}\) instead of \(\diag(\Sigma)^{-1}\) itself. This might have something
to do with the fact that the estimation of \(G_t\) happens online and the
\(\Cost(\param_k)\) are therefore highly correlated. Another reason might be
that the inverse of an estimator has different properties than the estimator
itself. Finally, the fact that only the diagonal is used might also be the reason, if
the preconditioner \(\diag(\Sigma)^{-1/2}\) is simply better when we restrict ourselves
to diagonal matrices.

\subsubsection{Proof of Extension~\ref{ext: geometric anisotropy rfd steps}}
\label{sec: proof of rfd steps under geometric anisotropy}
	Since the application of scale invariance provides no intuition, we provide a
	proof which retraces some of the steps of the original proof.

	Recall, that for an isotropic random function \(\rg\) we have the stochastic
	Taylor approximation
	\[
		\E[\rg(\param - \step) \mid \rg(x), \nabla\rg(x)]
		= \mu + \frac{\ikernel\bigl(\frac{\|\step\|^2}2\bigr)}{\ikernel(0)}(\rg(\param)-\mu)
		+ \frac{\ikernel'\bigl(\frac{\|\step\|^2}2\bigr)}{\ikernel'(0)}\langle \step, \nabla \rg(\param)\rangle
	\]
	This implies for a random function with geometric anisotropy \(\Cost(\param) =
	\rg(A\param)\) that
	\begin{align*}
		\E[\Cost(\param - \step) \mid \Cost(\param), \nabla\Cost(\param)]
		&=\E[\rg(A(\param - \step)) \mid \rg(A\param), \nabla\rg(A\param)]
		\\
		&= \mu + \frac{\ikernel\bigl(\frac{\|A\step\|^2}2\bigr)}{\ikernel(0)}(\rg(A\param)-\mu)
		- \frac{\ikernel'\bigl(\frac{\|A\step\|^2}2\bigr)}{\ikernel'(0)}\langle A\step, \nabla \rg(A\param)\rangle
		\\
		&= \mu + \frac{\ikernel\bigl(\frac{\|\step\|_{\Sigma}^2}2\bigr)}{\ikernel(0)}(\Cost(\param)-\mu)
		- \frac{\ikernel'\bigl(\frac{\|\step\|_{\Sigma}^2}2\bigr)}{\ikernel'(0)}\langle \step, \underbrace{A^T\nabla \rg(A\param)}_{=\nabla\Cost(\param)}\rangle
	\end{align*}
	with \(\Sigma:= A^T A\). As in the original proof, we now optimize over the
	direction first, while keeping the step size constant, although we now fix the
	step size with regard to the norm \(\|\cdot\|_\Sigma\) (which basically means that
	we still do the optimization in the isotropic space). Note that
	\[
		\max_{\step} \langle \step, \nabla \Cost(x)\rangle
		\quad \text{s.t.} \quad
		\|\step\|_\Sigma = \stepsize
	\]
	is equivalent to
	\[
		\max_{\step} \langle \step, \Sigma^{-1}\nabla \Cost(x)\rangle_\Sigma
		\quad \text{s.t.} \quad
		\|\step\|_\Sigma = \stepsize
	\]
	which is solved by
	\[
		\pm\stepsize\frac{\Sigma^{-1}\nabla\Cost(x)}{\|\Sigma^{-1}\nabla\Cost(x)\|_\Sigma}
	\]
	The remainder of the proof is exactly the same as in the original.

\subsection{Conservative RFD}\label{sec: conservative rfd}

In the first paragraph of Section~\ref{sec: rfd} we motivated the relation
between RFD and classical optimization with the observation, that gradient
descent is the minimizer of a regularized first order Taylor approximation
\[
	\tfrac1L \nabla\Cost(\param)
	= \argmin_\step T[\Cost(\param-\step)\mid \Cost(\param),\nabla\Cost(\param)] + \tfrac{L}2 \|\param\|^2.
\]
This regularized Taylor approximation is in fact an upper bound on our function
under the \(L\)-smoothness assumption \parencite{nesterovLecturesConvexOptimization2018},
i.e.
\[
	\Cost(\param-\step)
	\le T[\Cost(\param-\step)\mid \Cost(\param),\nabla\Cost(\param)]
	+ \tfrac{L}2 \|\step\|^2
\]
An improvement on of this upper bound compared to \(\Cost(x)\) therefore
guarantees an improvement of the loss. This guarantee was lost with the
conditional expectation (on purpose, as we wanted to consider the average case).
Losing this guarantee also makes convergence proofs more difficult as they
typically make use of this improvement.
In view of the confidence intervals of Figure~\ref{fig: visualize conditional
expectation}, it is natural to ask for a similar upper bound in the random
setting, where this can only be the top of an confidence
interval. This is provided in the following theorem
\begin{lemma}[An \(\gamma\)-upper bound]
	\label{lem: an 1-eps upper bound}
	We have
	\[
		\Pr\Bigl(
			\Cost(\param-\step)
			\le \E[\Cost(\param - \step)\mid \Cost(\param), \nabla\Cost(\param)]
			+ \rho_\gamma(\|\step\|^2)
		\Bigr) \ge \gamma 
	\]
	for
	\(
		\rho_\gamma(\stepsize^2)	
		:= \cdf^{-1}(\gamma)\sigma(\stepsize^2)
	\)
	with 
	\[
		\sigma^2(\stepsize^2):=
			\ikernel(0) - \frac{\ikernel\bigl(\frac{\stepsize^2}2\bigr)^2}{\ikernel(0)}
			- \frac{\ikernel'\bigl(\frac{\stepsize^2}2\bigr)^2}{-\ikernel'(0)}\stepsize^2
	\]
	where \(\cdf\) is the cumulative distribution function (cdf) of
	the standard normal distribution.
\end{lemma}

\begin{proof}
	Note that the conditional variance is with the usual argument about the
	covariance matrices (cf. the proof of Thoerem~\ref{thm: explicit rfd}) using
	Lemma~\ref{lem: cov of derivatives} and an application of
	Theorem~\ref{thm: conditional gaussian distribution} given by
	\[
		\sigma^2(\|\param\|^2)
		:=\Var[\Cost(\param-\step) \mid \Cost(\param), \nabla\Cost(\param)]
		= \ikernel(0) - \frac{\ikernel\bigl(\frac{\|\step\|^2}2\bigr)^2}{\ikernel(0)}
		- \frac{\ikernel'\bigl(\frac{\|\step\|^2}2\bigr)^2}{-\ikernel'(0)}\|\step\|^2.
	\]
	Since the conditional distribution is normal (by~Theorem~\ref{thm:
	conditional gaussian distribution}), we have
	\[
		\frac{
			\Cost(\param- \step) - \E[\Cost(\param - \step)\mid \Cost(\param), \nabla\Cost(\param)]
		}{
			\sigma(\|\param\|^2)
		}
		\sim \normal(0,1).
	\]
	But this implies the claim
	\begin{align*}
		&\Pr\Bigl(
			\Cost(\param-\step)
			\le \E[\Cost(\param - \step)\mid \Cost(\param), \nabla\Cost(\param)]
			+ \rho_\gamma(\|\step\|^2)
		\Bigr)
		\\
		&=\Pr\Bigl(
			\frac{
				\Cost(\param-\step)
				- \E[\Cost(\param - \step)\mid \Cost(\param), \nabla\Cost(\param)]
			}{\sigma(\|\param\|^2)}
			\le \cdf^{-1}(\gamma)
		\Bigr)
		\\
		&= \cdf(\cdf^{-1}(\gamma)) = \gamma.
	\end{align*}
	To avoid the Gaussian assumption, one could apply the Markov inequality instead, or
	another applicable concentration inequality.
\end{proof}

Using this upper bound, we obtain a natural conservative extension of RFD
\begin{extension}[\(\gamma\)-conservative RFD]
	Let \(\Cost\sim\normal(\mu, \ikernel)\) and \(\rho_\gamma(\stepsize^2) =
	\cdf^{-1}(\gamma) \sigma(\stepsize^2)\), where \(\sigma\) is the conditional
	standard deviation as defined in Lemma~\ref{lem: an 1-eps upper bound}.
	Then the conservative RFD step direction is given by
	\[
		\stepsize^* \frac{\nabla\Cost(\param)}{\|\nabla\Cost(\param)\|}
		= \argmin_\step \E[\Cost(\param-\step) \mid \Cost(\param), \nabla\Cost(\param)]
		+ \rho_\gamma(\|\step\|^2)
	\]
	and the \(\gamma\)-conservative RFD step size is given by
	\[
		\stepsize^* = 
		\argmin_{\stepsize}\frac{\ikernel\bigl(\frac{\stepsize^2}2\bigr)}{\ikernel(0)}
		(\Cost(\param)-\mu)
		-  \stepsize\frac{\ikernel'\bigl(\frac{\stepsize^2}2\bigr)}{\ikernel'(0)} \|\grad\Cost(\param)\|
		+ \rho_\gamma(\stepsize^2).
	\]
\end{extension}
\begin{proof}
	The proof is the same as in Theorem~\ref{thm: explicit rfd} with
	Lemma~\ref{lem: first stoch Taylor} replaced by Lemma~\ref{lem: an 1-eps
	upper bound}.
\end{proof}

Taking multiple steps should generally have an averaging effect, so we expect
faster convergence for almost risk neutral minimization of the conditional
expectation (i.e. \(\gamma\approx \frac12\)). Here \(\gamma\) is a natural
parameter to vary conservatism. In a software implementation it might be a
good idea to call this parameter `conservatism' and rescale it to be in
\([0,1]\) instead of \([\tfrac12, 1]\). But formulas look cleaner with \(\gamma\). 

In Bayesian optimization it is much more common to reverse this approach and
minimize a lower confidence bound (`conservatism' \(< 0\) or \(\gamma <
\tfrac12\)) in order to encourage exploration. But since
RFD is forgetful, this is not a good idea for RFD.

\begin{remark}[Conservative RFD coincides asymptotically with RFD in high dimension]
	\label{rem: high dimension}
	While conservative RFD might seem like a good approach to fix the instability
	of RFD under the isotropy assumption on some optimization problems, the
	variance generally vanishes in high dimension
	\parencite[see][]{benningGradientSpanAlgorithms2024} and conservative RFD
	coincides asymptotically with RFD. We therefore believe that the underlying
	issue is not an overly risk-affine algorithm but rather that distributional
	assumptions, in particular the stationarity assumption, are violated when
	instabilities occur (cf.~Section~\ref{sec: appendix: input invariance}).
\end{remark}

Nevertheless, conservative RFD might be a good approach for lower dimensional,
risk-sensitive applications.

\subsection{Beyond the Gaussian assumption}\label{sec: BlUE}

In this section we sketch how the extension beyond the Gaussian case using the
``best linear unbiased estimator'' BLUE
\parencite[e.g.][ch.~7]{johnsonAppliedMultivariateStatistical2007} works.

For this we recapitulate what a BLUE is. A \textbf{linear estimator}
\(\hat{Y}\) of \(Y\) using \(X_1,\dots,X_n\) is of the form 
\begin{equation*}
	\hat{Y}\in \linHull\{X_1,\dots,X_n\} + \real.
\end{equation*}
The set of \textbf{unbiased linear estimators} is defined as
\begin{align}\label{def: LUE}
	\LUE = \LUE{}[Y\mid X_1,\dots, X_n]
	&= \{ \hat{Y} \in \linHull\{X_1,\dots,X_n\} + \real : \E[\hat{Y}] = \E[Y]\}\\
	\nonumber
	&= \{ \hat{Y} + \E[Y] : \hat{Y} \in \linHull\{X_1-\E[X_1],\dots,X_n-\E[X_n]\}\}.
\end{align}
And the BLUE is the \textbf{best linear unbiased estimator}, i.e.
\begin{equation}\label{def: BLUE}
	\BLUE[Y\mid X_1,\dots, X_n] := \argmin_{\hat{Y}\in\LUE} \E[\|\hat{Y} - Y\|^2].
\end{equation}
Other risk functions to minimize are possible, but this is the usual one.

\begin{lemma}\label{lem: blue is cond. expectation}
	If \(X,Y_1,\dots, Y_n\) are multivariate normal distributed, then we
	have
	\begin{align*}
		\BLUE[Y\mid X_1,\dots,X_n]
		&= \E[Y\mid X_1,\dots, X_n]\\
		\bigg(&=
		\argmin_{\hat{Y}\in\{f(X_1,\dots, X_n) : f\text{ meas.}\}} \E[\|Y-\hat{Y}\|^2]
		\bigg).
	\end{align*}
\end{lemma}
\begin{proof}
	We observe that the conditional expectation of Gaussian
	random variables is linear (Theorem~\ref{thm: conditional gaussian
	distribution}). So as a linear function its \(L^2\) risk must be larger or
	equal to that of the BLUE. And as an \(L^2\) projection
	\parencite[Cor.~8.17]{klenkeProbabilityTheoryComprehensive2014} the
	conditional expectation was already optimal.
\end{proof}

If we now replace the conditional expectation with the BLUE, then all our theory
remains the same because the result in Theorem~\ref{thm: conditional gaussian distribution} 
remains the BLUE for general distributions
\parencite{johnsonAppliedMultivariateStatistical2007}.  Instead of minimizing
\[
	\min_\step\E[\Cost(\param-\step) \mid \Cost(\param), \nabla\Cost(\param)]	
\]
we can therefore always minimize
\[
	\min_\step \BLUE[\Cost(\param-\step) \mid \Cost(\param), \nabla\Cost(\param)]
\]
without the Gaussian assumption and all our results can be translated to this case.
The reader only needs to replace all mentions of Theorem~\ref{thm: conditional
gaussian distribution} with the BLUE equivalent and replace all ``idependence''
claims with ``uncorrelated''.

	\section{Input invariance}
\label{sec: appendix: input invariance}

In this section we generalize the notion of isotropy to non-stationary isotropy
and discuss why we believe this generalization is necessary. Recall that we
motivated isotropy as an invariant distribution with regard to isometric
transformations of the input.  In particular its distribution stays invariant
with regard to translations (also known as stationarity), which we do not
believe plausible for cost functions, because the cost at zero \(\Cost(0)\)
behaves fundamentally different from the cost of any other parameter vector.

In the following we will therefore generalize this notion to general input invariant
distributions. And we will discuss their applicability to machine learning after
we characterize the named categories.

\begin{definition}[Input Invariance]
    A random function \(\rf\) is \(\Phi\)-input invariant, if
	 \footnote{
		The input to a random function is somewhat ambiguous since it is a
		random variable, i.e. function from the probability space \(\Omega\)
		into function space, so its first input should be \(\omega\in
		\Omega\). Formally, the definition should therefore be: For
		all measurable sets of functions \(A\)
		\[
			\Pr_{\rf}(\phi_*^{-1}(A)) = \Pr_{\rf}(A) \quad \forall \phi \in \Phi
		\]
		where \(\phi_* : f\mapsto f\circ \phi\) denotes the pullback. But this
		is less helpful for an intuitive understanding.
	 }
	 \[
		\Pr_{\rf} = \Pr_{\rf \circ \phi} \qquad \forall \phi \in \Phi.
	 \]
	For certain sets of \(\Phi\) we give these \(\Phi\)-input invariant distributions
	names
    \begin{itemize}[topsep=0pt,parsep=0pt]
        \item If \(\Phi\) is the set of \emph{isometries}, we call \(\rf\)
        (stationary) \textbf{isotropic}.

        \item If \(\Phi\) is the set of \emph{translations}, we call \(\rf\)
        \textbf{stationary}.

        \item If \(\Phi\) is the set of \emph{linear isometries}, we call
        \(\rf\) \textbf{non-stationary isotropic}.
    \end{itemize}
    We further say a random function \(\rf\) is \(n\)-weakly \(\Phi\)-input invariant, if
    for all \(\phi\in \Phi\), all \(k\le n\) and all \(x_i\)
    \[
        \E[\rf(\phi(x_1))\cdot \dots \cdot \rf(\phi(x_k))] = \E[\rf(x_1)\cdot\dots\cdot\rf(x_k)].
    \]
    Since second moments fully determine Gaussian distributions, \(2\)-weakly input invariance
    is special, because it is equivalent to full input invariance in the
    Gaussian case. So an omitted \(n\) equals \(2\). ``Weakly isometry
    invariant'' naturally becomes ``weakly isotropic'', etc.
\end{definition}

While stationary and stationary isotropic random functions are well known
\citep[e.g.][]{rasmussenGaussianProcessesMachine2006,adlerRandomFieldsGeometry2007}, we
are not aware of research on non-stationary isotropy although we doubt the concept
is new. It turns out that the different notions of input isometry have simple
characterizations in terms of the covariance functions. We present these in
Theorem~\ref{thm: characterization of weak input invariances} of which the
stationary isotropic and stationary case are already well known.

\begin{restatable}[Characterization of Weak Input Invariances]{theorem}{characInputInv}
    \label{thm: characterization of weak input invariances}
    Let \(\rf:\real^\dims \to \real\) be a random function, then \(\rf\) is
    \begin{enumerate}[topsep=0pt,itemsep=0pt,partopsep=0pt]
        \item weakly stationary, if and only if there exists \(\mu\in\real\) and
        function \(\ikernel:\real^\dims \to \real\) such that for all \(x,y\)
        \[
            \mu_{\rf}(x) = \mu,
            \qquad \C_{\rf}(x,y) =
            \ikernel(x-y).
        \]

        \item weakly non-stationary isotropic, if and only if there exist
        functions \(\mu:\real_{\ge 0}\to\real\) and \(\kernel:D\to\real\) with
		  \(
			D = \{ \lambda\in \real_{\ge 0}^2 \times \real : |\lambda_3| \le 2\sqrt{\lambda_1\lambda_2}\}
			\subseteq\real^3.
			\)  
		  such that for all \(x,y\)
        \begin{align*}
            \mu_{\rf}(x) &= \mu\bigl(\tfrac{\|x\|^2}2\bigr)
            \\
            \C_\rf(x,y) &= \kernel\bigl(\tfrac{\|x\|^2}2, \tfrac{\|y\|^2}2, \langle x, y\rangle\bigr)
        \end{align*}
        
        \item\label{item: weak isotropy characterization} 
        weakly stationary isotropic, if and only if there exists \(\mu\in \real\) and a function
        \(\ikernel:\real_{\ge 0}\to \real\) such that for all \(x,y\)
        \[
            \mu_{\rf}(x) = \mu,
            \qquad
            \C_{\rf}(x,y) = \ikernel\bigl(\tfrac{\|x-y\|^2}2\bigr)
        \]
    \end{enumerate}
\end{restatable}
\begin{proof}
	The proof essentially follow as a corollary from a characterization
	of isometries (Proposition~\ref{prop: characterization of isometries}). For
	details see Subsec~\ref{sec: proof of theorem characterize input invariances}.
\end{proof}

Non-stationary isotropy is therefore a generalization of stationary isotropy.
It allows the zero parameter vector to have special meaning because the
distribution is only invariant to linear
isometries (i.e. rotations and reflections) which keep the zero in place.

It is important to highlight, that a geometric anisotropy (Section~\ref{sec:
geometric anisotropy}) retains stationarity, while breaking non-stationary
isotropy. A similar geometric generalization could also be applied to non-stationary
isotropy.

Another important observation is the fact, that non-stationary isotropy
coincides with stationary isotropy on the sphere. I.e. when \(\|x\|\) and
\(\|y\|\) are constant, the function
\[
	\kernel\bigl(\tfrac{\|x\|^2}2, \tfrac{\|y\|^2}2, \tfrac{\|x-y\|^2}2\bigr)
\]
only depends on \(\|x-y\|\) and the mean is also constant. In other words, we
have stationary isotropy on the sphere.

Isotropy might therefore `get by' as an assumption in machine learning,
as parameters are typically initialized on the sphere. This is because Glorot
intialization \parencite{glorotUnderstandingDifficultyTraining2010} samples
parameter entries independently, so their squared norm
\[
	\|\param\|^2 = \sum_{i=1}^\dims (\param^{(i)})^2
\]
is a sum of independent random variables which are normalized such that a law of large
numbers applies. Up to small variance their lengths are therefore all the same, and
are placed on a sphere at this radius.

If we leave this sphere, this equivalence stops being true. Weight normalization
\parencite{salimansWeightNormalizationSimple2016}, batch normalization
\parencite{ioffeBatchNormalizationAccelerating2015}, weight decay
\parencite[e.g.][]{goodfellowDeepLearning2016} or equivalently \(L^2\)
regularization, etc. might all contribute to keep this assumption intact.

But in the following section we will see, that even simple linear regressions
considered by researches investigating the average case behavior on quadratic
functions \parencite[e.g.][]{
	zhangWhichAlgorithmicChoices2019,
	pedregosaAccelerationSpectralDensity2020,
	lacotteOptimalRandomizedFirstOrder2020,
	deiftConjugateGradientAlgorithm2021,
	cunhaOnlyTailsMatter2022,
	paquetteHaltingTimePredictable2022,
	paquetteUniversalityConjugateGradient2022}, require non-stationary
isotropy. Moreover the covariance kernels suggested by investigations into
random neuronal networks
\parencite[e.g.][]{williamsComputationInfiniteNeural1998,choKernelMethodsDeep2009}
are also non-stationary isotropic but not stationary isotropic.

\subsection{Random linear regression}\label{sec: random regression}

In this section, we determine the distribution of the cost function induced by
a simple linear regression. For this we define the mean squared sample loss
\[
    \loss_i(\param) = (Y - f_\param(X))^2,
\]
where the random data \(X\) is mapped by the true relationship \(\rf\) to labels
\(Y = \rf(X)\) and
\[
    f_\param(x) = \langle x, \param\rangle
\]
is a linear model. If the true relationship \(\rf\) is also a random linear function
\(\rf(x) = \langle \signal, x\rangle\) with random signal
\(\signal\sim\normal(0,\identity)\) independent of input
\(X\sim\normal(0,\identity)\), then the cost
function is given by
\begin{align*}
    \Cost(\param)
    &= \E[\loss_i(\param)\mid \rf]
    = \E[\langle \signal - \param, X\rangle^2\mid \signal]
    \\
    &= (\signal - \param)^T\E[XX^T](\signal -\param)
    \\
    &=  \|\signal - \param\|^2
\end{align*}

\begin{lemma}
    The expectation and covariance of \(\Cost\) are given by
    \begin{align*}
        \E[\Cost(\param)] &= \text{const.} + \|\param\|^2
        \\
        \Cov(\Cost(\param), \Cost(\tilde{\param}))    
        &= \text{const.} + 4 \langle\param, \tilde{\param}\rangle
    \end{align*}
    In particular, the cost \(\Cost\) is non-stationary isotropic, but not
    stationary isotropic.
\end{lemma}

\begin{proof}
    Its expectation is given by
    \begin{align*}
        \E[\Cost(\param)]
        &= \E[\|\signal - \param\|^2]
        = \E[\|\signal\|^2] - 2\langle \underbrace{\E[\signal]}_{=0}, \param\rangle + \|\param\|^2
        \\
        &= \text{const.} + \|\param\|^2
    \end{align*}
    In particular it is not constant, but dependent on \(\|\param\|^2\), which means that we
    do not have stationary isotropy. But there is still hope for non-stationary isotropy, and this is 
    essentially true as can be seen by calculating
    \begin{align*}
        \Cov(\Cost(\param), \Cost(\tilde{\param}))    
        &= \E\Bigl[
            (\Cost(\param) - \E[\Cost(\param)])(\Cost(\tilde{\param}) - \E[\Cost(\tilde{\param})])
        \Bigr]
        \\
        &=\E\Bigl[
            (\|\signal\|^2 - \E\|\signal\|^2 - 2\langle \signal, \param\rangle)
            (\|\signal\|^2 - \E\|\signal\|^2 - 2\langle \signal, \tilde{\param}\rangle)
        \Bigr] 
        \\
        &= \begin{aligned}[t]
            &\Var(\|\signal\|^2 - \E\|\signal\|^2)
            \\
            &- 2\E[(\|\signal\|^2 - \E\|\signal\|^2)\langle \signal, \param\rangle]
            \\
            &- 2\E[(\|\signal\|^2 - \E\|\signal\|^2)\langle \signal, \tilde{\param}\rangle]
            \\
            &+ 4 \param^T\E[\signal\signal^T]\tilde{\param}
        \end{aligned}
        \\
        &= \Var(\|\signal\|^2 - \E\|\signal\|^2) + 4 \langle \param, \tilde{\param}\rangle
        \\
        &= \text{const.} + 4 \langle\param, \tilde{\param}\rangle
    \end{align*}
    because the terms in the middle are zero, e.g.
    \[
        \E[(\|\signal\|^2 - \E\|\signal\|^2)\langle \signal, \param\rangle]
        = \langle \underbrace{\E[\|\signal\|^2\signal]}_{=0}, \param\rangle
        - \E\|\signal\|^2\langle \underbrace{\E[\signal]}_{=0}, \param\rangle
    \]
    where the entries of \(\E[\|\signal\|^2\signal]\) are zero, because of
    independence and first moments being zero and third moments being zero.
\end{proof}

\subsection{Proof of Theorem~\ref{thm: characterization of weak input invariances}}
\label{sec: proof of theorem characterize input invariances}

\begin{prop}[Characterizing isometries]\label{prop: characterization of isometries}
	Let \(\cX\) be a vectorspace and \(x_i,y_i \in \cX\) for \(i=1,\dots, n\), then
	the following pairs of statements are equivalent
	\begin{enumerate}
		 \item \label{item: translation characterization}
		 \begin{enumerate}
			  \item\label{itm: translation charac. a}
			  \(x_i - x_j = y_i - y_j\) for all \(i,j\)

			  \item\label{itm: translation charac. b}
			  there exists a \textbf{translation} \(\phi\) with \(\phi(x_i) = y_i\) for all \(i\).
		 \end{enumerate}
	\end{enumerate}
	In the remainder we further assume \(\cX\) to be a Hilbertspace, 
	\begin{enumerate}[resume]
		 \item \label{item: linear isometry characterization}
		 \begin{enumerate}
			  \item \label{itm: lin. isometry charac. a}
			  \(\|x_i\|= \|y_i\|\) and \(\|x_i-x_j\| = \|y_i-y_j\|\) for all \(i,j\)
			  
			  \item \label{itm: lin. isometry charac. b}
			  there exists a \textbf{linear isometry} \(\phi\) with \(\phi(x_i) = y_i\)
			  for all \(i\).
		 \end{enumerate}

		 \item \label{item: affine isometry characterization}
		 \begin{enumerate}
			  \item\label{itm: affine isometry charac. a}
			  \(\|x_i - x_j\| = \|y_i - y_j\|\) for all \(i,j\)

			  \item\label{itm: affine isometry charac. b}
			  there exists an \textbf{(affine) isometry} \(\phi\) with \(\phi(x_i) = y_i\) for all \(i\).
		 \end{enumerate}
	\end{enumerate}
\end{prop}

\begin{proof}
	\begin{description}
		 \item[\eqref{itm: translation charac. a} \(\Rightarrow\) \eqref{itm: translation charac. b}:] 
		 we define
		 \[
			  \phi(x) := x + (y_0-x_0),
		 \]
		 which implies
		 \[
			  \phi(x_i) = x_i - x_0 + y_0 = (y_i-y_0) + y_0 = y_i.
		 \]
		 
		 \item[\eqref{itm: translation charac. b} \(\Rightarrow\) \eqref{itm: translation charac. a}:] 
		 Let \(\phi(x) = x + c\) for some \(c\). Then we immediately have
		 \[
			  y_i - y_j = \phi(x_i) - \phi(x_j) = x_i + c - (x_j + c) = x_i - x_j.
		 \]

		 \item[\eqref{itm: lin. isometry charac. a} \(\Rightarrow\) \eqref{itm: lin. isometry charac. b}:] 
		 By the polarization formula, for all \(i,j\)
		 \[
			  \langle x_i, x_j\rangle = \frac{\|x_i\|^2 + \|y_i\|^2 - \|x_i - x_j\|^2}{2}
			  = \langle y_i, y_j\rangle.
		 \]
		 We apply the Gram-Schmidt orthonormalization procedure to both \(x_i\) and
		 \(y_i\) such that
		 \[
			  U_{k_n} = \Span(u_1, \dots, u_{k_n}) = \Span(x_1,\dots, x_n)
		 \]
		 for orthonormal \(u_i\) where we skip \(x_m\) if it is already
		 in \(U_{k_{m-1}}\) (resulting in \(k_m = k_{m-1}\)), and similarly 
		 \[
			  V_{k_n} = \Span(v_1, \dots, v_{k_n}) = \Span(y_1,\dots, y_n).
		 \]
		 Since this procedure only uses scalar products, we inductively get
		 \[
			  \langle x_k, u_j\rangle = \langle y_k, v_j\rangle
			  \quad \forall k, j
		 \]
		 We now extend \(u_i\) and \(v_i\) to orthonormal basis of \(\cX\) and
		 define the linear mapping by its behavior on the basis
		 elements \(\phi: u_i \mapsto v_i\).
		 Mapping an orthonormal basis to an orthonormal basis is an isometry
		 and we have
		 \begin{align*}
			  \phi(x_k)
			  &= \phi\Bigl(
					\sum_{j=1}^k \langle x_k, u_j\rangle u_j
			  \Bigr)
			  \\
			  &= \sum_{j=1}^k \langle x_k, u_j\rangle \phi(u_j)
			  = \sum_{j=1}^k \langle y_k, v_j\rangle v_j 
			  \\
			  &= y_k.
		 \end{align*}

		 \item[\eqref{itm: lin. isometry charac. b} \(\Rightarrow\) \eqref{itm: lin. isometry charac. a}:] 
		 Isometries preserve distances by definition. This implies \(\|x_i - x_j\|
		 = \|y_i - y_j\|\). And linear functions map
		 \(0\) to \(0\), so we have 
		 \[
			  \|x_i\| = \|x_i - 0\| = \|\phi(x_i) - \phi(0)\| = \|y_i\|.
		 \]

		 \item[\eqref{itm: affine isometry charac. a} \(\Rightarrow\) \eqref{itm: affine isometry charac. b}:] 
		 We define
		 \[
			  \tilde{x}_i = x_i - x_0
		 \]
		 and similarly for \(y\). In particular, \(\tilde{x}_0 = \tilde{y}_0=0\). Since
		 \(\tilde{x}_i\) and \(\tilde{y}_i\) satisfy the requirements of
		 \ref{item: linear isometry characterization}, there exists a linear
		 isometry \(\tilde{\phi}\) with \(\tilde{\phi}(\tilde{x}_i) = \tilde{y}_i\).
		 Then the isometry
		 \[
			  \phi : x\mapsto \tilde{\phi}(x-x_0) + y_0
		 \]
		 does the job.

		 \item[\eqref{itm: affine isometry charac. b} \(\Rightarrow\) \eqref{itm: affine isometry charac. a}:] 
		 This is precisely the distance preserving property of Isometries.
		 \qedhere
	\end{description}
\end{proof}

\characInputInv*

\begin{proof}
	Starting from the mean and covariance function it is easy to check
	\(2\)-weak non-stationary isotropy. So we only need to check the other direction.

	The proof is essentially an application of Prop.~\ref{prop: characterization of isometries}.
	For brevity (and since the other two results are well known), we will only prove the
	weakly non-stationary isotropic case (the other two cases can be proven with minor adjustments to the proof).

	Without loss of generality, we will find the slightly different representation
	\[
		\E[\rf_\dims(x)] = \tilde{\mu}(\|x\|)
		\quad\text{and}\quad
		\C_{\rf_\dims}(x,y)	= \tilde{\kernel}(\|x\|, \|y\|, \langle x,y\rangle),
	\]
	where the domain of \(\tilde{\kernel}\) is given by \(\tilde{D} = \{
	\lambda\in \real_{\ge 0}^2 \times \real : |\lambda_3| \le
	\lambda_1\lambda_2\}\).	The representation of the theorem is then
	equivalent by a change to
	\[
		\mu(\lambda) := \tilde{\mu}\bigl(\tfrac{\lambda^2}2\bigr)
		\quad\text{and}\quad 
		\kernel(\lambda_1,\lambda_2,\lambda_2) := \tilde{\kernel}\bigl(\tfrac{\lambda_1^2}2, \tfrac{\lambda_2^2}2, \lambda_3\bigr).
	\]
	
	First we want to find \(\mu\). Let \(v\) be some vector (w.l.o.g. \(\|v\|=1\)). Then we define
	\[
		\mu(r) := \E[\rf_\dims(rv)]
	\]
	Now we need to show that this definition of \(\mu\) is an appropriate mean
	function. For this choose any \(x\in\cX\). Then for \(r=\|x\|\) there exists by
	Prop.~\ref{prop: characterization of isometries}~(\ref{item: linear
	isometry characterization}.) a non-stationary isometry \(\phi\) such that \(\phi(x) =
	rv\) (we use \(n=1\)). With \(1\)-weak non-stationary isotropy of \(\rf_\dims\) this implies
	\[
		\E[\rf_\dims(x)] = \E[\rf_\dims(rv)] = \mu(r) = \mu(\|x\|).
	\]
	Next we need to define \(\kernel(r_x, r_y, r_{xy})\). For this, choose two orthonormal vectors
	\(v,w\). For every \(r=(r_x,r_y, r_{xy})\in \tilde{D}\) we define
	\begin{align*}
		x^*(r) &= r_x v
		\\
		y^*(r) &= \frac{r_{xy}}{r_x} v + \sqrt{r_y^2 - \tfrac{r_{xy}^2}{r_x^2}} w.
	\end{align*}
	Where \(r\in \tilde{D}\) ensures \(|r_{xy}|\le r_x r_y\) and thus  \(r_y^2 -
	\frac{r_{xy}^2}{r_x^2}\ge 0\).
	Then we have
	\begin{equation}
		\label{eq: correct norm and angle}
		\|x^*(r)\| = r_x,
		\quad
		\|y^*(r)\| = r_y,
		\quad\text{and}\quad
		\langle x^*(r), x^*(y) \rangle = r_{xy},
	\end{equation}
	and define	
	\[
		\kernel(r_x, r_y, r_{xy}) := \C_{\rf_\dims}(x^*(r), y^*(r)).
	\]
	Again, we need to show that this kernel does the job. For this, choose any
	\(x,y\in \cX\). For
	\[
		r :=(\|x\|, \|y\|, \langle x,y\rangle),
	\]
	which is in \(\tilde{D}\) by the Cauchy-Schwarz inequality, the induced
	\(x^*(r)\) and \(y^*(r)\) satisfy by \eqref{eq: correct norm and
	angle}
	\[
		\|x^*(r)\| = \|x\|,
		\quad
		\|y^*(r)\| = \|y\|	
		\quad\text{and}\quad
		\|x^*(r) - y^*(r)\| = \|x-y\|.
	\]
	By Prop.~\ref{prop: characterization of isometries}~(\ref{item: linear
	isometry
	characterization}.) there therefore exists an isometry \(\phi\) such that
	\(\phi(x) = x^*(r)\) and \(\phi(y) = y^*(r)\). By \(2\)-weak input isotropy of
	\(\rf_\dims\) we conclude
	\[
		\C_{\rf_\dims}(x,y) \overset{\text{isotrop.}}= \C_{\rf_\dims}(x^*(r),y^*(r)) \overset{\text{def.}}= \kernel\bigl(\|x\|, \|y\|, \langle x,y\rangle\bigr).
		\qedhere
	\]
\end{proof}

	\section{Technical}


\subsection{Conditional Gaussian distribution}

For the following well known result we found a tidy proof giving insight into
the reason it is true, so we wrote it down for your convenience but do
not even expect this particular proof to be new.

\begin{theorem}[Conditional Gaussian distribution]
    \label{thm: conditional gaussian distribution}
    Let \(X\sim\normal(\mu,\Sigma)\) be a multivariate Gaussian vector where
    the covariance matrix is a block matrix of the form
    \[
        \mu = \begin{bmatrix}
            \mu_1\\ \mu_2
        \end{bmatrix}
        \quad \text{and} \quad
        \Sigma = \begin{bmatrix}
            \Sigma_{11} & \Sigma_{12}
            \\
            \Sigma_{21} & \Sigma_{22}
        \end{bmatrix},
    \]
    then assuming \(\Sigma_{11}\) is invertible, the conditional distribution of
    \(X_2\) given \(X_1\) is
	 \[
        X_2\mid X_1 \sim \normal(\mu_{2\mid 1}, \Sigma_{2\mid 1}),
    \]
    with conditional mean and variance
    \begin{align*}
        \mu_{2\mid 1} &:= \mu_2 + \Sigma_{21}\Sigma_{11}^{-1}(X_1-\mu_1)
        \\
        \Sigma_{2\mid 1} &:= \Sigma_{22} -  \Sigma_{21}\Sigma_{11}^{-1}\Sigma_{12}.
    \end{align*}
\end{theorem}
\begin{proof}
    Let \(\bar{X} := X-\mu\) be the centered version of \(X\).  There exists
    some lower triangular matrix \(L\) (even if \(\Sigma\) is only positive
    semidefinite only not uniquely) such that \(\Sigma = LL^T\) (i.e.\ the Cholesky
    Decomposition). We can then write without loss of generality
    \[
        X - \mu =:
        \begin{bmatrix}
            \bar{X}_1
            \\
            \bar{X}_2
        \end{bmatrix}
        = \begin{bmatrix}
            L_{11} & 0
            \\
            L_{21} & L_{22}
        \end{bmatrix}
        \begin{bmatrix}
            Y_1 \\ Y_2
        \end{bmatrix}
        = LY
    \]
    with independent standard normal \(Y_i\), i.e. \(Y\sim\normal(0, \identity)\).
    Since \(\Sigma_{11}\) is invertible, so is \(L_{11}\) and therefore the map
    from \(Y_1\) to \(X_1\). Conditioning on \(X_1\) is therefore equivalent to
    conditioning on \(Y_1\). But we have
    \[
        X_2 = \mu_2 + \bar{X}_2
        = \underbrace{\mu_2 + L_{21}Y_1}_{\text{conditional expectation}} + \underbrace{L_{22} Y_2}_{\text{conditional distribution}}
    \]
    So it follows that
    \[
        X_2\mid X_1 \sim  \normal(\mu_{2\mid 1}, \Sigma_{2\mid 1})
    \]
    with
    \begin{align*}
        \mu_{2\mid 1} &:=  \mu_2 + L_{21}Y_1
        \\
        \Sigma_{2\mid 1} &:= L_{22}L_{22}^T.
    \end{align*}
    What is left to do, is find a representation for the \(L_{ij}\)
    using the block matrices of \(\Sigma\). For this note
    \[
        \Sigma = LL^T = \begin{bmatrix}
            L_{11}L_{11}^T & L_{11}L_{21}^T
            \\
            L_{21}L_{11}^T & L_{22}L_{22}^T + L_{21}L_{21}^T
        \end{bmatrix}
    \]
    This implies
    \[
        L_{21}Y_1 = (L_{21}L_{11}^T L_{11}^{-T})(L_{11}^{-1}\bar{X}_1)
        = \Sigma_{21}\Sigma_{11}^{-1}(X_1-\mu_1)
    \]
    so we have the desired conditional expectation, and finally
    \begin{align*}
        L_{22}L_{22}^T
        &= \Sigma_{22} - L_{21}L_{21}^T
        \\
        &= \Sigma_{22} - \underbrace{L_{21}(L_{11}^T}_{=\Sigma_{21}}\underbrace{L_{11}^{-T}) (L_{11}^{-1}}_{=\Sigma_{11}^{-1}}\underbrace{L_{11})L_{21}^T}_{=\Sigma_{12}}.
        \qedhere
    \end{align*}
\end{proof}

\subsection{Covariance of derivatives}

By Swapping integration and differentiation we have for a centered random function \(\rf\)
\begin{align*}
	\Cov(\partial_{x_i} \rf(x), \rf(y))
	&= \E[\partial_{x_i} \rf(x) \rf(y)]
	= \partial_{x_i} \E[\rf(x)\rf(y)]
	\\
	&= \partial_{x_i} \C_{\rf}(x,y)
\end{align*}
So the covariance of a derivative of \(\rf\) with \(\rf\) is equal to a partial
derivative of the covariance function \citep[more details
in][]{adlerRandomFieldsGeometry2007}. Similarly other covariances can be calculated, e.g.
\[
	\Cov(\partial_{x_i}\rf(x), \partial_{y_i}\rf(y))
	= \partial_{x_i}\partial_{y_i} \C_{\rf}(x,y).
\]
For this reason the derivatives of the covariance function are interesting as they represent
the covariance of derivatives.

Applying this observation to isotropic covariance functions
\[
    \Cov(\rf(x),\rf(y)) = \ikernel\bigl(\tfrac{\|x-y\|^2}2\bigr)
\]
we obtain.
\begin{lemma}[Covariance of derivatives]
	\label{lem: cov of derivatives}
	Let \(\rf\sim\normal(\mu, \ikernel)\) and \(\step = x-y\), then
	\begin{equation*}
	\def\arraystretch{1.5}
	\begin{tabular}{c | c c}
		\(\Cov\) & \(\rf(y)\) & \(\partial_j\rf(y)\) \\	
		\hline
		\(\rf(x)\)
		&  \(\ikernel(\frac{\|\step\|^2}{2})\)
		& \(-\ikernel'(\frac{\|\step\|^2}2)\langle \step, e_j\rangle\)
		\\
		\(\partial_i\rf(x)\)
		&  \(\ikernel'(\frac{\|\step\|^2}{2})\langle \step, e_i\rangle\)
		& \(
			- \Bigl[
					\ikernel''(\frac{\|\step\|^2}2)\langle \step, e_j\rangle\langle \step, e_i\rangle
					+ \ikernel'(\frac{\|\step\|^2}2)\langle e_j, e_i\rangle
			\Bigr]
		\)
	\end{tabular}
	\end{equation*}
\end{lemma}

\subsection{Constrained linear optimization}

Let \(U\) be a vectorspace. We define the projection of a vector \(w\) onto \(U\) by
\[
	\proj{U}(w) := \argmin_{v\in U} \|v-w\|^2
\]
\begin{lemma}[Constrained maximiziation of scalar products]
	\label{lem: constrained maximiziation of scalar products}
	For a linear subspace \(U\subseteq \real^\dims\), we have	
	\begin{align}
		\max_{\substack{v\in U \\ \|v\|=\lambda}}\langle v, w\rangle 
		&= \lambda \|\proj{U}(w)\|
		\\
		\argmax_{\substack{v\in U\\ \|v\|=\lambda}} \langle v, w\rangle
		&= \lambda \frac{\proj{U}(w)}{\|\proj{U}(w)\|}
	\end{align}
\end{lemma}
Before we get to the proof let us note that this immediately results in the following
corollary about minimization.
\begin{corollary}[Constrained minimization of scalar products]
	\label{cor: constrained minimization of scalar products}
	\begin{align}
		\min_{\substack{v\in U \\ \|v\|=\lambda}}\langle v, w\rangle 
		&= -\lambda \|\proj{U}(w)\|
		\\
		\argmin_{\substack{v\in U\\ \|v\|=\lambda}} \langle v, w\rangle
		&= -\lambda \frac{\proj{U}(w)}{\|\proj{U}(w)\|}
	\end{align}
\end{corollary}
\begin{proof}[Proof of Corollary~\ref{cor: constrained minimization of scalar products}]
	The trick is to move one `\(-\)' outside from \(w=-(-w)\)
	\[
		\min_{\substack{v\in U \\ \|v\|=\lambda}}\langle v, w\rangle 
		= - \max_{\substack{v\in U \\ \|v\|=\lambda}} \langle v, -w \rangle
		= -\lambda \|\proj{U}(w)\|
	\]
	where we have used in the last equation that the projection is linear (we can move
	the minus sign out) and the norm removes the inner minus sign. The \(\argmin\) argument
	is similar.
\end{proof}

\begin{proof}[Proof of Lemma~\ref{lem: constrained maximiziation of scalar products}]
	\begin{enumerate}[label={\bf{Step \arabic*}:},wide]
		\item  We claim that
		\[
			v^* = \lambda \frac{\proj{U}(w)}{\|\proj{U}(w)\|}
		\]
		results in the value \(\langle v^*, w\rangle = \lambda \|\proj{U}(w)\|\).

		For this we consider
		\begin{align}
			\label{eq: projection def}
			\proj{U}(w)
			&= \argmin_{v\in U}\underbrace{\|v-w\|^2}_{= \|v\|^2 - 2\langle v,w\rangle + \|w\|^2}
			\\
			\nonumber
			&= \argmin_{v\in U} \underbrace{\|v\|^2 - 2\langle v, w\rangle}_{=: f(v)}
		\end{align}
		we know that \(t\mapsto f(t \langle w\rangle_U)\) is minimized at \(t=1\)
		by the definition of \(\langle w\rangle_U\). The first order condition implies
		\[
			0\overset!= \frac{d}{dt}	
			= 2t \|\proj{U}(w)\|^2 - 2\langle \proj{U}(w), w\rangle
		\]
		and thus
		\[
			1 = t^*  = \frac{\langle \proj{U}(w), w\rangle}{\|\proj{U}(w)\|^2}
		\]
		Multiplying both sides by \(\lambda\|\proj{U}(w)\|\) finishes this step
		\begin{equation}
			\label{eq: maximum can be achieved}
			\lambda \|\proj{U}(w)\|
			= \Bigl\langle
				\underbrace{\lambda \frac{\proj{U}(w)}{\|\proj{U}(w)\|}}_{=v^*}, w
			\Bigr\rangle.
		\end{equation}

		\item By \eqref{eq: maximum can be achieved}, we know that we can achieve
		the value we claim to be the maximum (and know the location \(v^*\) to do
		so). So if we prove that we can not exceed this value, then it is a
		maximum and \(v^*\) is the \(\argmax\). This would finish the proof. What
		remains to be shown is therefore
		\[
			\langle v, w\rangle \le \lambda \|\proj{U}(w)\| \qquad \forall v\in U: \|v\| = \lambda.
		\]
		Let \(v\in U\) with \(\|v\|=\lambda\). Then for any \(\mu\in \real\) we can plug
		\(\mu v\) into \(f\) from \eqref{eq: projection def} to get
		\begin{align*}
			\mu^2\lambda^2 - 2\mu \langle v, w\rangle
			&= f(\mu v)
			\\
			\overset{\eqref{eq: projection def}}&\ge f(\proj{U}w)
			= \|\proj{U}(w)\|^2 - 2\langle \proj{U}w , w\rangle
			\\
			&= - \langle \proj{U}w, w\rangle
		\end{align*}
		where the last equation follows from \eqref{eq: maximum can be achieved} with
		\(\lambda = \|\proj{U}w\|\). Reordering we get for all \(\mu\)
		\[
			\langle \proj{U}w, w\rangle + \mu^2\lambda^2 \ge 2\mu\langle v,w\rangle
		\]
		We now select \(\mu = \frac{\|\proj{U}w\|}{\lambda}>0\) and divide both
		sides by \(\mu\) to get
		\[
			2\langle v, w\rangle
			\le \underbrace{
				\Bigl\langle \underbrace{\frac{\proj{U}(w)}{\mu}}_{=v^*}, w\Bigr\rangle
			}_{
				 \overset{\eqref{eq: maximum can be achieved}}= \lambda\|\proj{U}w\|
			} + \lambda\|\proj{U}(w)\|
			= 2 \lambda\|\proj{U}w\|
		\]
		Dividing both sides by \(2\) yields the claim.
		\qedhere
	\end{enumerate}
\end{proof}

\end{document}